\documentclass[11pt]{article}

\title{\begin{center}
Enumeration of unicuspidal curves of any degree and genus on toric surfaces
\end{center}}

\author{Yaniv Ganor \and Eugenii Shustin\thanks{School of Mathematical Sciences, Tel Aviv University,
Ramat Aviv, 6997801 Tel Aviv, Israel. E-mail: ganory@gmail.com, shustin@tauex.tau.ac.il}}

\date{}

\usepackage{fancyhdr}
\usepackage[all]{xy}
\usepackage{graphicx}
\usepackage{amsmath}
\usepackage{amsthm}
\usepackage{amssymb}

\usepackage{latexsym,amsfonts,amscd,yhmath,epsfig,epic}

\usepackage{url}

\newcommand{\proofend}{\hfill$\Box$\bigskip}
\newcommand{\oa}{\overline{a}}
\newcommand{\Z}{\mathbb{Z}}\newcommand{\ob}{\overline{b}}
\newcommand{\Q}{\mathbb{Q}}\newcommand{\bx}{\boldsymbol{x}}
\newcommand{\R}{\mathbb{R}}\newcommand{\eps}{\varepsilon}
\newcommand{\C}{\mathbb{C}}
\newcommand{\K}{\mathbb{K}}
\newcommand{\PP}{\mathbb{P}}
\newcommand{\ux}{\underline{x}}
\newcommand{\uz}{\underline{z}}
\newcommand{\bp}{\boldsymbol{p}}\newcommand{\bw}{\boldsymbol{w}}
\newcommand{\ini}{\operatorname{ini}}
\newcommand{\Span}{\operatorname{Span}}\newcommand{\Def}{\operatorname{Def}}
\newcommand{\Ini}{\operatorname{Ini}}\newcommand{\Sing}{\operatorname{Sing}}
\newcommand{\const}{\operatorname{const}}\newcommand{\OG}{\operatorname{OG}}
\newcommand{\mt}{\operatorname{mult}}
\newcommand{\ord}{\operatorname{ord}}
\newcommand{\pr}{\operatorname{pr}}\newcommand{\Tr}{\operatorname{Tr}}
\newcommand{\bn}{\boldsymbol{n}}\newcommand{\conj}{\operatorname{Conj}}

\newcommand{\Trop}{\mathcal{T}} %
\DeclareMathOperator{\Tor}{Tor} %
\DeclareMathOperator{\Int}{Int} %
\DeclareMathOperator{\val}{val} %
\DeclareMathOperator{\conv}{conv} %

\newtheorem{theorem}{Theorem}[section]

\newtheorem{definition}[theorem]{Definition}

\newtheorem{lemma}[theorem]{Lemma}

\newtheorem{remark}[theorem]{Remark}

\usepackage{xcolor}

\begin{document}

\maketitle

\begin{abstract}
We enumerate complex curves on toric surfaces of any given degree and genus, having a single cusp and nodes as their singularities, and matching appropriately many point constraints. The solution is obtained via tropical enumerative geometry.
The same technique applies to enumeration of real plane cuspidal curves: We show that, for any fixed $r\ge1$ and $d\ge2r+3$, there exists a generic real $2r$-dimensional linear family of plane curves of degree $d$
in which the number of real $r$-cuspidal curves is asymptotically comparable
with the total number of complex $r$-cuspidal curves in the family, as $d\to\infty$.
\end{abstract}


\section*{Introduction}

Tropical geometry has been successfully applied to various enumerative problems starting with the foundational work \cite{M}, where Mikhalkin computed closed and open Gromov-Witten invariants of toric surfaces, or, more precisely, counted \emph{nodal} curves of any degree and genus on arbitrary toric surfaces both over the complex and the real field.

The \emph{goal} of this paper is to enumerate curves of any degree and genus that lie on a given toric surface, have nodes and a single \emph{cusp} as their singularities, and pass through an appropriate number of fixed points. We attack this problem with the techniques of tropical enumerative geometry \footnote{Particular cases when either a cusp is the only singularity, or $g=0$ were addressed via tropical geometry in \cite{Shustin2005} and \cite{YG}, though the treatment was incomplete in both sources.}. The case of curves with a cusp appears to be rather more complicated than that of purely nodal curves. Here we state our results in the following theorems:
\begin{itemize}\item Theorem \ref{tA} - the correspondence between cuspidal algebraic and cuspidal tropical curves (a fully detailed formulation appears in Theorems \ref{tuc1}, \ref{t2}, and Lemma \ref{O5}).
\item Theorem \ref{tB} - enumeration of cuspidal tropical curves (a fully detailed formulation appears in Theorem \ref{t4}).
 \item Theorem \ref{tC} - enumeration of real cuspidal curves (see Theorem \ref{tex}).
\end{itemize}

We quickly recall the idea of the tropical approach in the spirit of \cite{M} and \cite{Shustin2005} (details can be found in Section \ref{sec-pr} and references therein). Over the field of complex, locally convergent Puiseux series, an algebraic curve in a toric surface can be viewed as an equisingular family of curves over the punctured disc embedded into the trivial family of toric surfaces. This pair of families admits a flat extension to the disc's center so that both, the central surface and the central curve split into several components. This splitting can be encoded by an appropriate plane tropical curve, which is a metric graph with marked points properly mapped to the plane $\R^2$. Two issues have to be resolved:
\begin{itemize}\item Classification of the corresponding tropical curves,
\item and computation of the multiplicity of each tropical curve, i.e., the number of counted algebraic curves that tropicalize to the given tropical curve.
    \end{itemize}
In the purely nodal case \cite{M}, the tropical curves appear to be trivalent graphs, and the multiplicity of such a tropical curve is the product of Mikhalkin's multiplicities of all trivalent vertices. A similar statement holds in the cuspidal case, though the geometry of cuspidal tropical curves appears to be more complicated.

\begin{theorem}\label{tA}
Let $\Delta\subset\R^2$ be a nondegenerate convex lattice polygon, $X=\Tor(\Delta)$ the complex toric surface and ${\mathcal L}_\Delta$ the ample line bundle associated with $\Delta$. Given a nonnegative integer $g<|\Int(\Delta)\cap\Z^2|$, the family $V_{\Delta,g}(A_2)$ of curves of genus $g$ in the linear system $|{\mathcal L}_\Delta|$ having nodes and a single cusp is of pure dimension $N=-c_1({\mathcal L}_\Delta)K_X+g-2$. Then
\begin{enumerate}\item[(i)] one has
\begin{equation}\deg V_{\Delta,g}(A_2)=\sum_{T\in\Trop_{\Delta,g}(A_2,\bx)}\mt_c(T),\label{eA}\end{equation} where $\Trop_{\Delta,g}(A_2,\bx)$ is the set of plane cuspidal tropical curves, whose tropical degree induces $\Delta$, tropical genus $g$, that pass through a generic configuration $\bx$ of $N$ generic points in $\R^2$.
The set $\Trop_{\Delta,g}(A_2,\bx)$ is finite, and each tropical curve $T\in\Trop_{\Delta,g}(A_2,\bx)$ contains a tropical cuspidal fragment, while the rest is trivalent. The cuspidal multiplicity $\mt_c(T)$ is the product of the weight of the cuspidal fragment and of the Mikhalkin's weights of the trivalent vertices outside the cuspidal fragment.
\item[(ii)] If $\Delta$ is h-transverse, i.e., is a union of lattice triangles and trapezoids of horizontal width $1$
\footnote{An equivalent definition of $h$-transverse polygons requires that the primitive normal integral vectors to the sides of $\Delta$ must be of the form $(m,k)$ with $|k|\le1$.}, then there exists a configuration $\bx\subset\R^2$ such that each tropical curve $T\in\Trop_{\Delta,g}$ has one of the following three cuspidal fragments: either a rational four-valent vertex, or a rational flat trivalent vertex, or an elliptic edge. The multiplicities of each of these cuspidal fragments can be computed via explicit closed formulas.
\end{enumerate}
\end{theorem}

Note that the toric surfaces associated with h-transversal polygons include the plane, the quadric, all Hirzebruch surfaces and many other smooth and singular surfaces. Furthermore, the second part of Theorem \ref{tA} is tightly linked with the enumeration of tropical curves in the set $\Trop_{\Delta,g}(A_2,\bx)$
via a certain modification of Mikhalkin's lattice path algorithm \cite[Section 7.2]{M}.

\begin{theorem}\label{tB}
Let $\Delta$, $g$, $N$, and $\bx$ be as in Theorem \ref{tA}(ii). Fix $0<\eps\ll1$ such that $\lambda:\Delta\cap\Z^2\to\R$, given by $\lambda(i,j)=i-\eps j$ is injective. A $\lambda$-monotone lattice path of length $N$ in $\Delta$ is a sequence of points $P_0,...,P_N\in\Delta\cap\Z^2$ strictly ordered by the values of $\lambda$ and such that $\lambda(P_0)=\min\lambda$, $\lambda(P_n)=\max\lambda$. Then
$$\deg V_{\Delta,g}(A_2)=\sum_{T\in\Trop_{\Delta,g}(A_2,\bx)}\mt_c(T)=\sum_{\pi\in\Pi_{\lambda,N}(\Delta)}\mt_c(\pi),$$
where $\Pi_{\lambda,N}(\Delta)$ is the set of $\lambda$-monotone lattice paths of length $N$ in $\Delta$, and $\mt_c(\pi)$ is a nonnegative integer which can be computed via an explicit finite combinatorial algorithm.
\end{theorem}

We note that the existence of a ``cuspidal" modification of Mikhalkin's lattice path algorithm is not obvious, and the observation made in \cite[Proof of Proposition 2.6]{IKS1} appears to be crucial for the proof of such an existence.

An ultimate feature of the tropical enumerative geometry is that it allows one to count real curves. So far no enumerative invariant counting real cuspidal curves is known (except for the local case \cite[Section 4]{Shloc}). A reasonable question is to estimate the number of real cuspidal curves in properly chosen families. In the present paper we suggest the following bound.

\begin{theorem}\label{tC}
For every fixed $r\ge 1$ and $d\ge2r+3$, there exists a generic real $2r$-dimensional linear subsystem in $|{\mathcal O}_{\PP^2}(d)|$ that contains at least $c_r(d)$ real curves with $r$ real cusps as their only singularities, where $c_r(d)$ is a positive function satisfying
$$c_r(d)=\frac{(3d^2)^r}{r!}+O(d^{2r-1})\quad\text{as}\quad d\to\infty\ .$$
\end{theorem}

Note that $c_r(d)$ differs by a constant factor from the total number $\qquad$ $\frac{1}{r!}(12d^2)^r+O(d^{2r-1})$ of complex $r$-cuspidal curves in a generic
$2r$-dimensional linear subsystem of $|{\mathcal O}_{\PP^2}(d)|$.

Enumeration of plane curves with non-nodal singularities has attracted the attention of many experts. This enumerative problem can be traced back to the classically known formula, $12(d-1)(d-2)$, for the degree of the family of irreducible plane curves of degree $d$ with an ordinary cusp as its only singularity. Various ways to derive this formula can be found, for instance, in \cite[Example in page 3174]{A},  \cite[\S10]{Ka}, and \cite[Proposition 1.1]{Ke}; it can be derived tropically as well as we demonstrate in Section \ref{sec-ex}(1). Closed formulas for degrees of varieties parameterizing curves having several
singularities (including a cusp) with total Milnor number $\le 5$ can be found in \cite[\S10]{Ka}.
Rational curves with one cusp in the plane (or on a del Pezzo surface) were enumerated in
\cite{Zi} and \cite{Bi} via the study of the cuspidal stratum in the moduli spaces of stable maps of
rational curves to the plane (or to a del Pezzo surface). Enumeration of plane unicuspidal curves of positive genera was discussed in \cite[Section 3.1]{Va}. We point out that \emph{none} of the methods used in the cited works applies to enumeration of \emph{real} cuspidal curves. We also expect that the techniques used in the present paper will allow to enumerate curves on toric surfaces having more complicated singularities like $A_3$, $D_4$, or $E_6$. Another perspective is to convert the lattice path algorithm from Section \ref{sec-lp} into a Caporaso-Harris type formula in the presence of a cuspidal singularity. Namely, the tropical curves counted by the lattice path algorithm naturally split into \emph{floors} in the sense of \cite{BM} which encode degenerations occurring when a given divisor on the ambient surface splits off.

The structure of the paper is as follows. Section \ref{sec-pr} introduces necessary elements of tropical geometry and proves auxiliary statements. In Section \ref{sec-cor}, we prove the correspondence theorem between algebraic and tropical unicuspidal curves (Theorem \ref{tuc1} describes the tropical limits of unicuspidal algebraic curves, Theorem \ref{t2} gives a formula for the multiplicity of a unicuspidal tropical curve). In Section \ref{sec-lp}, we suggest a ``cuspidal" modification of Mikhalkin's lattice path algorithm (Theorem \ref{t4}). At last, Section \ref{sec-ex} is devoted to enumeration of real cuspidal curves (Theorem \ref{tex}). In Appendix we discuss multiplicities of the tropical cuspidal fragments of the two remaining types.

\smallskip

{\bf Acknowledgements.} The work at this paper has been supported by the Israeli Science Foundation grants no. 176/15 and 501/18. The first author was also partially supported by the European Research Council Advanced grant 338809, while the second author enjoyed a partial support from the Bauer-Neuman Chair in Real and Complex Geometry.
A part of the work has been performed during the stay of the second author at the Mittag-Leffler Institute, Stockholm, Sweden. The second author is very grateful to the MLI for hospitality and excellent working conditions. Finally, we express our gratitude to the anonymous referees for the careful reading of the manuscript and the criticism which helped us to correct mistakes, fill in gaps, and improve the presentation. Special thanks are due to D. Kerner, who explained to us the asymptotics of the count of multi-cuspidal plane curves, and R. Piene, who attracted our attention to the paper \cite{Qv}.

\section{Tropical curves and tropical limits}\label{sec-pr}

Here we remind some basic stuff in tropical geometry adapted to our setting, introduce notation, and
provide related auxiliary statements. Almost all the details can be found in
\cite{IMS,M,Shustin2005}.

\subsection{Definitions and notation}\label{sec-nota}
{\bf(1)} We use the complex field $\C$ and the field $\K$ of locally convergent complex Puiseux series.
For an element $a=\sum_{r\ge r_0}a_rt^r\in\K$, $a_{r_0}\in\C^*$, denote
\begin{equation}\val(a)=-r_0,\quad \ini(a)=a_{r_0},\quad \Ini(a)=a_{r_0}t^{r_0}\ .\label{emn14}\end{equation}
Let $F=\sum_{\omega\in P\cap\Z^n}a_\omega\uz^\omega\in\K[\uz]$, $\uz=(w_1,...,w_n)$, have Newton polytope $\Delta$. It
yields a {\it tropical polynomial}
$$N(\ux)=N_F(\ux)=\max_{\omega\in \Delta\cap\Z^n}(\langle\omega,\ux\rangle+\val(a_\omega)),\quad N:\R^n\to\R\ ,$$
and its Legendre dual, {\it valuation function} $\nu=\nu_N:\Delta\to\R$, whose graph defines a subdivision $\Sigma_\nu$
of $\Delta$ into linearity domains which all are convex lattice polytopes. One can write
$$F(\uz)=\sum_{\omega\in \Delta\cap\Z^n}(a^0_\omega+O(t^{>0}))t^{\nu(\omega)}\uz^\omega\ ,$$ where
$a^0_\omega\in\C$, $a^0_\omega\ne0$ for all $\omega$ vertices of the subdivision $\Sigma_\nu$. Given a face $\delta$ of
the subdivision $\Sigma_\nu$, we write
$$F^\delta(\uz)=\sum_{\omega\in\delta\cap\Z^n}a_\omega\uz^\omega,\quad\ini(F^\delta)(\uz)=
\sum_{\omega\in\delta\cap\Z^n}a^0_\omega\uz^\omega\in\C[\uz]\ ,$$
$$\Ini(F^\delta)(\uz)=\sum_{\omega\in\delta\cap\Z^n}a^0_\omega t^{\nu(\omega)}\uz^\omega
\in\K[\uz]\ .$$

\smallskip

{\bf(2)} All lattice polyhedra we consider lie in Euclidean spaces $\R^N$ with fixed integral lattices
$\Z^N\subset\R^N$. We consider these spaces as $\Z^N\otimes\R$ and denote them by $\R^N_\Z$.

For a convex lattice polygon $\Delta\subset\R^2_\Z$, we denote
by $\Tor(\Delta)$ the complex toric surface associated with $\Delta$
and by $\Tor^*(\Delta)\subset\Tor(\Delta)$ the big torus (the dense orbit of the torus action). Next, denote by ${\mathcal L}_\Delta$ the tautological line bundle over $\Tor(\Delta)$, by
$|{\mathcal L}_\Delta|$ the linear system generated by the non-zero global sections (equivalently, by the monomials
$x^iy^j$, $(i,j)\in \Delta\cap\Z^2$). We also use the notation ${\mathcal L}_\Delta(-Z)$ for
${\mathcal L}_\Delta\otimes{\mathcal J}_{Z/\Tor(\Delta)}$, where ${\mathcal J}_{Z/\Tor(\Delta)}$ is the ideal sheaf of a zero-dimensional subscheme $Z\subset\Tor(\Delta)$.

Let $\bn:\widehat C\to X$ be a non-constant morphism of a complete smooth irreducible curve $\widehat C$ to a toric surface $X$. We call it {\it peripherally unibranch} if, for any
toric divisor $D\subset X$, the divisor $\bn^*(D)\subset\widehat C$ is concentrated at one point.
Respectively we call the curve $C=\bn_*\widehat C$ {\it peripherally smooth},
if it is smooth at its intersection points with each toric divisor.

\smallskip
{\bf(3)}  For a vector $v\in\Z^n$, resp. a lattice segment $\sigma\subset\R^n_\Z$, we denote
its lattice length by $\|v\|_\Z$,
resp. $\|\sigma\|_\Z$. More generally, for an $m$-dimensional lattice polytope
$\Delta\subset\R^n_\Z$, $n\ge m$, we denote by $\|\Delta\|_\Z$ its $m$-dimensional lattice volume (i.e., the ratio
of the Euclidean $m$-dimensional volume of $\Delta$ and the minimal
Euclidean volume of a lattice simplex inside the affine $m$-dimensional subspace of $\R^n_\Z$ spanned by $\Delta$).

\smallskip{\bf(4)} We always take the standard basis in $\R^2_\Z$ and identify $\Lambda^2(\R^2_\Z)\simeq\R_\Z$ by
letting $$\oa\wedge\ob=\det\left(\begin{matrix}a_1&a_2\\ b_1&b_2\end{matrix}\right),
\quad\oa=(a_1,a_2),\ \ob=(b_1,b_2)\in\R^2_\Z\ .$$

\subsection{Plane tropical curves}\label{sec-tropcur}
{\bf(1)} A {\it plane tropical curve} is a pair $(\Gamma,h)$, where
\begin{itemize}\item
$\Gamma$ is a finite
connected metric graph, whose set $\Gamma^0$ of vertices is nonempty and does not contain
univalent vertices\footnote{We allow bivalent vertices, and we call a tropical curve trivalent if there are no vertices of valency $>3$.},
the set of edges $\Gamma^1$ contains a
subset $\Gamma^1_\infty\ne\emptyset$ consisting of edges isometric to $[0,\infty)$ (called ends), while
$\Gamma^1\setminus\Gamma^1_\infty$ consists of edges isometric to compact segments in $\R$ (called finite edges);
furthermore, each vertex $V\in\Gamma$ and each edge $E\in\Gamma^1$ are equipped with a nonnegative integral genus $g(V)$, resp. $g(E)$;
\item $h:\Gamma\to\R^2$ is a 
continuous map such that
$h$ is
affine-integral on each edge of $\Gamma$ in the length coordinate, and it is nonconstant on at least one edge of
$\Gamma$;
furthermore, at each vertex $V$ of $\Gamma$, the balancing condition holds
$$\sum_{E\in\Gamma^1,\ V\in E}
\oa_V(E)=0\ ,$$ where $\oa_V(E)$ is the image under the differential $D(h\big|_E)$
of the unit tangent vector to $E$ emanating from
its endpoint $V$.
\end{itemize}
We call $\oa_V(E)$ the {\it directing vector} of $E$ (centered at $V$).

Given a vertex $V\in\Gamma^0$ such that the $h$-images of the edges emanating from $V$ do not lie in one line, the directing vectors $\oa_V(E)$ over all edges $E$ incident to $V$, being positively rotated by $\frac{\pi}{2}$ form a nondegenerate convex lattice polygon
$P(V)$, which we call {\it dual} to $V$.
The multi-set $\deg(\Gamma,h)\subset\R^2_\Z$ of vectors $\{\oa_V(E)\ne0\ : \ E\in\Gamma^1_\infty,
\ V\in\Gamma^0,\ V\in E\}$ is called the {\it degree} of $(\Gamma,h)$. It is easy to see that $\deg(\Gamma,h)$
is nonempty and {\it balanced} (i.e., the vectors of $\deg(\Gamma,h)$ sum up
to zero. The degree $\deg(\Gamma,h)$ is called {\it primitive} if it consists of primitive vectors (i.e., vectors of lattice length $1$).
Positively rotated by $\frac{\pi}{2}$, the vectors of $\deg(\Gamma,h)$
can be combined into a convex lattice polygon $\Delta=\Delta(\Gamma,h)$, called
the {\it Newton polygon} of $(\Gamma,h)$. In this case, we say that the degree $\deg(\Gamma,h)$ {\it induces} the polygon $\Delta$. The degree $\deg(\Gamma,h)$ is called {\it nondegenerate} if $\dim \Delta(\Gamma,h)=2$.

To each edge $E\in\Gamma^1$ we assign the {\it weight} $$w(E)=\big\|\oa_V(E)\big\|_\Z\ .$$

An enhanced plane tropical curve is a triple $(\Gamma,h,\widehat g)$, where $\widehat g:\Gamma^0\cup\Gamma^1\to\Z$ is a nonnegative function. Define the {\it genus} of $(\Gamma, h,\widehat g)$ as
$$g(\Gamma,h,\widehat g)=b_1(\Gamma)+\sum_{E\in\Gamma^1}\widehat g(E)+\sum_{V\in\Gamma^0}\widehat g(V)\ .$$
For the sake of notation, we skip the symbol $\widehat g$ when the genus function $\widehat g$ vanishes.

A vertex $V\in\Gamma^0$ is called {\it flat} if
$$\dim\Span_\R\left\{\oa_V(E)\ :\ E\in\Gamma^1,\ V\in E\right\}\le1\ .$$

\smallskip{\bf(2)} In the preceding notation, the image $h(\Gamma)\subset\R^2$ is a
closed rational finite one-dimensional polyhedral complex without univalent and bivalent vertices.
It can be converted into a plane tropical curve as follows. Each edge
$e\subset h(\Gamma)$ is assigned a weight $w(e)$ which is the sum of the weights of the (non-contracted) edges of
$\Gamma$ intersecting $h^{-1}(x)$, where $x\in e$ is a generic point. We define a parametrization $\widehat h:\Gamma_h\to h(\Gamma)$, which is a homeomorphism and is affine-integral on each edge $E$ of $\Gamma_h$ with $\oa_V(E)=w(e)\cdot\oa$, where $e=h(E)$, and $\oa$ is a primitive integral vector parallel to $e$. The balancing condition evidently holds, and we call the resulting tropical curve $h_*(\Gamma)=(\Gamma_h,\widehat h)$
an {\it embedded plane tropical curve}.

We recall that the graph $\Gamma_h$ is a corner locus of a tropical polynomial
$N:\R^2\to\R$ with the Newton polygon $\Delta=\Delta(\Gamma,h)$. The Legendre dual function $\nu_N:\Delta\to\R$
defines a subdivision $\Sigma$ of $\Delta$ into convex lattice polygons, and this subdivision is completely
determined by $h_*(\Gamma)$. There is a duality reversing the incidence relation: The polygons of $\Sigma$ are in bijection with the vertices of $\Gamma_h$ so that the number of sides of a polygon in $\Sigma$ equals to the valency of the dual vertex of $\Gamma^0_h$, while the edges of $\Sigma$ are in bijection with the edges of $\Gamma_h$
so that the lattice length of an edge of $\Sigma$ equals the weight of the dual edge of $\Gamma_h$.

\smallskip
{\bf(3)}
An {\it enhanced marked plane tropical curve} is a tuple $(\Gamma,\bp,h,g)$, where $\bp$ is an ordered subset of
$n\ge1$ distinct points of $\Gamma$. Suppose that $h$ is injective on $\bp$. In this case, we introduce
a {\it reduced enhanced marked plane tropical curve} $(\Gamma_{red},\bp_{red},h_{red},\widehat g_{red})$
in the following manner.

First, we contract one by one the edges of $\Gamma$ on which $h$ is constant and modify the genus function as follows: when we contract an edge $E$ joining different vertices $V_1,V_2$, then the vertex obtained from gluing of $V_1$ and $V_2$ gets the genus $\widehat g(V_1)+\widehat g(V_2)+\widehat g(E)$; if we contract a loop $E$ incident to a vertex $V$, then $V$ gets the genus $\widehat g(V)+\widehat g(E)+1$. Any bivalent vertex obtained by contracting a rational edge should be removed.

Next, if there are no two edges $E',E''\in\Gamma^1$ incident to the same vertex and such that $h(E')\cap h(E'')$ is infinite (i.e., a compact segment or a ray), we define
$(\Gamma_{red},\bp_{red},h_{red},\widehat g_{red})=(\Gamma,\bp,h,\widehat g)$. Otherwise, we perform the following procedure until we end up
with a curve possessing the above property, which then will be denoted $(\Gamma_{red},\bp_{red},h_{red},\widehat g_{red})$.
Namely, given two edges $E',E''\in\Gamma^1$ incident to
a vertex $V\in\Gamma^0$ such that $h(E')\cap h(E'')$ is infinite, we perform the following operations
(cf. Figure \ref{fig-red}):
\begin{itemize}\item we respectively identify
$h^{-1}(h(E')\cap h(E''))\cap E'$ with $\qquad\qquad\qquad$ \mbox{$h^{-1}(h(E')\cap h(E''))\cap E''$} into one edge $\widehat E$,
\item set
$\oa_V(\widehat E)=\oa_V(E')+\oa_V(E'')$ and
\mbox{$w(\widehat E)=
w(E')+w(E'')$},
\item if $h(E')\subsetneq h(E'')$, $V'\ne V$ the second vertex of $E'$, $\widehat V$ the second vertex of $\widehat E$, and $\widehat E'$ is the closure of $E''\setminus h^{-1}(h(E')\cap h(E''))$, we set $\widehat g_{red}(\widehat E)=\widehat g(E')$, $\widehat g_{red}(\widehat V)=\widehat g(V')$, and $\widehat g_{red}(\widehat E')=\widehat g(E'')$,
    \item if $h(E')=h(E'')$ is a ray, we set $\widehat g_{red}(\widehat E)=\widehat g(E')+\widehat g(E'')$,
    \item if $h(E')=h(E'')$ is a compact segment, $V'\ne V''$ the second vertices of $E'$, $E''$, respectively, we set $$\widehat g_{red}(\widehat E)=\widehat g(E')+\widehat g(E''),\quad \widehat g_{red}(\widehat V)=\widehat g(V')+\widehat g(V'')\ ,$$
        \item if $h(E')=h(E'')$ is a compact segment, $V'$ a common second vertex of $E',E''$,
        we set $\widehat g_{red}(\widehat E)=\widehat g(E')+\widehat g(E'')+1$, $\widehat g_{red}(\widehat V)=\widehat g(V')$,
        \item each bivalent vertex should be removed, while the two edges incident to this vertex glue up into one edge, whose genus sums up the genera of the vertex and of the two glued edges,
        \item those vertices and edges which are not involved in the above procedure keep their genera.
\end{itemize}

\begin{figure}
\setlength{\unitlength}{1.0mm}
\begin{picture}(95,45)(-20,0)
\thinlines

\put(0,7){\line(1,1){5}}\put(0,17){\line(1,-1){5}}
\put(5,12){\line(1,0){20}}\put(25,12){\line(1,1){5}}
\put(25,12){\line(1,-1){5}}
\put(0,32){\line(1,1){5}}\put(0,42){\line(1,-1){5}}
\put(5,37){\line(1,0){10}}\put(15,38){\line(1,0){10}}\put(15,36){\line(1,0){10}}
\put(25,37){\line(1,1){5}}\put(25,37){\line(1,-1){5}}
\put(55,7){\line(1,1){5}}\put(55,17){\line(1,-1){5}}
\put(60,12){\line(1,0){20}}\put(70,7){\line(0,1){5}}\put(70,12){\line(1,1){5}}
\put(80,7){\line(0,1){5}}\put(80,12){\line(1,1){5}}
\put(55,32){\line(1,1){5}}\put(55,42){\line(1,-1){5}}
\put(60,38){\line(1,0){10}}\put(60,36){\line(1,0){20}}
\put(70,31){\line(0,1){4}}\put(70,38){\line(1,1){5}}
\put(80,31){\line(0,1){5}}\put(80,36){\line(1,1){5}}

\put(4,11){$\bullet$}\put(4,36){$\bullet$}\put(14,36){$\bullet$}\put(24,36){$\bullet$}
\put(24,11){$\bullet$}\put(59,11){$\bullet$}\put(69,11){$\bullet$}\put(79,11){$\bullet$}
\put(59,36){$\bullet$}\put(69,37){$\bullet$}\put(79,35){$\bullet$}
\put(14,23){$\Big\Downarrow$}\put(69,23){$\Big\Downarrow$}
\put(31,36){$\Gamma$}\put(86,35){$\Gamma$}\put(31,11){$\Gamma_{red}$}\put(86,11){$\Gamma_{red}$}
\put(7,38){$w$}\put(17,40){$w_1$}\put(17,33){$w_2$}\put(14,13){$w$}\put(12,9){$g=1$}
\put(62,40){$w_1$}\put(62,33){$w_2$}\put(62,13){$w$}\put(73,9){$w_2$}
\put(5,0){$w=w_1+w_2$}\put(60,0){$w=w_1+w_2$}

\end{picture}
\caption{Construction of $\Gamma_{red}$}\label{fig-red}
\end{figure}

It follows from the construction that $g(\Gamma)=g(\Gamma_{red})$. Denote by $\bp_{red}$ the image of $\bp$ in
$\Gamma_{red}$.

Suppose that either $\Gamma_{red}$ is trivalent, or all but one of its vertices are trivalent, while
the remaining one is four-valent. Suppose also that $\quad$ \mbox{$|\bp_{red}|=|\deg\Gamma_{red}|+b_1(\Gamma_{red})-1$}.
A curve $(\Gamma,\bp,h)$ is called {\it regular}, if
\begin{itemize}\item $\Gamma^0_{red}\cap\bp_{red}=\emptyset$ and $h_{red}$ is injective on $\bp_{red}$,
\item each connected component of $\Gamma_{red}\setminus\bp_{red}$ is either a trivalent tree containing exactly one end of $\Gamma_{red}$,
or is a tree with one four-valent vertex and exactly two ends of $\Gamma_{red}$, while the other vertices are trivalent,
or is a graph with $b_1=1$ containing exactly one end of $\Gamma_{red}$ and one four-valent vertex, which belongs to the cycle, while the remaining vertices are trivalent.
\end{itemize}
It is easy to see that, under the regularity condition, the closure $\overline K$ of each connected component $K$ of $\Gamma_{red}\setminus\bp_{red}$
possesses a unique {\it regular orientation} of all its edges, i.e., an orientation satisfying the following properties:
\begin{itemize}\item each end of $\overline K$ is oriented toward infinity,
\item each edge of $\overline K$ incident to a point $p\in\bp_{red}$ is emanating from $p$,
\item each vertex $V\in K\cap\Gamma_{red}^0$ is incident to exactly two incoming edges.
\end{itemize}
The regular orientation of $(\Gamma_{red},\bp_{red},h_{red})$ induces an orientation
on the edges of the closures of components of $\Gamma\setminus\bp$ so that among edges
incident to the same vertex $V\in\Gamma^0$ that are mapped to the same line, at most one is
incoming to $V$.

\smallskip{\bf(4)} We are going to describe geometry of plane marked tropical curves of a given genus
passing through a generic configuration of points, whose number is either maximal, or next to maximal.

\begin{lemma}\label{ltrop1} (i) Let $\Def(\Gamma,\bp,h,\widehat g)$ be the germ of the deformation space
of a plane marked enhanced tropical curve $(\Gamma,\bp,h,\widehat g)$, in which $(\Gamma,\bp)$ keeps its combinatorial type, $h$ keeps the differentials along all edges, and $\widehat g$ retains.
Then \begin{equation}\dim\Def(\Gamma_{red},\bp_{red},h_{red},\widehat g_{red})
\le|\deg(\Gamma,h)|+b_1(\Gamma_{red})-1+|\bp|\ .\label{erev1}\end{equation}

(ii) Suppose that $$n=|\bp|=|\deg(\Gamma,h)|+g(\Gamma)-1$$
and $h(\bp)\subset\R^2$ is a configuration of $n$ points in general position (cf. \cite[Section 4.2]{M}). Then
$(\Gamma,\bp,h,g)$ is trivalent, reduced, regular, $g(\Gamma)=b_1(\Gamma)$, and $\bp\cap\Gamma^0=\emptyset$.
\end{lemma}

{\bf Proof.} By \cite[Proposition 2.23]{M}, the (marked) embedded plane tropical curves having a vertex of valency $>3$ vary in a deformation space of dimension strictly less than the right-hand side of (\ref{erev1}). Furthermore, \cite[Corollary 2.24]{M} states that an embedded trivalent plane tropical curve $(\Gamma,\emptyset,h)$ varies in the space of dimension $|\deg(\Gamma,h)|+b_1(\Gamma)-1$, which proves (\ref{erev1}). At last, \cite[Lemma 4.20]{M} yields the statement of the second part of the lemma.
\proofend

\begin{lemma}\label{ltrop2}
Accepting the notation of Lemma \ref{ltrop1}(i), set also $g=g(\Gamma)$, $\Delta^*=\deg(\Gamma,h)$. Suppose that
\begin{itemize}\item $n=|\Delta^*|+g-2$,
\item $(\Gamma,h)$ either is not trivalent, or not reduced, or $g>b_1(\Gamma)$,
\item $h(\bp)\subset\R^2$ is a configuration of $n$ points in general position.
\end{itemize}

Then the following holds
\begin{enumerate}\item[(i)] either $(\Gamma,\bp,h,\widehat g)$ is trivalent, reduced, and regular, and $g=b_1(\Gamma)+1$,
\item[(ii)] or $(\Gamma,\bp,h,\widehat g)$ is reduced, regular, has one four-valent vertex, while the other
vertices are trivalent, and $g=b_1(\Gamma)$,
    \item[(iii)] or $(\Gamma,\bp,h,\widehat g)$ is not reduced, but $(\Gamma_{red},\bp_{red},h_{red})$ is trivalent, regular and
        satisfies $g=g(\Gamma_{red})=b_1(\Gamma_{red})+1$ and $|\deg(\Gamma_{red},h_{red})|=
        |\Delta^*|$,
        \item[(iv)] or $(\Gamma,\bp,h,\widehat g)$ is not reduced, but $(\Gamma_{red},\bp_{red},h_{red},\widehat g_{red})$ is trivalent, regular,
        satisfies $g=g(\Gamma_{red})=b_1(\Gamma_{red})$ and $|\deg(\Gamma_{red},h_{red})|=|\Delta^*|-1$,
        \item[(v)] or $(\Gamma,\bp,h,\widehat g)$ is not reduced, has at most one four-valent vertex and the others trivalent, and satisfies
        $g=b_1(\Gamma)$, the curve $(\Gamma_{red},\bp_{red},h_{red})$ is regular, all but one of its vertices are trivalent, one vertex is four-valent with exactly one pair of incident edges lying on one line,
        and it holds $g(\Gamma_{red})=b_1(\Gamma_{red})$.
\end{enumerate}
Furthermore, in all cases (i)-(iv), $\bp\cap\Gamma^0=\emptyset$, the intersection
of $\Gamma_{red}^0$ with the set of points $x\in\R^2$ such that $|h_{red}^{-1}(x)|>1$ is empty, and in the cases (i)-(iii), (v)
the $h$-images of distinct ends of $\Gamma$ do not lie on the same ray.
\end{lemma}

{\bf Proof.} It immediately follows from Lemma \ref{ltrop1} that $g\le b_1(\Gamma)+1$. Furthermore, Lemma \ref{ltrop1} yields that $\dim\Def(\Gamma,\emptyset,h)\le|\Delta^*|+g-2=n$. Since the germ of $\Def(\Gamma,\bp,h)$ is taken onto
an open neighborhood of $(\R^2)^n=\R^{2n}$, the marked point configuration $\bp$ must be in general position on $\Gamma$, whence $\bp\cap\Gamma^0=\emptyset$.

Suppose that $g=b_1(\Gamma)+1$. Then $n=|\Delta^*|+b_1(\Gamma)-1$. Again applying Lemma \ref{ltrop1} we get that $(\Gamma,\bp,h)$
is trivalent, reduced, and regular; hence, this case fits the conditions of item (i).

Suppose that $g=b_1(\Gamma)$, and $\Gamma$ is reduced, but not trivalent. Let us show that
$\Gamma$ cannot have a vertex of valency $v\ge5$, or a pair of vertices of valency $>3$.
For this purpose, we use the idea of the proof of \cite[Lemma 2.2]{Shustin2005}, but in a different setting. We
embed the deformation space $\Def(\Gamma,h)$ into
$\R^{e(\Gamma)}$, where $e(\Gamma)=|\Gamma^1|$, namely, for each edge $E\in\Gamma^1$
we take the value in $\R\simeq\R^2/\Span\left\{\oa_V(E)\right\}$ determined by the line containing
$h(E)$. Assuming that $\Gamma$ contains vertices as above, we shall show that
\begin{equation}
\dim\Def(\Gamma,h)\le|\Delta^*|+g-3\ .\label{emn11}\end{equation} Given a generic vector $\oa\in\R^2$, orient the edges of $\Gamma$
so that their $h$-images form acute angles with $\oa$. This, in particular defines a partial order on $\Gamma^0$,
which we extend to a linear order. By changing $\oa$ accordingly, we can suppose that, for some vertex $V^*\in\Gamma^0$ of valency
$v(V^*)\ge4$, there are at least two edges merging to $V^*$ and at least two edges emanating from $V^*$. Then we estimate
$\dim\Def(\Gamma,h)$ from above as follows: \begin{itemize}\item count all ends oriented from the infinity,
\item go through the set $\Gamma^0$ ordered as above and, at each vertex $V\in\Gamma^0$, we do not add new parameters
if there are at least two edges merging to $V$, and we add one new parameter corresponding to some of the emanating edges if there is only one edge merging to $V$.
\end{itemize} After that we perform the same estimation with respect to the orienting vector $-\oa$, sum up
these two bounds, and using the above assumption on the vertices of $\Gamma$ and the Euler characteristic relation
\begin{equation}2|\Delta^*|+2b_1(\Gamma)-2\ge|\Delta^*|+\sum_{V\in\Gamma^0}(v(V)-2)\ ,\label{emis1}\end{equation} obtain
$$2\dim\Def(\Gamma,h)\le|\Delta^*|+(|\Gamma^0|-1)$$
\begin{equation}=
|\Delta^*|+\sum_{V\in\Gamma^0}(v(V)-2)-\sum_{V\in\Gamma^0}(v(V)-3)-1\ .\label{euc12}\end{equation}
Hence, $$2\dim\Def(\Gamma,h)\overset{\text{(\ref{emis1}), (\ref{euc12})}}{\le}(2|\Delta^*|+2b_1(\Gamma)-2)-3=2|\Delta^*|+2b_1(\Gamma)-5$$
which yields (\ref{emn11}) and thereby excludes the above assumption on the
valency of vertices of $\Gamma$. That is, $\Gamma$ has one four-valent vertex while the others are trivalent.
The regularity follows immediately, since otherwise one would encounter either a bounded component of
$\Gamma\setminus\bp$, hence a restriction to the position of the points $h(\bp)$, or a simply connected component
of $\Gamma\setminus\bp$ with one four-valent vertex and one end, hence again a restriction to
the position of the points $h(\bp)$ in contradiction to the generality assumption.

If $(\Gamma,h)$ is not reduced, but $\Gamma_{red}$ is trivalent, then either $(\Gamma_{red},\bp_{red},h_{red})$ is obtained by
collapsing at least one cycle of $\Gamma$, thus, $b_1(\Gamma_{red})\le b_1(\Gamma)-1$, or by merging at least one pair of
ends of $\Gamma$, thus, $|\deg(\Gamma_{red},h_{red})|\le|\Delta^*|-1$. Then the above argument leads either to the case (iii), or to (iv).

Suppose that $(\Gamma,h)$ is not reduced, and $\Gamma_{red}$ is not trivalent. Again the above considerations leave the only possibility of the assertion (v).

\smallskip

Now we prove statement (2).
The condition $\bp\cap\Gamma^0=\emptyset$ is immediate.
The remaining conditions we demonstrate in the most difficult case (ii), while the other cases can be
settled in the same manner. So, assume that there are $V\in\Gamma^0$, $x\in\Gamma\setminus
\Gamma^0$ such that $h(V)=h(x)$. We identify $V$ and $x$ into a vertex $V'$ and obtain a curve
$\Gamma'$ with $b_1(\Gamma')=b_1(\Gamma)+1$, which either has a $5$-valent vertex and a $4$-valent one, or has a $6$-valent vertex, while the other vertices are trivalent. In both the cases we shall derive (\ref{emn11}).

In the former case, we apply the estimation procedure
with an orienting vector $\oa$ such that the $4$-valent vertex has two incoming and two outgoing incident edges.
Note that the $5$-valent vertex has at least two incoming and at least two outgoing edges. Then
$$2\dim\Def(\Gamma,h)\le2\dim\Def(\Gamma',h')$$
$$\overset{\text{cf. (\ref{euc12})}}{=}|\Delta^*|+\sum_{V\in(\Gamma')^0}(v(V)-2)-\sum_{V\in(\Gamma')^0}(v(V)-3)-2$$
$$=2|\Delta^*|+2g-5\quad\Longrightarrow\quad\dim\Def(\Gamma,h)\le|\Delta^*|+g-3\ .$$

In the latter case, we pick a generic point $O\in\R^2$ close to $h(V')$, order the vertices of $\Gamma'$ by the
growing distance of their $h$-images from $O$, and respectively orient the edges of $\Gamma'$ (in particular,
all ends are oriented to infinity).
Denote by $(\Gamma')^0_1$, $(\Gamma')^0_2$ the sets of trivalent vertices with one or two
emanating incident edges, respectively. Similarly to the previous estimations we have
$$\dim\Def(\Gamma',h')\le2+|(\Gamma')^0_2|\ .$$ Since
$$|(\Gamma')^0_1|+2|(\Gamma')^0_2|=|(\Gamma')^1|-6\ ,$$
$$|(\Gamma')^0_1|+|(\Gamma')^0_2|+1+|\deg(\Gamma',h')|-|(\Gamma')^1|=1-g'=-g\ ,$$ we get
$$|(\Gamma')^0_2|=|\deg(\Gamma',h')|+g-5\ ,$$ and hence
$$\dim\Def(\Gamma',h')\le3+(|\deg(\Gamma',h')|+g-5)=|\deg(\Gamma',h')|+g-3$$
as required.
\proofend

\begin{definition}\label{dtrop1}
Let $(\Gamma,h)$ be a plane tropical curve, and let $V\in\Gamma^0$ be a trivalent vertex, which is not incident to
contracted edges.
Define the Mikhalkin multiplicity of $V$ by
$$\mu(V)=\big|\oa_V(E_1)\wedge\oa_V(E_2)\big|\ ,$$
where $E_1,E_2\in\Gamma^1$ are some two incident to $V$ edges. Recall that by \cite[Proposition
6.17]{M} (see also \cite[Lemma 3.5]{Shustin2005}), $\mu(V)$ equals the number of peripherally unibranch rational curves in $|{\mathcal L}_{P(V)}|$ given by polynomials with Newton triangle $P(V)$ and coefficients $1$ at the vertices of $P(V)$.
\end{definition}

\subsection{Tropicalization of algebraic curves over a non-Archimedean field}\label{sec-trop}
\subsubsection{Embedded tropical limit}\label{sec-tl}
Let $\Delta^*\subset\Z^2\setminus\{0\}$ be a nondegenerate, primitive, balanced multiset which induces a convex lattice polygon $\Delta$, and let $C\in|{\mathcal L}_{\Delta}|_\K$ be a reduced,
irreducible curve of genus $g$, which does not hit intersection points of toric divisors. In particular, it can be given by an
equation
$$F(x,y)=\sum_{(i,j)\in\Delta\cap\Z^2}t^{\nu'(i,j)}(a^0_{ij}+O(t^{>0}))x^iy^j\ ,$$
where $a^0_{ij}\in\C$, $(i,j)\in \Delta\cap\Z^2$, and $a^0_{ij}\ne0$ if the coefficient
of $x^iy^j$ in $F$ does not vanish (for instance, when $(i,j)$ is a vertex of $\Delta$).
We then define a convex, piecewise-linear function $\nu:\Delta\to\R$, whose graph is the lower part of
the $\conv\{(i,j,\nu'(i,j)),\ (i,j)\in \Delta\cap \Z^2\}\subset\R^3_\Z$.
Via a parameter change $t\mapsto t^M$, we can make
$\nu(\Delta\cap\Z^2)\subset\Z$. Denote by $\Sigma$ the subdivision
of $\Delta$ into linearity domains of $\nu$, which are convex lattice polygons $\Delta_1,...,\Delta_m$. We then have
\begin{equation}F(x,y)=\sum_{(i,j)\in \Delta\cap\Z^2}t^{\nu(i,j)}(c^0_{ij}+O(t^{>0}))x^iy^j\ ,\label{e-new60}\end{equation}
where $c^0_{ij}\ne0$ for $(i,j)$ a vertex of some of the $\Delta_1,...,\Delta_m$. This data defines a flat family
$\Phi:{\mathfrak X}\to\C$, where
${\mathfrak X}=\Tor(\OG(\nu))$ and
$$\OG(\nu)=\{(i,j,c)\in\R^3_\Z\ :\ (i,j)\in \Delta,\ c\ge\nu(i,j)\}$$ is the overgraph of $\nu$, the central fiber
${\mathfrak X}_0=\Phi^{-1}(0)$ splits into the union of toric surfaces $\Tor(\Delta_k)$, $1\le k\le m$, and the other fibers are isomorphic to
$\Tor(\Delta)$. The evaluation of the parameter $t$ turns the given curve $C$ into an inscribed family of curves
\begin{equation}C^{(t)}\subset{\mathfrak X}_t,\quad C^{(t)}\in|{\mathcal L}_{\Delta}|,\quad t\in(\C,0)\setminus\{0\}\ ,
\label{euc2}\end{equation} (where $(\C,0)$ always means a sufficiently small disc in $\C$ centered at zero)
which closes up to a flat family over $(\C,0)$ with the central element
$$C^{(0)}=\bigcup_{k=1}^mC^{(0)}_k\ ,$$ where
$$C^{(0)}_k=\left\{F^{(0)}_k(x,y):=\sum_{(i,j)\in\Delta_k\cap\Z^2}c^0_{ij}x^iy^j=0\right\}\in|{\mathcal L}_{\Delta_k}|,
\ 1\le k\le m\ .$$
The function $\nu:\Delta\to\R$ defines an embedded plane tropical curve $\Tr(C)$ in the sense of Section \ref{sec-tropcur}, part (2). Its support is the closure of the valuation image of $C$.

We define the {\it embedded tropical limit} of $C$ to be the collection $\qquad\quad$
$(\Tr(C),\{(\Delta_k,C^{(0)}_k)\}_{k=1,...,m})$, where the pairs $(\Delta_1,C^{(0)}_1)$, ..., $(\Delta_m,C^{(0)}_m)$ are called
{\it limit curves}, cf. \cite[Section 2]{Shustin2005}\footnote{When it is clear which toric surfaces are associated with the curves $C^{(0)}_1,...,C^{(m)}$, we write $C^{(0)}_i$ instead of $(\Delta_i,C^{(0)}_i)$.}.

\subsubsection{Parameterized tropical limit}\label{sec-ptl}
{\bf(1)} In the notation of the preceding section, let $\bn:\widehat C\to C$ be the normalization, or, equivalently, the family
\begin{equation}\bn_t:\widehat C^{(t)}\to C^{(t)}\hookrightarrow{\mathfrak X}_t,\quad t\in(\C,0)\setminus\{0\}\ ,
\label{euc3}\end{equation} where
each $\widehat C^{(t)}$ is a smooth curve of genus $g$ (cf. \cite[Theorem 1, page 73]{Tei} or \cite[Proposition 3.3]{ChL}). We also assume that $\widehat C$
contains a configuration $\bw$ of $n$ marked points which form $n$ disjoint families $w_i(t)\in\widehat C^{(t)}$, $t\in(\C,0)\setminus\{0\}$,
projecting to disjoint families in $C^{(t)}$ (also denoted $w_i(t)$, no confusion will arise) that avoid singularities, and, furthermore,
the valuation image of $\bw$ consists of $n$ distinct points of $\Tr(C)$. More explicitly, we suppose that $w_i(t)=(\alpha_it^{r_i},\beta_it^{s_i})$ with generic $\alpha_i,\beta_i\in\C^*$ and with $r_i,s_i\in\Q$ chosen to that the points $(-r_i,-s_i)\in\R^2$, $i=1,...,n$, are in tropically general position (cf., \cite[Section 4.2]{M}).

The family
(\ref{euc3}) admits (after a suitable untwist $t\mapsto t^M$) a flat extension to the center $0\in(\C,0)$
with the central element $\bn_0:\widehat C^{(0)}\to{\mathfrak X}_0$, where $\widehat C^{(0)}$ is a connected nodal curve of arithmetic genus $g$ (see, for instance \cite[Theorem 1.4.1]{AV}), none of whose components is entirely mapped to a toric divisor $\Tor(e)$, $e\subset\partial\Delta$,
and such that $(\bn_0)_*\widehat C^{(0)}=C^{(0)}$. We assume also that the sections $w_i(t)$, $i=1,...,n$, close up at $n$ distinct points of $C^{(0)}$.

With the central fiber one can associate an enhanced plane marked tropical curve $(\Gamma,\bp,h,\widehat g)$ as defined in \cite[Section 2.2.1]{Ty}. In particular (all other details can be found in \cite[Section 2]{Ty}),
\begin{itemize}\item The vertices of $\Gamma$ bijectively correspond to components of $\widehat C^{(0)}$, the finite edges of $\Gamma$ bijectively correspond to the intersection points of distinct components of $\widehat C^{(0)}$, and the infinite edges of $\Gamma$ correspond to the points of $\widehat C^{(0)}$ mapped to toric divisors $\Tor(e)$ such that $e\subset\partial\Delta$;
\item the genus function $\widehat g$ assigns zero to all edges of $\Gamma$, while $\widehat g(V)$ for a vertex $V\in\Gamma^0$ takes the value $g(\widehat C_V^{(0)})$, where $\widehat C_V^{(0)}$ is the corresponding component of $\widehat C^{(0)}$.
\end{itemize}
Thus, we obtain the
{\it parameterized tropical limit} (briefly, {\it PTL}) of $\qquad$ \mbox{$(\bn:\widehat C\to C,\bw)$} to be the pair consisting of the above enhanced plane marked tropical curve $(\Gamma,\bp,h,\widehat g)$ and the marked parameterized
complex curve $\quad$ $(\bn_0:\widehat C^{(0)}\to C^{(0)},\bw(0))$ (in ${\mathfrak X}_0$).

\smallskip{\bf(2)} Our purpose is to explicitly describe parameterized tropical limits that can be bijectively associated with the curves $\bn:\widehat C\to C\hookrightarrow\Tor_\K(\Delta)$ such that $\bw\subset C$ and $C\in V_{\Delta,g}(A_2)$.
To achieve this goal, we introduce now a simplified PTL (which is quite close to the embedded tropical limit as in Section \ref{sec-tl}). Later, in Section \ref{sec-cor}, we extend (in a constructive way) the simplified PTL to a {\it modified PTL}, which will contain all the information on the PTL as defined in part (1) above and will associate with each contracted component $\widehat C^{(0)}_i$ of $\bn_0:\widehat C^{(0)}\to C^{(0)}$ a non-constant map of $\widehat C^{(0)}_i$ to a certain toric surface.

To obtain the simplified PTL, we first contract all components of $\widehat C^{(0)}$ that are mapped to points by $\bn_0$, and we obtain a map of the quotient $\bn'_0:(\widehat C^{(0)})'\to C^{(0)}$. Then we construct an enhanced plane tropical curve $(\Gamma',\bp',h',\widehat g')$ starting with $(\Gamma,\bp,h,\widehat g)$:
\begin{itemize}\item we contract one-by-one the edges of $\Gamma$ incident to the vertices dual to the contracted components of $\widehat C^{(0)}$, and we each time modify the genus function along the recipe of Section \ref{sec-tropcur}, part (3), second paragraph;
\item in case a bivalent vertex of the resulting graph does not correspond to a component of
$(\widehat C^{(0)})'$, we remove this vertex gluing the two incident edges into one edge and adding the genus of the removed vertex to the genus of the obtained edge.
\end{itemize}
We make the following comments on the simplified PTL:
\begin{itemize}\item $g(\Gamma',\bp',h',\widehat g')=g(\Gamma,\bp,h,\widehat g)$ by construction;
\item $(\Gamma',\bp',h',\widehat g')$ in general differs from $(\Gamma_{red},\bp_{red},h_{red},\widehat g_{red})$
(and from $\Tr(C)$ as introduced in Section \ref{sec-tl};
for instance, it may contain fragments like the two upper ones in Figure \ref{fig-red}.
\end{itemize}

\section{Correspondence between unicuspidal algebraic curves and unicuspidal tropical curves}\label{sec-cor}

Throughout this section, $\Delta^*$ always means a nondegenerate, primitive, balanced multi-set in $\Z^2\setminus\{0\}\subset\R^2_\Z$, and $\Delta\subset\R^2_\Z$ the induced convex lattice polygon. Note that the primitivity of $\Delta^*$ allows one to uniquely recover $\Delta^*$ out of $\Delta$, in particular, $|\delta^*|=|\partial\Delta\cap\Z^2|$.

\subsection{Statement of the enumerative problem}\label{sec-prob}
Define the arithmetic genus by $p_a(\Delta)=|\Int(\Delta)\cap\Z^2|$.
Assuming that $p_a(\Delta)\ge1$, for any $0\le g<p_a(\Delta)$, denote by $V_{\Delta,g}(A_2)$
the family of reduced, irreducible curves $C\in|{\mathcal L}_{\Delta}|$ on the toric surface $\Tor(\Delta)$ that have
genus $g$ and such that the normalization $\nu:\widehat C\to C\hookrightarrow\Tor(\Delta)$ is an immersion everywhere but at one point $z\in\widehat C$, while the germ $(\widehat C,z)$ is mapped to a singular local branch of multiplicity $2$.

\begin{lemma}\label{luc1}
Suppose that $|\Delta^*|\ge5$ and there exists a subset $Q\subset(\partial \Delta\cap\Z^2)$ such that
${\mathcal Q}=\conv(Q)$ is a nondegenerate quadrangle without parallel edges.
Then $V_{\Delta,g}(A_2)\ne\emptyset$ is of pure dimension $|\partial\Delta\cap\Z^2|+g-2$, and a generic element of any component
of $V_{\Delta,g}(A_2)$ is a curve
with one ordinary cusp and $p_a(\Delta)-g-1$ nodes.
\end{lemma}

Under the hypotheses of Lemma \ref{luc1}, we can pose

\smallskip
{\bf Problem.} {\it For any $0\le g<p_a(\Delta)$, compute $\deg V_{\Delta,g}(A_2)$.}

\smallskip

For the case of $\Tor(\Delta)$ a smooth projective surface and $p_a(\Delta)-g$ obeying a certain
upper bound (for example, $p_a(\Delta)-g<2(m-1)$ for degree $m$ plane curves), the formulas for $\deg V_{\Delta,g}(A_2)$ can, in principle, be obtained
by the methods of \cite{Ka} (see \cite[Section 10.2]{Ka} with examples covering the domain
$p_a(\Delta)-d\le4$). At the other extreme, the plane rational curves of any degree $d$ with one cusp and $\frac{d(d-3)}{2}$ nodes
have been enumerated in \cite{Zi}.

We solve the problem for a wide class of toric surfaces, including the plane, quadric, toric del Pezzo, and Hirzebruch surfaces, for any divisor class, satisfying conditions of Lemma \ref{luc1}, and for all genera.

The solution is given via tropical geometry in the style of \cite{M} and \cite{Shustin2005}: namely, using the Lefschetz principle, we pass to an equivalent problem over
the field $\K$ of locally convergent complex Puiseux series,
choose $|\partial\Delta\cap\Z^2|+g-2$ points in a special position in $(\K^*)^2\subset\Tor_\K(\Delta)$, prove the correspondence theorem,
which describes the plane tropical curves (called further on {\it unicuspidal})
obtained via the tropicalization of the counted algebraic curves, and compute the
multiplicities of these tropical curves (i.e., the number of algebraic curves tropicalizing to the given tropical curve).
Moreover, we provide a suitable analogue of the Mikhalkin's lattice path algorithm allowing one, in a purely combinatorial way, to
enumerate the unicuspidal tropical curves and compute their multiplicities.

\smallskip

{\bf Proof of Lemma \ref{luc1}.}
First of all notice that any curve in $V_{\Delta,g}(A_2)$ admits a deformation into
a curve (also belonging to $V_{\Delta,g}(A_2)$) that has an ordinary cusp and $p_a(\Delta)-g-1$ nodes
as its singularities. We derive this from \cite[Theorem 1.1]{Shustin1998} which confirms the existence of the above deformation in the case of curves in the plane under the condition
\begin{equation}c\ge2g-d+2,\label{e-new55}\end{equation}
where $d$ is the degree of the curves under consideration, $c$ the degree of the dual curves. Introducing the number of virtual cusps of a reduced, irreducible plane curve $C$ (cf., \cite{Ku}) to be
$$k_{vir}(C)=\sum_{z\in\Sing(C)}(\varkappa(C,z)-2\delta(C,z))$$
(the definition of $\varkappa$- and $\delta$-invariants can be found, for instance, in \cite[Section I.3.4]{GLS}),
we can convert inequality (\ref{e-new55}) into
\begin{equation}k_{vir}\le3d-4=-KC-4,\label{e-new56}\end{equation} where $K$ is the canonical divisor class of the plane.
The proof of  \cite[Theorem 1.1]{Shustin1998} can easily be extended to the case of curves on an arbitrary toric surface with the same sufficient condition (\ref{e-new56}).
In our situation, $k_{vir}=1$ and $-KC=|\partial\Delta\cap\Z^2|\ge5$; hence, the condition holds.

We can triangulate $\Delta\setminus{\mathcal Q}$ by diagonals
(i.e., chords joining vertices of $\Delta$) that may intersect only at vertices of $\Delta$. The resulting subdivision of $\Delta$ is convex, i.e., lifts to a graph of a convex piecewise-linear function.
We, furthermore, orient the adjacency tree so that ${\mathcal Q}$ is a source, then extend the induced partial order on the polygons of the subdivision up to some linear order. Notice that, for each triangle, at most two edges of the adjacency graph are incoming.

By Lemma \ref{lquad1} below, there exists a rational peripherally unibranch curve $C_Q\in|{\mathcal L}_Q|$, which is
peripherally smooth, has a unique singular branch in $\Tor^*({\mathcal Q})$, and this branch is of type $A_2$.
Following the chosen above linear order, we extend $C$ to a collection of rational peripherally
unibranch curves for each triangle of the subdivision, each time matching the points on the common toric divisors
determined by the preceding curves. By the patchworking theorem \cite[Theorem 2.4]{Sh} (see also Section \ref{sec-pw} below) there exists a rational curve $\PP^1\to C_{\Delta}\in|{\mathcal L}_{\Delta}|$, which is
immersed everywhere but at one point of $\PP^1$ that is mapped to a
singular branch of multiplicity $2$
in the torus $\Tor^*(\Delta)$. Hence, $V_{\Delta,0}(A_2)\ne\emptyset$.

It follows from \cite[Theorem in Section 6.1]{GK} that
\begin{itemize}\item for any $1\le g<p_a(\Delta)$, a generic curve $C\in V_{\Delta,0}(A_2)$ admits a deformation into a curve $C'\in V_{\Delta,g}(A_2)$ by smoothing out any $g$ nodes, while keeping the rest of singularities,
\item for any $0\le g<p_a(\Delta)$, the family of curves of genus $g$ in $|{\mathcal L}_{\Delta}|$ having one cusp
and $p_a(\Delta)-1-g$ nodes as its singularities, is smooth of expected dimension $|\partial\Delta\cap\Z^2|+g-2$.
\end{itemize} Indeed, the sufficient condition for that asserted in \cite[Theorem in Section 6.1]{GK} requires the number of cusps to be
less than $-KC=|\partial\Delta\cap\Z^2|\ge5$. \proofend

\subsection{Tropicalization of nodal curves and modifications}\label{nod-mod}
To help the reader in better understanding of the next sections devoted to cuspidal curves, we quickly recall here the well-understood example of nodal curves, see \cite{M} as well as \cite{Shustin2005}, \cite[Chapter 2]{IMS}. Namely, we shall describe simplified parameterized tropical limit, then discuss modifications which allow one to obtain the
{\it modified parameterized tropical limit} in the sense of Section \ref{sec-ptl}, and, finally, we explain how to compute the multiplicity of the underlying parameterized tropical curve.

{\bf(1)} Let $\Delta^*\subset\Z^2\setminus\{0\}$ be a nondegenerate, primitive, balanced multiset with induces a convex lattice polygon $\Delta$, and let $C\in|{\mathcal L}_{\Delta}|_\K$ an irreducible nodal curve of genus $0\le g\le p_a(\Delta)$, which does not hit intersection points of toric divisors on $\Tor_\K(\Delta)$. Suppose that $C$ passes through a configuration $\bw$ of $n=|\partial\Delta\cap\Z^2|+g-1$ distinct points in $(\K^*)^2\subset\Tor_\K(\Delta)$, and the valuation image of $\bw$ is a tropically generic configuration $\bx$ of $n$ distinct points in $\R^2$.

Then the induced subdivision of $\Delta$ (see Section \ref{sec-tl}) consists of triangles and parallelograms. Each limit curve $(\Delta_i,C^{(0)}_i)$, corresponding to a triangle $\Delta_i$ is nodal, rational, and is peripherally unibranch (see Section \ref{sec-nota}, part (2)). Each limit curve $(\Delta_i,C^{(0)}_i)$, corresponding to a parallelogram $\Delta_i$, is given by a polynomial
\begin{equation}(\alpha x^a+\beta y^b)^p(\gamma x^c+\delta y^d)^q,\label{e-new58}\end{equation} $$\gcd(a,b)=\gcd(c,d)=1,\ ad-bc\ne0,\ \alpha\beta\gamma\delta\ne0.$$
The simplified tropical limit includes a plane marked tropical curve $(\Gamma',\bp',h')$, while all the components of $(\widehat C^{(0)})'$ are rational. The graph $\Gamma'$ has tri- and bivalent vertices. The trivalent vertices $V\in\Gamma^{\prime,0}$ correspond to the triangles $T=P(V)$ of the subdivision of $\Delta$ and to the normalization maps $\bn'_0:(\widehat C^{(0)}_i)'=\PP^1\to C^{(0)}_i\subset\Tor(T)$. The bivalent vertices occur in pairs $V_1,V_2\in\Gamma^{\prime,0}$ such that $h'(V_1)=h'(V_2)\in\R^2$, and each pair $(V_1,V_2)$ corresponds to a parallelogram $P$ in the subdivision of $\Delta$ and to the two maps (see (\ref{e-new58}))
$$\bn'_0:(\widehat C^{(0)}_i)'=\PP^1\to\{\alpha x^a+\beta y^b=0\},\quad\bn'_0:(\widehat C^{(0)}_i)'=\PP^1\to\{\gamma x^c+\delta y^d=0\},$$ which are $p$- and $q$-multiple covers ramified at the intersection points of the image with the toric divisors of $\Tor(P)$, respectively.

\begin{figure}
\setlength{\unitlength}{1.0mm}
\begin{picture}(140,110)(0,0)
\thinlines

\put(40,45){\vector(0,1){20}}\put(40,45){\vector(1,0){45}}
\put(95,45){\vector(0,1){20}}\put(95,45){\vector(1,0){45}}
\put(15,85){\vector(0,1){25}}\put(0,85){\vector(1,0){35}}
\put(60,85){\vector(0,1){25}}\put(45,85){\vector(1,0){35}}

\dashline{1}(0,25)(30,25)
\dashline{1}(0,55)(30,55)
\dashline{1}(95,115)(105,115)\dashline{1}(125,115)(135,115)
\dashline{1}(95,100)(135,100)

\dottedline{1}(95,60)(105,60)

\thicklines

\put(5,25){\line(2,-1){10}}\put(15,20){\line(2,1){10}}\put(15,10){\line(0,1){10}}
\put(10,10){\line(0,1){15}}\put(20,10){\line(0,1){15}}

\put(2,50){\line(1,1){5}}\put(2,60){\line(1,-1){5}}\put(7,55){\line(1,0){18}}
\put(25,55){\line(1,1){5}}\put(25,55){\line(1,-1){5}}\put(9,60){\line(1,-2){5}}\put(22,61){\line(-1,-3){4}}

\put(40,55){\line(2,-1){10}}\put(40,55){\line(2,1){10}}
\put(50,50){\line(0,1){10}}\put(50,50){\line(3,1){10}}\put(50,60){\line(3,1){10}}\put(60,53.3){\line(0,1){10}}
\put(60,53.3){\line(2,-1){10}}\put(60,63.3){\line(2,-1){10}}\put(70,48.3){\line(0,1){10}}
\put(70,48.3){\line(1,1){5}}\put(70,58.3){\line(1,-1){5}}

\put(95,45){\line(0,1){10}}\put(95,55){\line(2,1){10}}\put(95,45){\line(1,0){40}}\put(135,45){\line(0,1){5}}
\put(100,45){\line(1,3){5}}\put(110,45){\line(1,3){5}}\put(120,45){\line(-1,3){5}}\put(130,45){\line(-1,3){5}}
\put(125,60){\line(1,-1){10}}\put(105,60){\line(1,0){20}}

\put(5,95){\line(1,1){10}}\put(5,95){\line(1,-1){10}}\put(15,105){\line(2,-3){9.3}}\put(15,85){\line(3,2){9}}
\put(15,85){\line(0,1){20}}

\put(50,85){\line(0,1){10}}\put(55,85){\line(1,4){5}}\put(65,85){\line(-1,4){5}}\put(50,95){\line(1,1){10}}
\put(60,105){\line(2,-3){10}}\put(70,85){\line(0,1){5}}\put(50,85){\line(1,0){20}}

\put(100,110){\line(1,1){5}}\put(100,120){\line(1,-1){5}}\put(105,115){\line(1,0){20}}
\put(125,115){\line(2,3){5}}\put(125,115){\line(3,-2){7}}
\put(105,100){\line(2,-1){10}}\put(125,100){\line(-2,-1){10}}\put(115,95){\line(0,-1){10}}

\put(13,75){\rm{(a)}}\put(58,75){\rm{(b)}}\put(113,75){\rm{(c)}}
\put(13,35){\rm{(d)}}\put(58,35){\rm{(e)}}\put(113,35){\rm{(f)}}
\put(13,0){\rm{(g)}}

\put(16,105){$\ell$}\put(61,105){$\ell$}\put(92,59){$\ell$}
\put(104,41){$k-1$}\put(117,41){$k+1$}
\put(52,81){$-1$}\put(64,81){$1$}
\put(3,26){$V_1$}\put(24,26){$V_2$}
\put(6,51){$V_1$}\put(22,51){$V_2$}
\put(105,116){$V_1$}\put(121,116){$V_2$}
\put(103,101){$V_1$}\put(124,101){$V_2$}\put(116,91){$V_{mod}$}
\put(8,94){$\Delta_i$}\put(16,94){$\Delta_j$}

\end{picture}
\caption{Modifications of nodal curves}\label{fig-new59}
\end{figure}

{\bf(2)} The simplified tropical limit does not care enough information. For example, if a common point $z\in(\widehat C^{(0)}_i)'\cap(\widehat C^{(0)}_j)'$ is mapped to a point on the toric divisor $\Tor(e)$, where $e=\Delta_i\cap\Delta_j$ is a common edge of the lattice length $\ell>1$, then $z$ is a singular point of type $A_{2\ell-1}$. The following operation, which we
call {\it modification}\footnote{In \cite[Chapter 2]{IMS} and \cite{Shustin2005} these operations are called refinement, but change here the name to be consistent with the tropical counterpart called "modifications" (\cite[Chapter 5]{MR}).} resolves this singularity and shows that its deformation in the family $\bn_t:\widehat C^{(t)}\to C^{(t)}$, $t\in(\C,0)$,  yields $m-1$ extra nodes of the general member $C^{(t)}$, $t\ne0$. So, let $e=\Delta_i\cap\Delta_j$ be the common edge of two triangles $\Delta_i,\Delta_j$ of the subdivision of $\Delta$, and $\|e\|_\Z=\ell>1$. Applying a suitable automorphism of $\Z^2$, we can place $e$ on the vertical axis (see Figure \ref{fig-new59}(a)). We also can assume that
the piecewise linear function $\nu:\Delta\to\R$ defining the subdivision of $\Delta$ vanishes along $e$, and is positive outside $e$. Thus, the Laurent polynomial $F(x,y)$ given by (\ref{e-new60}) reads
\begin{equation}F(x,y)=(y-\alpha)^\ell+O(t^{>0}),\quad\alpha\in\C^*.\label{e-new66}\end{equation}
There exists a (non-toric) coordinate change of the form
\begin{equation}(x,y)=(x_1,y_1+\alpha+\zeta(t)),\quad\zeta(0)=0,\label{e-new61}\end{equation} such that, in the Laurent polynomial $F^{(1)}(x,y):=F(x,y+\alpha+\zeta(t))$, the coefficient at the monomial $y^{\ell-1}$ vanishes.
This polynomial defines a subdivision of its Newton polygon containing the fragment depicted in Figure \ref{fig-new59}(b). The triangle $T_{mod}=\conv\{(-1,0),(1,0),(0,\ell)\}$ is formed by the Newton diagrams $[(-1.0),(0,\ell)]$ and $[(1,0),(0,\ell)]$ of the point of tangency of the limit curves $C^{(0)}_i$ and $C^{(0)}_j$ to the toric divisor $\Tor([(0,0),(0,\ell)])$. The limit curve $C_{mod}\subset\Tor(T_{mod})$ is rational, has $m-1$ nodes as its only singularities, and intersects transversally the toric divisors $\Tor([(-1,0),(0,\ell)])$ and $\Tor([1,0),(0,\ell)])$.
Thus, we obtain an extra limit curve and its normalization $\PP^1\to C_{mod}$ as a new component of the parameterizing curve $\widehat C^{(0)}$, which is mapped to a point on ${\mathfrak X}_0$. Furthermore, one can make the image non-trivial by a weighted blowing-up on the point $z\in{\mathfrak X}$ with the exceptional divisor $\Tor(T_{mod})$ (cf., \cite[Section 2]{ST}).

The Mikhalkin's tropical modification can be seen as follows (cf., \cite[Section 2]{ST}): the edge $E$ of $\Gamma'$, dual to $e=\Delta_i\cap\Delta_j$ is replaced by the fragment with a trivalent vertex $V_{mod}$ shown in Figure \ref{fig-new59}(c) (which in turn is dual to the triangle $T_{mod}$). From the tropical modification one can easily restore the metric on the modified part of the parameterizing graph.

A slightly more complicated picture appears when the two vertices $V_1,V_2\in\Gamma^{\prime,0}$ are joined by a segment divided into several edges by bivalent vertices (following \cite{Shustin2005}), we call such a segment an {\it extended edge} of the tropical curve, while finite edges joining vertices of valency $\ge3$ are called {\it ordinary}). Referring for all details to \cite[Section 3.6 and Figure 3]{Shustin2005} and \cite[Section 2.5.8]{IMS}), we shortly explain here what is the modification along the extended edge of weight $\ell>1$. The planar image
of an extended edge is a linear segment transversally crossed by several other segments (see Figure \ref{fig-new59}(d)), which is dual to the fragment of the subdivision of $\Delta$ consisting of two triangles joined by a chain of parallelograms (see Figure \ref{fig-new59}(e)). We have a sequence of parallel edges $E_1,...,E_r$ with $E_1\subset \Delta_i$, $E_r\subset\Delta_j$, and the truncations of the polynomial $F(x,y)$ to $E_1,...,E_r$ all are the products of the same factor $(y-\alpha)^\ell$ and a monomial. Thus, the coordinate change (\ref{e-new61}) yields a polynomial $F(x_1,y_1+\alpha+\zeta(t))$ that induces the fragment of the subdivision of its Newton polygon as depicted in Figure \ref{fig-new59}(f), while the Mikhalkin's tropical modification results in the fragment shown in Figure \ref{fig-new59}(g). As a consequence, we resolve all non-nodal singularities on the toric divisors $\Tor(E_1),...,\Tor(E_r)$. In particular, one observes a new trivalent vertex $V_{mod}$ dual to the triangle $T_{mod}$ and an extra rational nodal limit curve $C_{mod}\subset\Tor(T_{mod})$, which all are exactly as in the preceding case. In what follows, we always assume the absence of bivalent vertices in the considered fragments of tropical curves, since adding bivalent vertices does not affect the tropicalization and patchworking statements treated in the following sections.

{\bf(3)} Given the tropical part of the simplified PTL $(\Gamma',\bp',h')$ with the configuration $\bx=h'(\bp')\subset\R^2$ in tropically general position, and the point configuration $\bw(0)\subset{\mathfrak X}_0$,
in accordance with this data, one can restore the limit curves $C^{(0)}_i\subset\Tor(\Delta_i)$, $1\le i\le m$, (see Section
\ref{sec-tl}) and the limit curves $C_{mod}\subset\Tor(T_{mod})$ modifications along all ordinary and extended edges of $\Gamma'$ of weight $>1$ in finitely many ways. The number of these ways is called the multiplicity of the tropical curve
$(\Gamma',\bp',h')$\ \footnote{It appears to be independent of the choice of the configuration $\bw(0)$.}. The patchworking statement as in \cite{M,Shustin2005,IMS,Sh} or \cite{Ty} claims that,
for any collection of limit curves $\{(\Delta_i,C^{(0)}_i),(T_{mod},C_{mod})\}$
there exists a unique curve $C\in|{\mathcal L}_\K(\Delta)|$ of genus $g$, passing through the configuration $\bw\subset(\K^*)^2$.

The curve $(\Gamma',\bp',h')$ possesses a regular orientation (cf., Section \ref{sec-tropcur}, part (3)), which defines a partial order of the set of trivalent vertices of $\Gamma'$, and the multiplicity is computed as the product of the following factors:
\begin{itemize}\item a trivalent vertex $V\in\Gamma^{\prime,0}$ with two incoming incident edges $E_1,E_2$ contributes $\mu(V)(w(E_1)w(E_2))^{-1}$ (which is the number of suitable limit curves $C^{(0)}_i\subset\Tor(\Delta_i)$, $\Delta_i=P(V)$, see \cite[Lemma 3.5]{Shustin2005});
\item an ordinary or extended edge of $\Gamma'$ of weight $\ell>1$ contributes the factor $\ell$ (which is the number of suitable limit curves $C_{mod}\subset\Tor(T_{mod})$ coming from the corresponding modification, see \cite[Lemma 3.9]{Shustin2005});
    \item a marked point $p\in\bp$ lying on an ordinary or extended edge of weight $\ell>1$ contributes the factor $\ell$, which comes from the correction of the shift (\ref{e-new61}) by the term of higher order in $t$ (see \cite[Section 2.5.9]{IMS}).
\end{itemize}
In what follows, we consider tropical curves containing a certain cuspidal fragment, while the remaining part is exactly as in the nodal case, and we will compute the multiplicities of the special fragments using the above rule for rational trivalent vertices and rational edges of weight $\ell>1$.

\subsection{Cuspidal tropical curves}\label{sec-corr}

\begin{definition}\label{def-ctc}
A plane regular marked tropical curve $(\Gamma,\bp,h)$ of degree $\Delta^*$ and genus $g\ge0$ is called {\it cuspidal} if
it contains one of the following cuspidal fragments:
\begin{enumerate}\item[(A)] a small neighborhood of a four-valent vertex of genus $0$, whose incident edges are mapped by $h$ to four distinct rays,
\item[(B)] a small neighborhood of a flat trivalent vertex of genus $0$,
\item[(C)] an open bounded edge of genus $1$ (i.e., an edge without endpoints),
\item[(D)] a small neighborhood of either a flat trivalent vertex of genus $0$ and a four-valent vertex of genus $0$ connected by two edges of genus $0$, or of a pair of two four-valent vertices of genus $0$ connected by two edges of genus $0$,
\item[(E)] a small neighborhood of a trivalent vertex of genus $1$,
\end{enumerate}
while the remaining part is reduced and trivalent with the zero genus at
all edges and vertices. Furthermore, we equip the cuspidal fragment with the orientation of its edges induced by the regular orientation of $(\Gamma,\bp,h)$.
\end{definition}

\begin{theorem}\label{tuc1}
Let
\begin{itemize}\item $\Delta^*\subset\Z^2\setminus\{0\}$ be a nondegenerate, primitive, balanced multiset
satisfying the hypotheses of Lemma \ref{luc1}, and $\Delta\subset\R^2_\Z$ the induced convex lattice polygon,
\item $0\le g<p_a(\Delta)-1$,
$n=|\partial\Delta\cap\Z^2|+g-2$, \item $(\bn:\widehat C\to C\hookrightarrow\Tor_\K(\Delta),\bw)$ an
$n$-marked parameterized curve such that $C\in V_{\Delta,g}(A_2)$ and $\bx=\val(\bn(\bw))$ is a configuration
of $n$ distinct points in $\R^2$ in general position.\end{itemize} Then the simplified parameterized tropical limit of $C$ (in the sense of Section \ref{sec-ptl})
consists of a regular
cuspidal tropical curve $(\Gamma,\bp,h,\widehat g)$ of degree $\Delta^*$ and genus $g$ such that
$h(\bp)=\bx$, and of the pair $(\bn'_0:(\widehat C^{(0)})'\to C^{(0)},\bw(0))$ such that, for each trivalent vertex $V$ of $\Gamma$ outside the
cuspidal tropical fragment, the corresponding limit curve $(\Tor(P(V)),C^{(0)}_i)$
is rational, nodal, peripherally unibranch and smooth along the toric divisors.
Furthermore, let $z_c\in C$ be the cuspidal singular point, $z_c(t)\in{\mathfrak X}_t$, $t\in(\C,0)$, the
corresponding section of the family ${\mathfrak X}\to(\C,0)$ (see Section \ref{sec-tl}). Then
$z_c(0)\in{\mathfrak X}_0$ belongs to the union of the subvarieties $\Tor(S)\subset{\mathfrak X}$, where
$S$ runs over the cells of the subdivision of $\Delta$ dual to the cells of the reduced cuspidal fragment.
\end{theorem}

{\bf Proof.} We derive the statement from Lemma \ref{ltrop2}. Consider the simplified PTL of $C$.
In general, the corresponding plane marked tropical curve
$(\Gamma,\bp,h,\widehat g)$ has genus $g'\le g$ and degree $\deg(\Gamma,h)$ either equal to $\Delta^*$, or
obtained from $\Delta^*$ by replacing certain disjoint submultisets of $\Delta^*$ by sums of their vectors.
Observe that the case when $g'=g$, $\deg(\Gamma,h)=\Delta$, and $(\Gamma,\bp,h,\widehat g)$ has only bivalent and non-flat trivalent vertices, all of genus zero, and all edges of genus zero, is not possible, since it cannot be the limit of a cuspidal curve (cf., Section \ref{nod-mod}).
Thus, in view of the general position of the configuration
$\bx$, $(\Gamma,\bp,h)$ must be of one of the types listed in Lemma \ref{ltrop2}(i-v):
\begin{itemize}\item In case (i), we obtain either an elliptic non-flat trivalent vertex (cuspidal fragment E), or
an elliptic edge (cuspidal fragment C).
\item In case (ii) we obtain a rational four-valent vertex (cuspidal fragment A).
\item In case (iii) we obtain a flat cycle (cuspidal fragment D).
\item In case (iv) we obtain a flat trivalent vertex incident to two ends of $\Gamma$
(cuspidal fragment B).
\item In case (v) we obtain a flat rational trivalent vertex incident to at least two bounded edges
(cuspidal fragment B).
\end{itemize}

Finally, we notice that, it follows from \cite[Lemmas 3.5 and 3.9]{Shustin2005}
that the rational limit curves, corresponding to rational trivalent vertices of $(\Gamma,\bp,h)$ outside the cuspidal fragment, are purely nodal and that the deformation of the singularities $A_{2\ell-1}$ at the intersection points of limit curves with a common toric divisor, dual to an edge of $(\Gamma,\bp,h,\widehat g)$ outside the cuspidal fragment, also yields only nodal singular points. Hence, the tropicalization of the cuspidal singularity of $C$ lies in the
cuspidal tropical fragment, except maybe for the case of Lemma \ref{ltrop2}(iv).
Below, in Lemma \ref{lmod} we show that this case does not occur.
\proofend

\begin{figure}
\setlength{\unitlength}{1.0mm}
\begin{picture}(150,80)(0,0)
\thicklines

\put(0,65){\line(1,1){10}}\put(0,65){\line(1,-1){10}}
\put(10,55){\line(1,1){10}}\put(10,75){\line(1,-1){10}}
\put(10,55){\line(0,1){20}}\put(35,55){\line(1,2){10}}
\put(35,55){\line(1,0){29}}\put(45,75){\line(1,-2){7}}
\put(52,61){\line(2,-1){12}}

\put(85,55){\line(1,0){10}}\put(85,55){\line(1,2){5}}
\put(95,55){\line(-1,2){5}}\put(115,55){\line(1,0){10}}
\put(115,55){\line(0,1){5}}\put(125,55){\line(0,1){5}}
\put(115,60){\line(1,0){10}}\put(115,60){\line(-1,1){5}}
\put(125,60){\line(-1,1){5}}\put(110,65){\line(1,0){10}}
\put(110,65){\line(1,2){5}}\put(120,65){\line(-1,2){5}}

\put(0,30){\line(1,-2){10}}\put(10,10){\line(1,0){20}}
\put(0,30){\line(3,-2){30}}\put(40,10){\line(1,2){10}}
\put(50,30){\line(1,-1){20}}\put(40,10){\line(1,0){30}}
\put(80,10){\line(1,2){10}}\put(80,10){\line(1,0){40}}
\put(90,30){\line(3,-2){30}}

\thinlines
\put(30,55){\vector(1,0){40}}\put(45,55){\vector(0,1){25}}
\put(80,55){\line(-1,1){5}}\put(80,55){\line(1,0){20}}
\put(100,55){\line(0,1){5}}\put(110,55){\line(-1,1){5}}\put(110,55){\line(1,0){20}}
\put(130,55){\line(0,1){5}}
\dashline{1}(75,60)(75,67)\dashline{1}(100,60)(100,67)
\dashline{1}(105,60)(105,67)\dashline{1}(130,60)(130,67)

\put(0,10){\vector(1,0){35}}\put(0,10){\vector(0,1){25}}
\put(40,10){\vector(1,0){35}}\put(40,10){\vector(0,1){25}}
\put(80,10){\vector(1,0){45}}\put(80,10){\vector(0,1){25}}

\dottedline{1}(50,10)(50,30)\dottedline{1}(40,30)(50,30)
\dottedline{1}(90,10)(90,30)\dottedline{1}(80,30)(90,30)

\put(4,64){$\Delta_i$}\put(12,64){$\Delta_j$}\put(54,70){$P$}
\put(88,68){$T$}\put(108,73){$T$}
\put(9,6){$s$}\put(25,6){$s+2$}\put(2,30){$r$}\put(36,29){$r'$}
\put(49,6){$1$}\put(69,6){$3$}\put(77,29){$2$}\put(89,6){$1$}
\put(114,6){$2k+2$}\put(25,64){$\Longrightarrow$}

\put(26,45){(a)}\put(100,45){(b)}\put(17,0){(c)}\put(55,0){(d)}\put(95,0){(e)}

\end{picture}
\caption{Modifications and the cuspidal fragment of type B}\label{fig-mod}
\end{figure}

\begin{lemma}\label{lmod}
Under the hypotheses of Theorem \ref{tuc1}, a tropicalization of a curve $C\in V_{\Delta,g}(A_2)$
passing through $\bw$ cannot be a cuspidal tropical curve $(\Gamma,\bp,h)$ of degree $\Delta^*$ and genus $g$ that contains a flat trivalent vertex incident to two ends of $\Gamma$.
\end{lemma}

{\bf Proof.}
We argue on the contrary, assuming that some curve $C\in V_{\Delta,g}(A_2)$ tropicalizes to $(\Gamma,\bp,h)$ as described in the lemma.

Since the degree $\Delta^*$ is primitive, the ends incident to the flat trivalent vertex have weight $1$, and hence they are mapped by $h$ onto the same ray in $\R^2$, while the third edge has weight $2$, and thus must be bounded.

Note that the curve $(\Gamma_{red},\bp_{red},h_{red})$ is trivalent of genus $g$, without flat vertices, and
having degree $\deg(\Gamma_{red},h_{red})$ obtained from $\Delta^*$ by replacing two equal primitive vectors $\oa,\oa$ by one vector $2\oa$. Since $n=|\deg(\Gamma_{red},h_{red})|-1+g$, and $h_*(\Gamma)=
\Tr(C)$ passes through $n$ points in general position, by \cite[Proposition 4.3]{M} (or
\cite[Lemma 2.2]{Shustin2005}) the dual subdivision of $\Delta$ consists of triangles and parallelograms.
The boundary of $\Delta$ is subdivided into segments of lattice length $1$ and one segment
$\sigma$ of lattice length $2$, and there is a unique triangle
$T$ (dual to the trivalent vertex of $\Gamma_{red}$ incident to the (unique) double end of
$(\Gamma_{red},h_{red})$) either containing $\sigma$ as a side, or having a parallel,
congruent side joined with $\sigma$ via a sequence of parallelograms (see Figure
\ref{fig-mod}(b)). It follows from Theorem \ref{tuc1} that the limit curves for all triangles different from
$T$ are rational, peripherally unibranch and smooth along the toric divisors, while
the limit curve $C_T\subset\Tor(T)$ for the triangle $T$ is rational and unibranch along the two toric divisors associated with the inclined sides of $T$ (see Figure \ref{fig-mod}(b)). Placing (via a suitable automorphism of $\Z^2$) $T$ on the plane
$\R^2_\Z$ as shown in Figure \ref{fig-mod}(c), we obtain a parametrization of $C_T$ in the form
\begin{equation}x=\alpha t^r,\quad y=\beta t^sR(t),\quad \alpha,\beta\in\C^*,\ \deg R=2\ ,\label{emis2}
\end{equation} and hence $C_T
\subset\Tor(T)$ is immersed outside its two intersection points with the toric divisors
$D_1=\Tor([(s,0),(0,r)])$, $D_2=\Tor([(s+2,0),(0,r)])$. Hence, the tropical limit $z_c(0)$ of the cuspidal singularity belongs to one of the above toric divisors.

If $C_T$ is smooth at the intersection points with the toric divisors $D_1,D_2$, then we obtain the desired contradiction. Indeed, a local modification at $z_c(0)\in D_1$ along $D_1$
would consist of the triangle with the vertices $(-1,0)$, $(1,0)$, and $(0,\gcd(r,s))\}$ and a rational, peripherally smooth, nodal curve in the corresponding toric surface (see
\cite[Lemma 3.9]{Shustin2005}), and none of the points of that curve can be a limit of a cuspidal singularity.

Suppose that $C_T$ is singular at the intersection point $z_c(0)$ with $D_1$.
Consider the function $f(t)=y^r/x^s$ regular in a neighborhood of $z_c(0)$. Since $\beta\ne0$, one has $\ord_t(f(t)-f(0))\le\deg R(t)=2$. Therefore, the multiplicity of the singular point $z_c(0)$ is two, i.e., it is of type $A_{2k}$.
The local modification at the point $z_c(0)$
along $D_1$ results (after a horizontal shift) in the triangle $T_1=\conv\{(0,0),(3,0),(1,r')\}$, where $r'=\gcd(r,s)$ (see Figure \ref{fig-mod}(d)) and a rational limit curve
$C_{T_1}\subset\Tor(T_1)$, unibranch along the toric divisors associated with the inclined sides of
$T_1$. As above one can show that $C_{T_1}$ is immersed outside the intersection
point $z_1\in \Tor([(1,\gcd(r,s)),(3,0)])$. If $\gcd(r,s)$ is odd then the lattice length of the segment $[(1,\gcd(r,s)),(3,0)]$ is $1$; hence
$C_{T_1}$ is smooth at $z_1$, which again yields the desired contradiction. Suppose that
$\gcd(r,s)$ is even. Then $C_{T_1}$ can have at $z_1$ a singularity $A_{2k}$ ($k\ge1$) transversal to the
toric divisor. A local modification at $z_1$ along $\Tor([(1,r'),(3,0)])$ leads to the triangle
$T_2=\conv\{((0,0),(1,2),(2k+2,0)\}$ (see Figure \ref{fig-mod}(e)) and a rational curve $C_{T_2}\subset\Tor(T_2)$. We claim that $C_{T_2}$ is immersed, and hence the desired contradiction. Indeed, its affine part $C_{T_2}\cap\C^2$ is an image of
$\C^*$, which is naturally projected onto the torus $\C^*$ in the horizontal axis, and then, by Riemann-Hurwitz, no ramifications
are possible. Since the edges $[(0,0),(1,2)]$ and $[(1,2),(2k+2,0)]$ have the integral length $1$, the remaining two points of $C_{T_2}$ (corresponding to $0$ and $\infty$ in $\PP^1\setminus\C^*$) are smooth. That is, again we come to a contradiction to the existence of a limit of a cuspidal singularity. \proofend

\subsection{Multiplicities of tropical cuspidal fragments}\label{sec-mult}

Now we compute the multiplicity of each cuspidal tropical curve $(\Gamma,\bp,h,\widehat g)$ as in Theorem \ref{tuc1}.
In each case we show that the multiplicity is the product of some factor $\mu_c$ corresponding to the cuspidal tropical fragment and the factors related to the rational trivalent vertices and rational ordinary and extended edges of $\Gamma$ outside the tropical fragment (the latter factors being as listed in Section \ref{nod-mod}, the last paragraph of part (3). Furthermore, we consider only cuspidal tropical curves with cuspidal fragments of type A, B, and C.
In Section \ref{sec-lp} we show that,
in most interesting cases, tropical curves with cuspidal fragments of type D and E do not appear at all.
The treatment of types D and E can be found in Appendix.

\subsubsection{Multiplicity of a cuspidal fragment of type A}\label{sec-cfa}
Consider a cuspidal fragment $(\Gamma_c,h_c)$ of type A, that is, a small neighborhood of an unmarked four-valent vertex $V$.
Denote by $E_1,E_2$ the edges incident to $V$ and regularly oriented towards $V$ (while the two other edges $E_3,E_4$ are oriented outwards). Denote by $\Delta(\Gamma_c,h_c)$ the multiset of vectors
$\oa_V(E_i)$, $i=1,2,3,4$ and let $Q$ be the quadrangle induced by $\Delta(\Gamma_c,h_c)$. Denote also by $D_i$ the toric divisor
in $\Tor(Q)$ corresponding to $E_i$, $i=1,2,3,4$.

\begin{lemma}\label{lquad1}
Let
$M\subset |{\mathcal L}_Q|$ be the family of peripherally unibranch rational curves $C\in|{\mathcal L}_Q|$ having a singular local branch centered in $\Tor^*(Q)$. Then

(i) If the $h$-images of some two edges of $\Gamma$ lie on the same line, then $M=\emptyset$;
in this case we set
$\mu_c(\Gamma_c,h_c)=0$. Otherwise $M$ is smooth of dimension $2$ and
consists of peripherally smooth curves with a unique singular local branch centered in
the big torus $\Tor^*({\mathcal Q})$; furthermore, this singular branch is of type
$A_2$.

(ii) Let $M\ne\emptyset$, and let $z_1\in D_1$, $z_2\in D_2$ be generic points. Then the intersection of $M$ with the linear system $|{\mathcal L}_Q(-z_1-z_2)|
\subset|{\mathcal L}_Q|$ is transversal.

(iii) The number of intersection points $M\cap|{\mathcal L}_Q(-z_1-z_2)|$ in item (2)
equals $\mu_c(\Gamma_c,h_c)\cdot (w(E_1)w(E_2))^{-1}$.
where
\begin{equation}\mu_c(\Gamma_c,h_c)=|\oa_V(E_1)\wedge\oa_V(E_2)|\ .\label{emuc1}\end{equation}\end{lemma}

{\bf Proof.}
Up to an automorphism of $\Z^2$, we can assume that $Q$ has vertices $(p,0),(q,0),(0,m),(r,s)$, where $m,p,q,r,s>0$,
$q>p\ge m$, $q>r+s$. Let $\bn:\PP^1\to C\in M$. Denoting an affine coordinate on $\C\subset\PP^1$ by $t$, we can suppose that
$$\bn^{-1}\Tor([(p,0),(0,m)])=0,\quad \bn^{-1}\Tor([(p,0),(q,0)])=1,$$
$$\bn^{-1}\Tor([(0,m),(r,s)])=\infty,\quad\bn^{-1}\Tor([(r,s),(q,0)])=\tau\in\C\setminus\{0,1\}\ .$$ Then, in the
corresponding affine coordinates $x,y$ in the torus $\Tor^*(Q)$,
the curve $C$ is parameterized by
\begin{equation}
x=\frac{\alpha t^m}{(t-\tau)^s},\quad y=\frac{\beta t^p(t-1)^{q-p}}{(t-\tau)^{q-r}},
\quad\alpha,\beta\in\C^*\ .
\label{eqquad2}\end{equation}
These formulas define a three-dimensional family, parameterized by $\alpha,\beta,\tau$, which, for a generic $\tau$, yields a curve immersed in
$\Tor^*(Q)$: Indeed,
the singular branch condition in $\Tor^*(Q)$ for some $t=t^*\in\C^*\setminus\{0,1,\tau\}$ reads
$$\frac{dx}{dt}(t^*)=\frac{dy}{dt}(t^*)=0\quad\Longleftrightarrow\quad \begin{cases}&(m-s)t^*-m\tau=0\\
&(qs-m(q-r))t^*-(ps-m(q-r))=0\end{cases}$$ which is a nontrivial condition on
$\tau$. Furthermore, if $Q$ has a pair of parallel edges, i.e., either $m=s$, or
$ps-m(q-r)=0$, one obtains
from the above relations $\tau=0$ contrary to the assumption $\tau\in\C^*\setminus\{0,1\}$, and hence, $M=\emptyset$. Similarly, if $Q$ is degenerate, i.e., the points $(0,m)$, $(r,s)$, and $(q,0)$ are collinear, then $qs-m(q-r)=0$, that is $\tau=\infty$, again a contradiction.
Otherwise, we get uniquely defined values in $\C^*\setminus\{0,1\}$
\begin{equation}t^*=\frac{ps-m(q-r)}{qs-m(q-r)},\quad\tau=\frac{(m-s)(ps-m(q-r))}{m(qs-m(q-r))}
\ .\label{eqquad1}\end{equation} Furthermore, in the latter case,
one can easily check that the intersection points of $C$ with the toric divisors are smooth. For example, the smoothness at the intersection point with the tori divisor $\Tor([(p,0),(q,0)])$ reduces to the relation
$$\frac{dx}{dt}\bigg|_{t=1}=\frac{(m-s)-m\tau}{(1-\tau)^{s+1}}\ne0$$ with the latter inequality coming from
(\ref{eqquad1}). The same argument works for the other toric divisors via an $\Z^2$-automorphism placing another edge of ${\mathcal Q}$ on the horizontal axis.

Thus, $M\simeq(\C^*)^2$.
If $q>1$ is the order of the singular branch $B$ of $C\in M$, then
\cite[Inequality (5)]{IKS}\footnote{This, in fact, is a reformulation of
\cite[Theorem 2]{GuS}.} reads (in the notation of \cite{IKS}, $Z=\emptyset$ and $W$ is the set of
intersection points of $C$ with the toric divisors):
$$-CK_{\Tor(Q)}\ge2+(-CK_{\Tor(Q)}-4)+(q-1)+1\quad\Longrightarrow\quad q\le2\ ,$$
and hence $q=2$ as asserted in the lemma. Let us show that this branch $B$ is of type $A_2$, arguing on the contrary. Suppose that $C=C_{\alpha_0,\beta_0}\in M$ is given by (\ref{eqquad2}) with coefficients
$\alpha_0,\beta_0$, and its singular branch
$B$ is centered at some point $z\in\Tor^*(Q)$ and is of type $A_{2k}$, $k\ge2$. The action of $(\C^*)^2$ on $M$
defines a section of the bundle $\PP T\Tor^*(Q)$ that to each point $z'\in\Tor^*(Q)$ assigns the tangent line
to the singular branch $B'$ centered at $z'$ of the corresponding curve $C'\in M$. Let $(L,z)\subset(\Tor^*(Q),z)$ be the germ of an integral curve to this section, and let
$$\alpha=\varphi(\tau),\quad\beta=\psi(\tau),\quad \tau\in(\C,0),\ \varphi(0)=\alpha_0,\ \psi(0)=\beta_0\ ,$$
be the corresponding curve germ in the space of parameters $\alpha,\beta$. Consider the one-dimensional
family $C_{\alpha,\beta}\subset M$, $\alpha=\varphi(\tau)$, $\beta=\psi(\tau)$, $\tau\in(\C,0)$, and apply to it
\cite[Inequality (5)]{IKS} taking into account that the order of the branch $B$ is $2$ and the
intersection multiplicity of $B$ with $(L,z)$ is at least $4$:
$$-CK_{\Tor(Q)}\ge2+(-CK_{\Tor(Q)}-4+2)+1=-CK_{\Tor(Q)}+1\ ,$$
a contradiction.

Let us fix the intersection points of $C\in M$ with the toric divisors $\Tor(e_1)$ and
$\Tor(e_2)$, where $e_1=[(p,0),(q,0)]$, $e_2=[(p,0),(0,m)]$. This means that the polynomial defining $C$ has truncations
$$(x^{p_1}+ay^{m_1})^d\quad\text{and}\quad x^p\left(\frac{x}{b}+1\right)^{q-p}$$ on the edges $e_1$, $e_2$
of $Q$,
respectively, with given constants $a,b\in\C^*$ and $d=\gcd(m,p)$, $m_1=\frac{m}{d}$, $p_1=\frac{p}{d}$. Plugging the expressions (\ref{eqquad2}) and (\ref{eqquad1}) into these truncations, we obtain
$$m_1(q-p)=\frac{|\vec{e}_1\wedge\vec{e}_2|}
{\|e_1\|_\Z\cdot\|e_2\|_\Z}$$ solutions for the pair
$(\alpha,\beta)\in(\C^*)^2$ in accordance with the assertion of the lemma. Similarly we treat the other choices of pairs of edges of $Q$.
\proofend

\subsubsection{Multiplicity of a cuspidal fragment of type B}\label{sec-cfb}
{\bf(1)} Suppose that $(\Gamma,\bp,h)$ contains a cuspidal fragment
$(\Gamma_c,h_c)$ of type B, that is, a small neighborhood of a flat trivalent vertex $V$. Denote the weights of the incident edges by $m,m_1,m_2$, where $m=m_1+m_2$ (see Figure \ref{ftrap1}(a)). It follows from Lemma
\ref{lmod} that the edge of weight $m$ and at least one of the other edges are bounded
(see Figure \ref{ftrap1}(a)). It will be convenient to extend the cuspidal fragment up to one shown in Figure \ref{ftrap1}(a). We shall compute the multiplicity of the extended fragment and then divide the answer by the product of the factors contributed by the extra trivalent vertices, the result will be the multiplicity of the original cuspidal fragment. Furthermore, we will consider only the lower fragment shown in Figure \ref{ftrap1}(a),
since, in the upper fragment, $m_1=1$, and one can append a trivalent vertex of Mikhalkin's multiplicity $1$
without affecting the multiplicity of the whole extended fragment.
The $h$-image of this fragment is dual in the subdivision of $\Delta$ to the union of
three polygons: triangles $T_1$, $T_2$ dual to the trivalent vertices $h(V_1),h(V_2)$ and a trapeze
$T$ dual to $h(V_3)$ (see Figure \ref{ftrap1}(b)). Let $e_1=T_1\cap T$ and $e_2=T_2\cap T$ be the parallel bases of $T$. The corresponding limit curves $C_i\subset\Tor(T_i)$, $i=1,2$, are rational, nodal,
peripherally unibranch and smooth, while the limit curve $C\subset\Tor(T)$ splits into two irreducible components, $m_1C'$, where $C'\simeq\PP^1$ intersects only the toric divisors $\Tor(e_1),\Tor(e_2)$, and
$C''$ rational, nodal, peripherally unibranch and smooth, disjoint from $\Tor(e_1)$ and
intersecting $\Tor(e_2)$ at the same point $z$ as $C'$.
Note that, in the simplified PTL, it is $m_1$-multiply covered by a component $\widehat C'\simeq\PP^1$ of the parameterizing curve $(\widehat C^{(0})'$ (see Section \ref{sec-ptl}, and it, in fact, corresponds to a bivalent vertex on the edge of weight $m_1$ (see Figure \ref{ftrap1}(a)), but we do not mark it in order to simplify pictures, no confusion will arise.

\begin{figure}
\setlength{\unitlength}{1.0mm}
\begin{picture}(140,185)(0,0)
\thinlines
\put(10,179){\line(1,0){9}}\put(21,179){\line(1,0){9}}
\put(20,181){\line(1,0){10}}\put(30,180){\line(1,0){10}}
\put(20,181){\line(0,-1){6}}\put(20,181){\line(-1,1){5}}
\put(40,180){\line(1,1){5}}\put(40,180){\line(1,-1){5}}
\put(10,164){\line(1,0){9}}\put(21,164){\line(1,0){9}}
\put(20,166){\line(1,0){10}}\put(30,165){\line(1,0){10}}
\put(20,166){\line(0,-1){6}}\put(20,166){\line(-1,1){5}}
\put(40,165){\line(1,1){5}}\put(40,165){\line(1,-1){5}}
\put(10,164){\line(-1,-1){5}}\put(10,164){\line(-1,1){5}}

\put(55,170){\line(1,1){15}}\put(55,170){\line(1,-1){5}}
\put(60,165){\line(0,1){10}}\put(60,165){\line(1,0){10}}
\put(70,165){\line(0,1){20}}\put(70,165){\line(1,1){10}}
\put(80,175){\line(-1,1){10}}

\put(29,178.7){\Large{$\bullet$}}\put(29,163.7){\Large{$\bullet$}}
\put(13,161){$m_1$}\put(13,176){$m_1$}\put(25,168){$m_2$}\put(23,183){$m_2$}
\put(34,167){$m$}\put(34,182){$m$}
\put(29,160){$V$}\put(29,175){$V$}
\put(10,166){$V_1$}\put(37,161){$V_2$}\put(19,168){$V_3$}
\put(56,161){$T_1$}\put(64,161){$T$}\put(73,161){$T_2$}
\put(25,150){\rm{(a)}}\put(65,150){\rm{(b)}}

\put(5,110){\vector(0,1){35}}\put(5,110){\vector(1,0){35}}
\put(50,110){\vector(0,1){35}}
\put(50,110){\vector(1,0){35}}
\put(50,110){\line(2,3){20}}\put(60,125){\line(2,-3){10}}
\put(70,140){\line(1,-3){10}}
\put(5,110){\line(1,1){10}}
\put(15,120){\line(1,2){10}}
\put(25,140){\line(1,-3){10}}
\dottedline{1}(60,110)(60,125)
\dottedline{1}(70,110)(70,140)
\dottedline{1}(105,110)(105,120)
\dottedline{1}(25,110)(25,140)
\dottedline{1}(50,125)(60,125)
\dottedline{1}(50,140)(70,140)
\dottedline{1}(5,120)(15,120)
\dottedline{1}(5,140)(25,140)
\dottedline{1}(15,110)(15,120)

\put(1,139){$m$}\put(46,139){$m$}\put(91,139){$m$}
\put(0,119){$m_1$}\put(45,124){$m_1$}
\put(14,106){$1$}\put(24,106){$2$}\put(34,106){$3$}\put(59,106){$1$}\put(69,106){$2$}\put(79,106){$3$}
\put(20,100){\rm{(d)}}\put(65,100){\rm{(e)}}

\put(90,119){$m_1$}
\put(110,100){\rm{(f)}}

\put(90,160){\vector(0,1){35}}\put(90,160){\vector(1,0){35}}
\put(90,160){\line(2,3){20}}
\put(110,190){\line(1,-3){10}}
\dottedline{1}(110,160)(110,190)\dottedline{1}(90,175)(100,175)\dottedline{1}(90,190)(110,190)
\dottedline{1}(100,160)(100,175)
\put(99,156){$1$}\put(109,156){$2$}\put(119,156){$3$}
\put(86,189){$m$}\put(85,174){$m_1$}\put(105,150){\rm{(c)}}

\put(95,110){\vector(0,1){35}}\put(95,110){\vector(1,0){35}}
\put(95,110){\line(1,1){10}}\put(105,120){\line(1,2){10}}
\put(115,140){\line(1,-3){10}}\put(105,120){\line(1,-1){10}}
\dottedline{1}(95,120)(105,120)\dottedline{1}(115,110)(115,140)\dottedline{1}(105,110)(105,120)
\dottedline{1}(95,140)(115,140)\put(114,106){$2$}\put(104,106){$1$}\put(124,106){$3$}

\thicklines

\put(20,74){\line(1,0){9}}\put(31,74){\vector(1,0){9}}
\put(40,76){\vector(-1,0){10}}\put(40,75){\vector(1,0){10}}
\put(30,76){\vector(0,-1){6}}\put(25,81){\vector(1,-1){5}}
\put(50,75){\vector(1,-1){5}}\put(55,80){\vector(-1,-1){5}}
\put(15,79){\vector(1,-1){5}}\put(15,69){\vector(1,1){5}}
\put(28,59){\vector(-1,0){8}}\put(40,59){\line(-1,0){8}}
\put(30,61){\vector(1,0){10}}\put(40,60){\vector(1,0){10}}
\put(30,55){\vector(0,1){6}}\put(25,66){\vector(1,-1){5}}
\put(55,65){\vector(-1,-1){5}}\put(50,60){\vector(1,-1){5}}
\put(20,59){\vector(-1,-1){5}}\put(15,64){\vector(1,-1){5}}
\put(29,44){\vector(-1,0){9}}\put(40,44){\line(-1,0){9}}
\put(40,46){\vector(-1,0){10}}\put(50,45){\vector(-1,0){9}}
\put(30,46){\vector(0,-1){6}}\put(25,51){\vector(1,-1){5}}
\put(55,50){\vector(-1,-1){5}}\put(55,40){\vector(-1,1){5}}
\put(20,44){\vector(-1,-1){5}}\put(15,49){\vector(1,-1){5}}

\put(39,43.7){\Large{$\bullet$}}\put(39,73.7){\Large{$\bullet$}}
\put(39,58.7){\Large{$\bullet$}}
\put(34,30){\rm{(g)}}\put(76,30){\rm{(h)}}\put(107,30){\rm{(i)}}

\put(80,55){\vector(0,-1){10}}
\put(80,55){\vector(-1,-1){10}}\put(65,50){\vector(1,-1){5}}\put(70,45){\vector(0,-1){5}}
\put(70,65){\vector(1,-1){10}}\put(90,58.3){\vector(-3,-1){10}}

\put(110,47){\vector(0,-1){10}}
\put(110,47){\vector(-1,-1){7}}
\put(103,54){\vector(1,-1){7}}\put(120,50.3){\vector(-3,-1){10}}
\put(110,67){\vector(0,-1){10}}
\put(110,67){\vector(-1,-1){7}}
\put(103,74){\vector(1,-1){7}}\put(120,70.3){\vector(-3,-1){10}}
\put(110,87){\vector(0,-1){10}}
\put(110,87){\vector(-1,-1){7}}
\put(103,94){\vector(1,-1){7}}\put(120,90.3){\vector(-3,-1){10}}
\put(106,43){\Large{$*$}}\put(105,69.5){\Large{$*$}}\put(114,87.5){\Large{$*$}}

\put(19,14){\vector(-1,0){9}}\put(21,14){\vector(1,0){9}}
\put(30,16){\vector(-1,0){10}}\put(30,15){\vector(1,0){10}}
\put(20,16){\vector(0,-1){6}}\put(15,21){\vector(1,-1){5}}
\put(45,20){\vector(-1,-1){5}}\put(40,15){\vector(1,-1){5}}
\put(10,14){\vector(-1,-1){5}}\put(5,19){\vector(1,-1){5}}
\put(23,12.8){\Large{$*$}}
\put(29,13.7){\Large{$\bullet$}}\put(29,0){\rm{(j)}}
\put(64,14){\vector(-1,0){9}}\put(75,14){\line(-1,0){9}}
\put(65,16){\vector(1,0){10}}\put(75,16){\vector(-1,0){10}}
\put(75,15){\vector(1,0){10}}
\put(65,16){\vector(0,-1){6}}\put(60,21){\vector(1,-1){5}}
\put(90,20){\vector(-1,-1){5}}\put(85,15){\vector(1,-1){5}}
\put(55,14){\vector(-1,-1){5}}\put(50,19){\vector(1,-1){5}}
\put(68,14.8){\Large{$*$}}
\put(74,13.7){\Large{$\bullet$}}\put(74,0){\rm{(k)}}
\put(109,14){\vector(-1,0){9}}\put(120,14){\line(-1,0){9}}
\put(120,16){\vector(-1,0){10}}\put(130,15){\vector(-1,0){9}}\put(120,15){\vector(1,0){10}}
\put(110,16){\vector(0,-1){6}}\put(105,21){\vector(1,-1){5}}
\put(135,20){\vector(-1,-1){5}}\put(130,15){\vector(1,-1){5}}
\put(100,14){\vector(-1,-1){5}}\put(95,19){\vector(1,-1){5}}
\put(125,13.8){\Large{$*$}}
\put(119,13.7){\Large{$\bullet$}}\put(119,0){\rm{(l)}}

\end{picture}
\caption{Cuspidal flat vertex and its modification}\label{ftrap1}
\end{figure}

Notice that the lengths of the horizontal edges in the fragment depicted in Figure \ref{ftrap1}(a) do not affect the
multiplicity we are looking for (specific rational values of the multiplicities cannot jump in a continuous deformation of these lengths).
Hence, we can restrict our attention to the cases $m_1=m_2$ and $m_1<m_2$.

Assuming that the convex piecewise linear function $\nu:\Delta\to\R$ vanishes along the trapeze $T$, we get that the truncations of the polynomial $F(x,y)$ to the edges $e_1=T_1\cap T$ and $e_2=T\cap T_2$ are
$$x^iy^j(y-\alpha)^{m_1}+O(t^{>0})\quad\text{and}\quad x^ky^l((y-\alpha)^m+O(t^{>0}),$$
respectively. Thus, performing the coordinate change (\ref{e-new61}), we obtain The Laurent polynomial $F^{(1)}(x,y):=F(x_1,y+\alpha+\zeta(t))$ which indices a subdivision of its Newton polygon with the (somehow subdivided) domain $\widetilde P$ of one of the shapes shown in Figure \ref{ftrap1}(c,d) according as $m_1=m_2$ or $m_1<m_2$. The incline sides of $\widetilde P$ are the Newton diagrams of the tangency points of the limit curves in $\Tor(T_1),\Tor(T_2),\Tor(T)$ with the toric divisors $\Tor(e_1),\Tor(e_2)$.
It is easy to verify that the induced piecewise linear function
$\widetilde\nu$ on $\widetilde P$ is such that its graph restricted to the union of the edges
$$[(0,0),(1,m_1)],\quad [(1,m_1),(2,m)],\quad [(2,m),(3,0)]$$ of $\widetilde P$ does not lie in one plane.
Hence, $\widetilde P$ is necessarily subdivided, and the point $(1,m_1)$ is a vertex of the subdivision.

The local tropical limit
associated with $\widetilde P$
meets the following requirements: the union
$C^{(0)}_{\widetilde P}$ of the limit curves has arithmetic genus zero and deforms into a rational curve with
three holes, with a cusp and (possibly) some nodes in the deformation along the family $\{C^{(t)}\}_{t\in(\C,0)}$. It follows that
\begin{itemize}\item an interior integral point of the polygon $\widetilde P$ cannot be a vertex of the subdivision, since otherwise the union of the limit curves would be of a positive arithmetic genus,
\item a possible subdivision of $\widetilde P$ indices the union of toric surfaces, whose common toric divisors intersect with the limit curves only at one point, again since otherwise the union of the limit curves would have a positive arithmetic genus,
    \item the subdivision cannot consist of only triangles or of triangles and a trapeze, since, in such a case, no cuspidal singularity appears in the deformation, see Lemma \ref{lquad1}(i) and \cite[Proposition
6.17]{M} (or \cite[Lemma 3.5]{Shustin2005}); for example, the subdivision into two triangles with the common edge $[(1,m_1),(3,0)]$ is not possible, since the rational curve with the lower Newton triangle cannot have a cuspidal singularity as explained in the proof of Lemma \ref{lmod} (see Figure \ref{fig-mod}(d) and the related argument).
\end{itemize}
Furthermore, we can chose $\zeta(t)$ in the shift formula (\ref{e-new61}) in order to eliminate the monomial $xy^{m_1-1}$ from $F^{(1)}(x,y)$. This, in particular yields that the subdivision of $\widetilde P$ cannon contain the edge $[(1,0),(1,m_1)]$, since otherwise the limit curves with cross the toric divisor $\Tor([(1,0),(1,m_1)])$ in at least two points, and hence the union of all limit curves related to the domain $\widetilde P$ will have a positive arithmetic genus in contradiction to the rationality of the considered tropical cuspidal fragment.

The latter observation yields that the only possible subdivisions are shown in Figure \ref{ftrap1}(e,f).

{\bf(2)} We are now in a position to determine all possible multiplicities of the considered fragment.
Namely, the multiplicity we are looking for counts in how many ways one can find a tuple of limit curves corresponding to the trivalent vertices $V_1,V_2,V_3$ and to a pair of a trivalent and a four-valent vertex of the modification.
The answer depends on the induced regular orientation of the fragment without marked points as shown in
Figure \ref{ftrap1}(g), or of the fragment with a marked point as shown in Figure \ref{ftrap1}(j,k,l), as well as of the fragment shown in Figure \ref{ftrap1}(h) dual to a suitable subdivision in the modification shown in Figures \ref{ftrap1}(e,f).

Suppose that no horizontal edge contains a marked point. Then possible orientations are shown in
Figure \ref{ftrap1}(g) (here orientations of non-horizontal edges may be chosen other way,
provided that precisely two edges are oriented towards a non-flat trivalent vertex).
Denote by $m'_1,m'_2,m'_3,m'_4$ the weights of the non-horizontal edges oriented inside the fragment.
The induced orientation on the fragment dual to the subdivision of $\widetilde P$ in all the cases is as shown
in Figure \ref{ftrap1}(h). Recall that the limit curves of the tropical limit can be restored in a recursive procedure
obeying the above orientation of edges in the
following way (cf. \cite[Section 3.7]{Shustin2005}). The edges merging to a tri- or four-valent vertex indicate the toric divisors in the toric surface,
associated with the dual triangle or quadrangle, where the intersection points with the considered limit curve are already fixed by the
previously constructed limit curves. Now, applying \cite[Lemmas 3.5 and 3.9]{Shustin2005} (cf., the factors described in the last paragraph of Section \ref{nod-mod}, part (2)) and Lemma \ref{lquad1},
we obtain the multiplicities
\begin{equation}\frac{\mu(V_1)\mu(V_2)\mu(V_3)}{m'_1m'_2m'_3m'_4}\cdot\mu_c(\Gamma_c,h_c)\ ,\label{e-new73}\end{equation}
where $\mu_c(\Gamma_c,h_c)$ equals
\begin{equation}\frac{(m+m_2)m_1}{mm_2},\quad \frac{m+m_2}{m},\quad \frac{m+m_2}{m_2}
\label{emuc}\end{equation}
for the upper, middle, and lower oriented fragments shown in Figure \ref{ftrap1}(g),
respectively. For example, the first value in (\ref{emuc}) is obtained as follows:
\begin{itemize}\item the multiplicative contributions of the trivalent vertices $V_1,V_2,V_3$ in the upper fragment in Figure \ref{ftrap1}(g) are $\mu(V_1)(m'_1m'_2)^{-1}$, $\mu(V_2)(m'_3m)^{-1}$, and $\mu(V_3)(m'_4m_2)^{-1}$, respectively,
\item the multiplicative contribution of the four-valent vertex in the modification shown in Figure \ref{ftrap1}(h)
is $|(m_2,-1)\wedge(-m,-1)|=m+m_2$ (cf. Lemma \ref{lquad1}(iii)), while the contribution of the trivalent vertex is $m_1$ (cf. \cite[Lemma 3.9]{Shustin2005}).
\end{itemize}
Thus, we obtain the equality (cf., (\ref{e-new73}))
$$\frac{\mu(V_1)\mu(V_2)\mu(V_3)}{m'_1m'_2m'_3m'_4}\cdot\mu_c(\Gamma_c,h_c)=
\frac{\mu(V_1)}{m'_1m'_2}\cdot\frac{\mu(V_2)}{m'_3m}
\cdot\frac{\mu(V_3)}{m'_4m_2}\cdot(m+m_2)\cdot m_1$$
$$\Longrightarrow\quad\mu_c(\Gamma_c,h_c)=\frac{(m+m_2)m_1}{mm_2}.$$

Suppose now that the fragment contains a marked point on a horizontal edge of the cuspidal fragment
$(\Gamma_c,h_c)$ (shown by asterisk in Figures
\ref{ftrap1}(j,k,l)). Here, three of the non-horizontal edges are oriented inward,
and we denote their weights by $m'_1,m'_2,m'_3$.
As compared to the preceding case, the computation of the multiplicity requires also
the modification
of the marked point condition (cf. \cite[Section 2.5.9]{IMS}). We comment on it in more detail.
Namely, this modification requires a shift
\begin{equation}(x,y)\mapsto(x,y+b+\zeta(t)),\quad b\in\C,\ \zeta(0)=0\ ,\label{emn13}\end{equation}
that annulates the monomial $x^2y^{m-1}$ in the polynomial $\widehat F(x,y):=F(x,y+b+\zeta)$, and the determination of $\zeta(t)$ depends on the algebraic preimage $w=(\xi(t),b)$
The shift (\ref{emn13}) moves $w$ to the point
$\widehat w=(\xi(t),
-\zeta(t))$, whose tropical image appears on one of the edges of the four-valent vertex
of the fragment shown in Figure \ref{ftrap1}(h): in accordance with the original
position shown in Figures \ref{ftrap1}(j,k,l), the new position is shown by asterisk in Figure \ref{ftrap1}(i).
The condition to pass through $\widehat w$ yields
\begin{equation}\widehat F(\xi(t),-\zeta(t))=0
\ .\label{emn15}\end{equation} The minimal exponent of $t$ in the left-hand side of (\ref{emn15})
must occur twice, and it happens at the two vertices of one of the inclined sides of the quadrangle
(see Figure \ref{ftrap1}(d)), namely, for the side dual to the
edge containing the marked point (see Figure \ref{ftrap1}(i)). Let
$$a_{1,m_1}t^\alpha+\text{h.o.t},\quad a_{2,0}t^\beta+\text{h.o.t.},
\quad a_{2,m}+\text{h.o.t.},\quad a_{3,0}t^\gamma+\text{h.o.t.}$$
be the coefficients of $\widehat F$ at the vertices of $Q$, where $a_{1,m_1},a_{2,0},a_{2,m},a_{3,0}\in\C^*$
are determined by the limit curves associated with the subdivision shown in Figure \ref{ftrap1}(d).
Relation (\ref{emn15}) yields one of the equations
\begin{equation}\begin{cases}&a_{1,m_1}t^\alpha(-\zeta(t))^{m_1}+a_{2,0}
t^\beta\xi(t)+\text{h.o.t.}=0,\\
&a_{1,m_1}t^\alpha+a_{2,m}\xi(t)(-\zeta(t))^{m_2}+\text{h.o.t.}=0,\\
&a_{2,m}(-\zeta(t))^m+a_{3,0}t^\gamma\xi(t)+\text{h.o.t.}=0\end{cases}\label{efixed}\end{equation}
in accordance with the case shown in Figures \ref{ftrap1}(j,k,l), respectively. Correspondingly, we
obtain $m_1$, $m_2$, or $m$ solutions for $\Ini(\zeta)$. Combining this result with the count of the limit
curves related to the subdivision shown in Figure \ref{ftrap1}(d), we finally obtain the multiplicity
$$\frac{\mu(V_1)\mu(V_2)\mu(V_3)}{m'_1m'_2m'_3}\cdot\mu_c(\Gamma_c,h_c)\ ,$$
where $\mu_c(\Gamma_c,h_c)$ is again given by (\ref{emuc})
in accordance with the cases shown in Figures \ref{ftrap1}(j,k,l), respectively.

\subsubsection{Multiplicity of a cuspidal fragment of type C}\label{sec-cfc}
Let $(\Gamma,\bp,h,\widehat g)$ contain a cuspidal fragment $(\Gamma_c,\bp_c,h_c)$ of type C, that is, an open bounded edge of genus $1$.
We extend $(\Gamma_c,\bp_c,h_c)$ to a fragment $(\Gamma',\bp',h')$
by adding the endpoints
$V_1,V_2\in\Gamma^0$ of this edge and other edges of $\Gamma$ incident to $V_1,V_2$ (see
Figure \ref{fig-cfc}(a)). As in the preceding section, we compute the multiplicity of the extended fragment and then divide by the multiplicative contribution of the two added vertices.

By Theorem \ref{tuc1}, the vertices $V_1,V_2$ are trivalent, and we can suppose that the
triangles $T_1=P(V_1)$ and $T_2=P(V_2)$ of the subdivision of the Newton polygon $P$ of $(\Gamma,\bp,h)$ share a common side $S$ dual to the edge $e=h(\Gamma_c)$,
see Figure \ref{fig-cfc}(b). By Theorem \ref{tuc1}, the limit curves $C_1\subset\Tor(T_1)$, $C_2\subset\Tor(T_2)$ are
nodal rational, smooth and unibranch along the toric divisor $\Tor(e)$. Thus, performing a modification
as in Section \ref{sec-cfb}, we obtain a new fragment dual to the lattice triangle
$T=\conv\{(0,0),(0,2),(m,1)\}$ (cf. \cite[Section 3.5]{Shustin2005} and \cite[Section 2.5.8]{IMS}), and the limit curve $C\subset\Tor(T)$ must have a singularity more complicated than a node (a tropical limit of the cusp) in the torus $\Tor^*(T)$ and
genus $\le1$. By \cite[Lemma 3.9]{Shustin2005}, a rational curve with Newton
triangle $T$ must be nodal, and hence, $C$ is elliptic. Furthermore, $C$ cannot be immersed, since a singularity of an immersed curve may deform into a cuspidal curve 
only with a jump of the genus: indeed, the conservation of genus assumes that the deformation of each singularity is $\delta$-constant, and hence equinormalizable according to \cite[Theorem 1, page 73]{Tei} and \cite[Proposition 3.3]{ChL} (see also \cite[Theorem II.2.56]{GLS}), but then all local branches stay immersed in a deformation, that is, no cusp pops up.

\begin{figure}
\setlength{\unitlength}{1.0mm}
\begin{picture}(130,80)(0,0)
\thicklines
\put(10,50){\line(1,1){10}}\put(10,70){\line(1,-1){10}}
\put(20,60){\line(1,0){20}}\put(40,60){\line(1,1){10}}
\put(40,60){\line(1,-1){10}}\put(60,60){\line(1,1){10}}
\put(60,60){\line(1,-1){10}}\put(70,50){\line(0,1){20}}
\put(70,50){\line(1,1){10}}\put(70,70){\line(1,-1){10}}
\put(90,50){\line(1,0){20}}\put(90,50){\line(1,2){10}}
\put(110,50){\line(-1,2){10}}

\put(20,20){\vector(-1,-1){10}}\put(10,30){\vector(1,-1){10}}
\put(45,20){\vector(-1,0){25}}\put(55,10){\vector(-1,1){10}}
\put(55,30){\vector(-1,-1){10}}

\put(80,20){\vector(-1,-1){10}}\put(70,30){\vector(1,-1){10}}
\put(95,20){\vector(-1,0){15}}\put(95,20){\vector(1,0){10}}\put(105,20){\vector(1,-1){10}}
\put(115,30){\vector(-1,-1){10}}

\thinlines
\put(90,50){\vector(0,1){30}}\put(90,50){\vector(1,0){30}}
\dottedline{1}(100,50)(100,70)\dottedline{1}(100,70)(90,70)

\put(95,18.7){\Large{$*$}}\put(85,69){$m$}\put(109,46){$2$}\put(89,46){$0$}
\put(64,59){$T_1$}\put(72,59){$T_2$}\put(20,22){$V_1$}\put(41,22){$V_2$}
\put(80,22){$V_1$}\put(101,22){$V_2$}\put(95,16){$p$}

\put(30,0){(d)}\put(90,0){(e)}\put(30,40){(a)}\put(68,40){(b)}\put(98,40){(c)}

\end{picture}
\caption{Cuspidal elliptic edge and its modification}\label{fig-cfc}
\end{figure}

\begin{lemma}\label{ledge1}
Let $m\ge3$, and $T$ be the triangle with vertices $(0,0)$, $(0,2)$, and $(m,1)$.
\begin{enumerate}
\item[(i)] If an elliptic curve $C\subset\Tor(T)$, defined by a polynomial 
with the Newton polygon $T$, has at least one singular local branch, then such a singular local branch is unique, and it is of order $2$.
\item[(ii)] Suppose that the curve $C$ from item (i) has one ordinary cusp and $m-3$ nodes. Then the equisingular family $M\subset|{\mathcal L}_T|$, containing $C$, is the smooth $3$-dimensional orbit of the action of the group of transformations $G:=\{(x,y)\mapsto(\alpha x+a,\beta y),\ \alpha,\beta\in\C^*,\ a\in\C\}$.
\item[(iii)] The family $M$ has a finite transversal intersection with the subspace \mbox{$\Lambda(z_1,z_2)\subset|{\mathcal L}_T|$} defined by the conditions to pass through fixed points $z_1\in\Tor([(0,2),(m,1)])$, $z_2\in\Tor([(0,0),(m,1)])$, and by the vanishing of the coefficient of $x^{m-1}y$ in the defining polynomial.
\end{enumerate}
\end{lemma}

{\bf Proof.} Note, first, that whatever the singularities of $C$ are, the equisingular family containing $C$ contains also the corresponding $3$-dimensional orbit $M$ of the $G$-action.
Next, we note that the intersection multiplicity of $C$ with any line $x=\const$ is $2$; hence, the singularities of $C\in M$ are of types $A_k$, $k\ge1$. Suppose that $C$ has two singularities $A_{2i},A_{2j}$ and some other singularities. 
We claim that in such a case, the dimension of the germ
at $C$ of the equisingular family does not exceed
\begin{equation}
\dim|{\mathcal L}_T|-\delta(C)-2=m+2-(m-2)-2=2\ ,\label{eqedge}\end{equation} where $\delta(C)$ is the total $\delta$-invariant of $C$. Indeed, following
the proof of \cite[Theorem 1.1]{Shustin1998}, 
we choose projections of the germ of $|{\mathcal L}_T|$ at $C$ to versal
deformation bases of the singular points of $C$ and consider the germ at $C$ of the family $M_1\subset
|{\mathcal L}_T|$, which is the intersection of pull-backs of the equiclassical
strata in the versal deformations of $A_{2i}$ and $A_{2j}$, and the equigeneric strata in the versal deformations of the rest of singularities (see the definitions and main properties of the equiclassical and equigeneric strata in \cite[Sections 1, 4, and 5]{DH} and in \cite{Di}). The germ at $C$ of the equisingular stratum is contained in
$M_1$. According to \cite[Section 4.3]{DH} and \cite[Section 2]{Shustin1998} the tangent cone to $M_1$
can be identified with $H^0(\widehat C,{\mathcal O}_{\widehat C}(D-D_0-D_{ec}))$, where $\widehat C\to C$ is the
normalization, $D$ is the pull-back of $c_1({\mathcal O}_C\otimes{\mathcal L}_T)$, $D_0$ is the double-point divisor, $D_{ec}$ is the pull-back of the singular points of type $A_{2i}$ and $A_{2j}$. Here $g(\widehat C)=1$, and
$$\deg(D-D_0-D_{ec})=C^2_{\Tor(T)}-2\delta(C)-2=2m-2(m-2)-2=2$$ $$>2g(\widehat C)-2=0\ ,$$ and hence,
(\ref{eqedge}) follows from Riemann-Roch. However, then the equisingular family cannot contain the $3$-dimensional orbit $M$, which is a contradiction proving claim (i).


Under the assumption of claim (ii), the equisingular and the equiclassical families coincide, and according to \cite[Theorem 6.1(ii,iii) and Example 6.4(1)]{GK}, the germ at $C$ of the equisingular family is smooth and $3$-dimensional; hence, it coincides with $M$. Moreover, the germ of $M$ at $C$ is a transversal intersection in $|{\mathcal L}_T|$ of the germs of the (smooth) equisingular families corresponding to the individual nodes and to the cusp of $C$. The transversality of the intersection of $M$ with the space $\Lambda(z_1,z_2)$ immediately follows
from the fact that $M$ is an orbit of the $G$-action. Indeed, the intersection comes is described by the system of equations
$$\frac{\alpha^m\beta a^0_{m1}}{a^0_{00}}=\xi_0,\quad\frac{\alpha^m\beta a^0_{m1}}{\beta^2a^0_{02}}=\eta_0,\quad \alpha^{m-1}\beta a^0_{m-1,1}+ma\alpha^m\beta a^0_{m1}=0,$$
in the parameters $\alpha,\beta,a$ of the group $G$ with the fixed nonzero coefficients $a^0_{00},a^0_{m1},a^0_{m-1,1},a^0_{02}$ and $\xi_0,\eta_0$ determined by the given curve $C$ and the points $z_1,z_2$, respectively. Clearly, all solutions to the system are simple. \proofend

%

\begin{definition}\label{dcfc}
For a pair of integers $1\le p\le q$, introduce the convex lattice polygons
$$\sigma_{p,q}(n)=\conv\left\{(i,j)\in\Z^2\ :\ i,j\ge0,\ pi+qj\le n\right\},\quad n\ge1\ .$$
For a pair of convex polygons $P_1,P_2\subset\R^2_\Z$, denote by $\langle P_1,P_2\rangle$ their mixed area.
Put $\theta(1)=\theta(2)=0$, $\theta(3)=2$ and, for $k\ge2$,
$$\theta(2k+1)=2(2k+1)(\langle\sigma_{2,3}(k),\sigma_{2,3}(k+1)\rangle+
\langle\sigma_{1,2}(k-1),\sigma_{1,2}(k)\rangle)\ ,$$
$$\theta(2k)=4k(\langle\sigma_{2,3}(k-2),\sigma_{2,3}(k-1)\rangle+
\langle\sigma_{1,2}(k-1),\sigma_{1,2}(k)\rangle)\ .$$
\end{definition}

\begin{lemma}\label{ledge2}
Let $m\ge1$, and let $T$ be the triangle with vertices $(0,0)$, $(0,2)$, and $(m,1)$.
There are exactly
$\theta(m)$
polynomials $F(x,y)=y^2-2yf(x)+1$ with Newton polygon $T$ and such that
\begin{itemize}\item $f(x)$ is a monic polynomial of degree $m$ with the vanishing coefficient of $x^{m-1}$,
\item the curve defined by $F$ in the toric surface $\Tor(T)$ is elliptic and has a singular local branch.
\end{itemize} Moreover, each of the above curves has one ordinary cusp and $(m-3)$ nodes.
\end{lemma}

{\bf Proof.}
The equation $F(x,y)=0$ yields $y=f(x)\pm\sqrt{f(x)^2-1}$, and hence a singular branch of $C=\{F=0\}$
corresponds to a root of $f(x)^2-1$ of multiplicity $2s+1$, $1\le s<m/2$, while other singularities correspond
to roots of even multiplicity. Hence,
$f(x)^2-1=(x-\xi)^{2s+1}Q(x)^2R_3(x)$, where $\xi\in\C$, $\deg Q=m-s-2$, $\deg R_3=3$, where
$Q(\xi)R_3(\xi)\ne0$. Note that $R_3(x)$
has no multiple roots, since otherwise the curve $C$ would be rational, but a rational curve with Newton polygon
$T$ must be nodal \cite[Lemma 3.5]{Shustin2005}. Observe also that $\gcd(Q,f)=1$.
We will canonically normalize $F(x,y)$ by substitution of $x+\xi$ for $x$.

There are no solutions for $m\le2$. If $m=3$, we obtain
two solutions $f(x)=x^3\pm1$.

Suppose that $m>3$. Since $\gcd(f-1,f+1)=1$, we have for $m=2k+1$, $k\ge2$,
\begin{equation}\text{either}\quad \begin{cases}f(x)-1=x^{2s+1}S_{k-s}(x)^2&\\
f(x)+1=T_{k-1}(x)^2R_3(x)&\end{cases}\label{eqedge1}\end{equation}
\begin{equation}\text{or}\quad \begin{cases}
f(x)+1=x^{2s+1}S_{k-s}(x)^2&\\
f(x)-1=T_{k-1}(x)^2R_3(x)&\end{cases}\label{eqedge2}\end{equation}
\begin{equation}\text{or}\ \begin{cases}f(x)-1=x^{2s+1}S_{k-s-1}(x)^2R_2(x)&\\
f(x)+1=T_k(x)^2R_1(x)&\end{cases}\label{eqedge3}\end{equation}
\begin{equation}\text{or}\quad \begin{cases}
f(x)+1=x^{2s+1}S_{k-s-1}(x)^2R_2(x)&\\
f(x)-1=T_k(x)^2R_1(x)&\end{cases}\label{eqedge4}\end{equation}
where $R,S,T$ stand for monic polynomials with the subindex designating the degree, pairwise coprime in each system
(\ref{eqedge1})-(\ref{eqedge4}). Similarly for $m=2k$, $k\ge2$, we get
\begin{equation}\text{either}\quad \begin{cases}f(x)-1=x^{2s+1}S_{k-s-1}(x)^2R_1(x)&\\
f(x)+1=T_{k-1}(x)^2R_2(x)&\end{cases}\label{eqedge5}\end{equation}
\begin{equation}\text{or}\quad \begin{cases}
f(x)+1=x^{2s+1}S_{k-s-1}(x)^2R_1(x)&\\
f(x)-1=T_{k-1}(x)^2R_2(x)&\end{cases}\label{eqedge6}\end{equation}
\begin{equation}\text{or}\ \begin{cases}f(x)-1=x^{2s+1}S_{k-s-2}(x)^2R_3(x)&\\
f(x)+1=T_k(x)^2&\end{cases}\label{eqedge7}\end{equation}
\begin{equation}\text{or}\quad \begin{cases}
f(x)+1=x^{2s+1}S_{k-s-2}(x)^2R_3(x)&\\
f(x)-1=T_k(x)^2&\end{cases}\label{eqedge8}\end{equation}
where $R,S,T$ stand for monic polynomials with the subindex designating the degree, pairwise coprime in each of the systems
(\ref{eqedge5})-(\ref{eqedge8}).

Let us analyze system (\ref{eqedge1}). Differentiating both the equations, we obtain
\begin{equation}x^{2s}S_{k-s}((2s+1)S_{k-s}+2xS'_{k-s})=T_{k-1}(2T'_{k-1}R_3+T_{k-1}R'_3)\ .\label{erev2_1}\end{equation}
Taking into account that the factors in each side of this equation are pairwise coprime, we derive that $s=1$ and that none of $S_{k-1}$ and $T_{k-1}$ has multiple roots; hence, all solutions represent curves with an ordinary cusp and $(m-3)$ nodes. Thus, by Lemma \ref{ledge1}(ii,iii), system (\ref{eqedge1}) describes a (zero-dimensional) subscheme of the intersection $M\cap\Lambda(z_1,z_2)$, which then is finite and consists of simple points (i.e., zero-dimensional schemes of length $1$). That is, all the solutions to system (\ref{eqedge1}) which we recover below are isolated and simple.

It follows from (\ref{erev2_1}) and the relation $s=1$ that
\begin{equation}
\begin{cases}&3S_{k-1}+2xS'_{k-1}=(2k+1)T_{k-1}\\
&2T'_{k-1}R_3+T_{k-1}R'_3=(2k+1)x^2S_{k-1}\end{cases}\label{eqedge9}\end{equation}
Plugging $S_{k-1}=x^{k-1}+\alpha_1x^{k-2}+...+\alpha_{k-1}$, $T_{k-1}=x^{k-1}+
\beta_1x^{k-2}+...+\beta_{k-1}$, we obtain from the former equation in (\ref{eqedge9})
\begin{equation}\beta_i=\frac{2k+1-2i}{2k+1}\alpha_i,\quad i=1,...,k-1\ ,\label{erev2_2}\end{equation} which we
together with $R_3=x^3+\gamma_1x^2+\gamma_2x+\gamma_3$ plug to the second equation and,
subsequently equating the coefficients of $x^{k-1},...,x^2,x,1$, obtain the following system of
quasihomogeneous equations, where each variable $\alpha_i,\gamma_i$ has weight $i$, $i\ge1$:
\begin{equation}\frac{4i(2k+1-i)}{2k+1}\alpha_i=\sum_{\renewcommand{\arraystretch}{0.6}
\begin{array}{c}
\scriptstyle{r\ge1,\ s\ge0}\\
\scriptstyle{r+s=i}\end{array}}a_{rs}\gamma_r\alpha_s,\quad i=1,...,k-1\ ,\label{eqedge11}\end{equation}
\begin{equation}0=\sum_{\renewcommand{\arraystretch}{0.6}
\begin{array}{c}
\scriptstyle{r\ge1,\ s\ge0}\\
\scriptstyle{r+s=i}\end{array}}a_{rs}\gamma_r\alpha_s,\quad i=k,k+1\ ,\label{eqedge12}\end{equation}
where $a_{rs}>0$ for all $r,s$, and we set $\alpha_0=1$, $\alpha_i=0$ as $i\ge k$. Finally, we add one more equation coming from the comparison of the constant terms in (\ref{eqedge1})
\begin{equation}2=\frac{9}{(2k+1)^2}\alpha_{k-1}^2\gamma_3\ .\label{eqedge13}\end{equation}
Substituting expressions (\ref{eqedge11}) to (\ref{eqedge12}) and (\ref{eqedge13}) several times, we end up with the quasihomogeneous equations
\begin{equation}\sum_{\renewcommand{\arraystretch}{0.6}
\begin{array}{c}
\scriptstyle{r,s,t\ge0}\\
\scriptstyle{r+2s+3t=i}\end{array}}a_{rst}\gamma_1^r\gamma_2^s\gamma_3^t=0,\quad i=k,k+1\ ,
\label{eqedge10}\end{equation}
\begin{equation}\sum_{\renewcommand{\arraystretch}{0.6}
\begin{array}{c}
\scriptstyle{r,s,t\ge0}\\
\scriptstyle{r+2s+3t=2k+1}\end{array}}a_{rst}\gamma_1^r\gamma_2^s\gamma_3^t=2\ ,
\label{eqedge14}\end{equation}
with all coefficients $a_{rst}$ positive. Consider solutions to the system (\ref{eqedge10}) as
intersections of two curves on a suitable toric surface. Notice that there are no intersections on toric divisors, since otherwise one would get a one-parameter family of solutions to the system
(\ref{eqedge10}), (\ref{eqedge14}) contrary to Lemma \ref{ledge1}(iii). Hence all intersection lie in the big torus.
Passing to the variables
$\gamma_2\gamma_1^{-2},\gamma_3\gamma_1^{-3}$, we obtain in (\ref{eqedge10}) two polynomial equations with Newton polygons $\sigma_{2,3}(k)$ and $\sigma_{2,3}(k+1)$, respectively, which yield $\langle\sigma_{2,3}(k),\sigma_{2,3}(k+1)\rangle$ solutions (counting multiplicities) by Bernstein-Koushnirenko theorem
\cite{B}. Each of these solutions yields
$(2k+1)$ values for $\gamma_1$ from (\ref{eqedge14}). Note that formulas (\ref{erev2_2}) and (\ref{eqedge11}) yields a scheme-theoretical isomorphism between the solutions to system (\ref{eqedge10})-(\ref{eqedge14}) and the solutions to system (\ref{eqedge1}). Since the latter ones are simple (as we noticed above), the former ones are simple too.

In the same manner, we treat systems (\ref{eqedge2}), (\ref{eqedge7}), and (\ref{eqedge8}), obtaining respectively $(2k+1)[\sigma_{2,3}(k),\sigma_{2,3}(k+1)]$ solutions in the first case and $\quad$ \mbox{$2k\langle\sigma_{2,3}(k-2),\sigma_{2,3}(k-1)\rangle$} solutions in each of the two remaining cases.

The same approach we explore for system (\ref{eqedge3}). Equating the derivatives of the right-hand sides, we get
$$x^{2s}Q_{k-s-1}(((2s+1)R_2+xR'_2)S_{k-s-1}+2xR_2S'_{k-s-1})=T_k(2T'_kR_1+T_kR'_1)$$
as well as $s=1$ and the simplicity of the roots of $S_{k-2}$ and $T_k$.
Furthermore,
$$\begin{cases}&((2s+1)R_2+xR'_2)S_{k-s-1}+2xR_2S'_{k-s-1}=(2k+1)T_k\\
&2T'_kR_1+T_kR'_1=(2k+1)x^2Q_{k-2}\end{cases}$$
Substituting
$$S_{k-2}=x^{k-2}+\alpha_1x^{k-3}+...+\alpha_{k-2},\quad T_k=x^k+\beta_1x^{k-1}+...+\beta_k\ ,$$
$$R_2=x^2+\gamma_1x+\gamma_2,\quad R_1=x+\delta_1$$ to the above system, we, first, obtain
$$\beta_i=\frac{2k+1-2i}{2k+1}\alpha_i+\sum_{\renewcommand{\arraystretch}{0.6}
\begin{array}{c}
\scriptstyle{r,s\ge1\ge0}\\
\scriptstyle{r+s=i}\end{array}}a_{rs}\alpha_r\gamma_s,\quad i=1,...,k\ ,$$ where $a_{rs}>0$ for all $r,s$,\and we set $\alpha_r=0$ as $r\ge k-1$ and $\gamma_s=0$ as $s\ge3$, and then
\begin{equation}\frac{4i(2k+1-i)}{2k+1}\alpha_i=\sum_{\renewcommand{\arraystretch}{0.6}
\begin{array}{c}
\scriptstyle{r,s,t\ge0,\ t<i}\\
\scriptstyle{r+s+t=i}\end{array}}a_{rst}\gamma_r\delta_s\alpha_t,
\quad i=1,...,k-2\ ,\label{eqedge15}\end{equation}
\begin{equation}0=\sum_{\renewcommand{\arraystretch}{0.6}
\begin{array}{c}
\scriptstyle{r,s,t\ge0,\ t<i}\\
\scriptstyle{r+s+t=i}\end{array}}a_{rst}\gamma_r\delta_s\alpha_t,
\quad i=k-1,k\ ,\label{eqedge16}\end{equation} where $a_{rst}>0$ for all $r,s,t$, and we set $\delta_0=\gamma_0=\alpha_0=1$,
$\delta_r=0$ as $r\ge2$, $\gamma_s=0$ as $s\ge3$, $\alpha_t=0$ as $t\ge k-1$. Iteratively plugging the right-hand sides of (\ref{eqedge15}) to (\ref{eqedge16}), we expel all 
$\alpha_i$'s and come to
the system
$$0=\sum_{\renewcommand{\arraystretch}{0.6}
\begin{array}{c}
\scriptstyle{r,s,t\ge0,\ t<i}\\
\scriptstyle{r+2s+t=i}\end{array}}b_{rst}\gamma_1^r\gamma_2^s\delta_1^t,\quad i=k-1,k,\quad b_{rst}>0
\ \text{for all}\ r,s,t\ ,$$ which together with the constant term relation in (\ref{eqedge3})
$$2=\beta_k^2\delta_1$$ gives $2k[\sigma_{1,2}(k-1),\sigma_{1,2}(k)]$ solutions in total.
Similarly we solve the systems (\ref{eqedge4}), (\ref{eqedge7}), and (\ref{eqedge8}), completing the proof of the lemma. \proofend

We are now in a position to determine the multiplicity $\mu_c(\Gamma_c,\bp_c,h_c)$; more precisely, we count
how many compatible triples $(C_1,C_2,C)$, where $C_1,C_2$ are rational, peripherally
smooth and unibranch, $C$ as in Lemma \ref{ledge2}, satisfy appropriate initial data. Compatibility means that
the intersection points of $C$ with the toric divisors
$\Tor([(0,0),(1,m)])$, $\Tor([(1,m),(2,0)])$ are determined by $C_1,C_2$, respectively. The initial data
depend in the regular orientation of the fragment induced from $(\Gamma,\bp,h)$.

If $\Gamma_c$ does not contain
a marked point, i.e., $\bp_c=\emptyset$, then the regular orientation of the fragment $(\Gamma',h')$ is as shown in Figure \ref{fig-cfc}(d) (up to exchange of $V_1,V_2$).
Denote by $m_1,m_2,m_3$ the
weights of the outer edges oriented towards $V_1$ or $V_2$. It then follows from \cite[Lemma 3.9]{Shustin2005} (cf., the factors described in the last paragraph of Section \ref{nod-mod}, part (2)) and Lemma \ref{ledge2} that the number of required triples $(C_1,C_2,C)$ equals
\begin{equation}\frac{\mu(V_1)\mu(V_2)}{m_1m_2m_3}\cdot\mu_c(\Gamma_c,h_c),\quad
\text{where}\quad\mu_c(\Gamma_c,h_c)=
\frac{\theta(m)}{m}\ .
\label{ecfc1}\end{equation}

If $\Gamma_c$ contains a marked point, then the regular orientation looks like in
Figure \ref{fig-cfc}(e). Denote by $m_1,m_2$ the weights of the outer edges
oriented towards $V_1$ or $V_2$. In addition to \cite[Lemma 3.9]{Shustin2005}
and Lemma \ref{ledge2}, we apply the result of \cite[Section 2.3.9]{IMS} that
takes into account the condition to pass through the given marked point. It follows that the number of required triples $(C_1,C_2,C)$ equals
\begin{equation}\frac{\mu(V_1)\mu(V_2)}{m_1m_2}\cdot\mu_c(\Gamma_c,\bp_c,h_c),\quad
\text{where}\quad\mu_c(\Gamma_c,\bp_c,h_c)=
\frac{\theta(m)}{m}\ .
\label{ecfc2}\end{equation}

\subsection{Patchworking}\label{sec-pw}
The following statement is complementary to Theorem \ref{tuc1}.

\begin{theorem}\label{t2}
Let $\Delta^*\subset\Z^2\setminus\{0\}$ be a nondegenerate, primitive, balanced multiset, and $\Delta$ the induced convex lattice polygon. Let $Q\subset\partial \Delta\cap\Z^2$ be such that
${\mathcal Q}=\conv(Q)$ is a nondegenerate quadrangle without parallel edges (cf. Lemma \ref{luc1}), $0\le g<p_a(\Delta)-1$,
$n=|\partial\Delta\cap\Z^2|+g-2$, and $\bw$ a configuration of $n$ distinct points in $(\K^*)^2$ such that
$\bx=\val(\bn(\bw))$ is a set
of $n$ distinct point in $\R^2$ in general position. Let $(\Gamma,\bp,h,\widehat g)$ be a plane cuspidal $n$-marked
tropical curve of degree $\Delta^*$ and genus $g$ such that $h(\bp)=\bx$ and its cuspidal tropical
fragment $(\Gamma_c,\bp_c,h_c)$ is of type A, B, or C. Assume in addition that, if $(\Gamma,\bp,h)$ has a cuspidal fragment of type B, then the flat trivalent vertex is incident to at most one end of $\Gamma$ (see Lemma \ref{lmod}).
Then the number of curves $C\in V_{\Delta,g}(A_2)$ passing through $\bw$ and tropicalizing to $(\Gamma,
\bp,h,\widehat g)$ equals
\begin{equation}\mu_c(\Gamma_c,\bp_c,h_c)\cdot\prod_{V\in\Gamma^0\setminus\Gamma_c}\mu(V)\ ,
\label{epw}\end{equation}
where the value $\mu_c(\Gamma_c,\bp_c,h_c)$ should be appropriately chosen from formulas
(\ref{emuc1}), (\ref{emuc}), (\ref{ecfc1}), or (\ref{ecfc2}).
\end{theorem}

{\bf Proof.}
We closely follow the ideas of \cite{Shustin2005}. Recall that the number (\ref{epw}) counts in how many ways
the given curve $(\Gamma,\bp,h,\widehat g)$ and the point configuration $\bw\subset{\mathfrak X}$ can be enhanced to a
PTL, that is, the tropical limit which matches the given configuration
$\bw$ and includes the modification along the ordinary and extended edges of $\Gamma$ of weight $>1$ as well as modification of the cuspidal tropical fragment as defined in Sections \ref{nod-mod},
\ref{sec-cfa}, \ref{sec-cfb}, and \ref{sec-cfc}.

Then we apply the patchworking statement in \cite[Theorem 2.4]{Sh} to conclude that an enhanced tropical limit as above lifts to a unique family of curves $C^{(t)}\subset{\mathfrak X}_t\simeq\Tor(\Delta)$, $t\in(\C,0)\setminus\{0\}$, (cf. Sections \ref{sec-tl} and \ref{sec-ptl}) such that
$C^{(t)}\in V_{\Delta,g}(A_2)$ and $C^{(t)}\supset\bw(t)$ as $t\in(\C,0)\setminus\{0\}$.

We only comment on the following issues specific to the considered case.
\begin{itemize}\item
The limit curves in \cite[Theorem 2.4]{Sh} assumed to have semiquasihomogeneous singularities on the toric divisors, but the proof literally extends to the case of Newton nondegenerate singularities (as those which occur in Section \ref{sec-cfb}).
\item In Lemma \ref{lquad1} and consequently, in modifications in Section \ref{sec-cfb}, the deformation in a neighborhood of the singular branch of the limit curve along the family $C^{(t)}$, $t\in(\C,0)$, in not equisingular, but equiclassical\footnote{For precise definition of the equiclassical deformations, see \cite[Page 433, item (3), and Section 6]{DH}.}, that is, preserving a cuspidal singular branch and the total $\delta$- and $\varkappa$-invariant in a neighborhood of each singular point
in the torus $\Tor^*(Q)$, $Q$ being the quadrangle from either Lemma \ref{lquad1}, or in the
subdivisions shown in Figure \ref{ftrap1}(e,f). This means, for instance, that the limit curve may have a cuspidal singular branch and also several smooth branches centered at the same point. We notice that the whole equiclassical family is smooth in this case by \cite[Theorem 27]{Di}, and
the curves $C^{(t)}$, $t\ne0$ have nodes and one cusp as their singularities
by Lemma \ref{luc1} and due to the general position of the configurations $\bw(t)
\subset\Tor(\Delta)$, $t\ne0$.
\item The required transversality conditions for \cite[Theorem 2.4]{Sh} follow from \cite[Lemmas 3.5 and 3.9]{Shustin2005}, Lemmas \ref{lquad1}(ii) and \ref{ledge1}(iii) precisely like in \cite[Section 5.3]{Sh0} or in \cite[Section 2.5.10]{IMS}. The two following facts are crucial:

- the existence of a regular orientation on the marked curve $(\Gamma,\bp)$ and on the
tropical curves that appear in the modifications
yields a partial order on the set of all limit curves including modifications, which then can be extended to a linear order;

- moving along the set of all limit curves in the aforementioned linear order, each time the next
limit curve, matching the conditions imposed by all the preceding limit curves, is determined up to a finite choice; moreover, each choice corresponds to a transverse intersection of a certain equisingular (or equiclassical) family having an expected dimension with the linear system coming from the imposed conditions (see \cite[Lemmas 3.5 and 3.9]{Shustin2005}, Lemmas \ref{lquad1}(ii) and \ref{ledge1}(iii)). \proofend
\end{itemize}

\section{Lattice path algorithm}\label{sec-lp}

Theorems \ref{tuc1} and \ref{t2} reduce the enumeration of unicuspidal algebraic curves in a toric surface
to enumeration of unicuspidal plane tropical curves passing though an appropriate configuration of points in the plane, provided, that the counted algebraic curves tropicalize to unicuspidal tropical curves having a cuspidal fragment of type A, B, or C. Here we present a version of Mikhalkin's lattice path algorithm \cite[Section 7.2]{M} adapted to our setting and solving the stated enumerative problem of computing $\deg V_{\Delta,g}(A_2)$ in the following sense: we exhibit a configuration of points in $\R^2$ such that the unicuspidal plane tropical curves
of a given degree and genus passing through the chosen configuration have a cuspidal fragment of type A, B, or C, and we reduce the enumeration of these tropical curves to a finite combinatorial problem of enumeration of lattice paths and related subdivisions of the Newton polygon.
\medskip

Let $\Delta^*\subset\Z^2\setminus\{0\}$ be a nondegenerate, primitive, balanced multiset, and $\Delta$ an induced convex lattice polygon. Let $Q\subset(\partial\Delta\cap\Z^2)$ be such that
${\mathcal Q}=\conv(Q)$ is a nondegenerate quadrangle without parallel edges (cf., Lemma \ref{luc1}) and such that the Newton polygon $\Delta$ is {\it h-transversal}, which, in particular, means that the intersection of $\Delta$ with any vertical line $x=i$, $i\in\Z$, is either empty, or an integral point, or a lattice segment. Let $0\le g<p_a(\Delta)-1$ and $n=|\partial\Delta\cap\Z^2|+g-2$.

Fix a linear functional
\begin{equation}\lambda:\R^2_\Z\to\R,\quad\lambda(x,y)=x-y\eps,\quad0<\eps\ll1\ .\label{elf}\end{equation}
It defines a linear order on the set of integral points $\Delta\cap\Z^2$:
$$(i,j)\prec(i',j')\quad\Longleftrightarrow\quad\lambda(i,j)<\lambda(i',j')\quad
\Longleftrightarrow\quad i<i'\ \text{or}\ \begin{cases}i=i',&\\
j>j'&\end{cases}$$ A $\lambda$-monotone lattice path (briefly, $\lambda$-path) of length $m\ge1$ in $\Delta$ is a broken line with vertices
$(\omega_0,...,\omega_m)\subset \Delta\cap\Z^2$ such that $\omega_i\prec\omega_{i+1}$ for all $i=0,...,m-1$,
and $$\omega_0=\omega_{\min}:=\min(\Delta\cap\Z^2),
\quad\omega_m=\omega_{\max}:=\max(\Delta\cap\Z^2)\ .$$
Introduce also the linear functional $\lambda^\perp(x,y)=\eps x+y$. For a given $\lambda$-path $\pi$,
any line $L_a=\{\lambda=a\}$ with $\lambda(\omega_{\min})\le a\le\lambda(\omega_{\max})$ splits into two rays
\begin{eqnarray}&L_a^+(\pi)=\{q\in L_a\ :\ \lambda^\perp(q)\ge\lambda^\perp(L_a\cap\pi)\},\nonumber\\
&L_a^-(\pi)=\{q\in L_a\ :\ \lambda^\perp(q)\le\lambda^\perp(L_a\cap\pi)\}.\nonumber\end{eqnarray}
Set
$$\Pi^+(\pi)=\bigcup_{\lambda(\omega_{\min})\le a\le\lambda(\omega_{\max})}L_a^+(\pi),\quad
\Pi^-(\pi)=\bigcup_{\lambda(\omega_{\min})\le a\le\lambda(\omega_{\max})}L_a^-(\pi)\ .$$

The following algorithm constructs a subdivision of $\Delta$ into convex lattice polygons:
\begin{enumerate}\item[(i)] Choose a $\lambda$-path $\pi_0$ of length $n$.
\item[(ii)] Define two lattice paths $\pi^+_0=\pi^-_0=\pi_0$, then construct two finite sequences of
$\lambda$-paths $\{\pi^+_i\}_{i\ge0}$ and $\{\pi^-_i\}_{i\ge0}$ following the recipe described in the next item.
\item[(iii)] Given a $\lambda$-path $\pi^+_i$, $i\ge0$, with vertices
$(v_0,v_1,...,v_m)$, take
\begin{eqnarray}&k=\min\{1\le s\le m\ :\ v_{s-1},v_s,v_{s+1}\ \text{are not collinear,}\nonumber\\
&\qquad\qquad\qquad\qquad\qquad \text{and}\ [v_{s-1},v_{s+1}]\subset\Pi^+(\pi^+_i)\}\ ,\nonumber\end{eqnarray}
then perform one of the following operations
\begin{enumerate}\item[(a)] either define $\pi^+_{i+1}$ by the sequence of vertices
$(v_0,...,v_{k-1},v_{k+1},...,v_m)$;
\item[(b)] or define $\pi^+_{i+1}$ by the sequence of vertices
$(v_0,...,v_{k-1},v'_k,v_{k+1},...,v_m)$, where
$\conv\{v_{k-1},v_k,v_{k+1},v'_k\}$ is a parallelogram and $v'_k\in \Delta$;
\item[(c)] or define $\pi^+_{i+1}$ by the sequence of vertices
$(v_0,...,v_{k-1},v'_k,v_{k+1},...,v_m)$, where $\conv\{v_{k-1},v_k,v_{k+1},v'_k\}$ is a
nondegenerate quadrangle without parallel sides and $v'_k\in \Delta$;
\item[(d)] or define $\pi^+_{i+1}$ by the sequence of vertices
$(v_0,...,v_{k-1},v'_k,v_{k+1},...,v_m)$, where $\conv\{v_{k-1},v_k,v_{k+1},v'_k\}$ is a
trapeze and $v'_k\in \Delta$.
\end{enumerate} Add either the triangle $P_i^+=\conv\{v_{k-1},v_k,v_{k+1}\}$ in case (a), or
the quadrangle $P_i^+=\conv\{v_{k-1},v_k,v_{k+1},v'_k\}$ in cases (b)-(d) to the list of elements of the subdivision. If such $k$ does not exist, then $\pi_i^+$ is called terminal, and the sequence
$\{\pi_s^+\}_{0\le s\le i}$ is completed.
\item[(iv)] The same operations we perform with $\pi_i^-$ replacing everywhere plus by minus.
\end{enumerate}
The outcome of the above procedure is called an admissible subdivision of $\Delta$ if
\begin{itemize}\item the union of all obtained polygons $P_i^\pm$, $i\ge0$, equals $\Delta$,
and all the points of $\partial \Delta\cap\Z^2$ are vertices of that subdivision;
\item at most one polygon $P_i^\pm$, $i\ge0$, is not a triangle or a parallelogram;
\item if all the polygons $P_i^\pm$, $i\ge0$, are triangles or parallelograms, then there exists
an edge of the subdivision of lattice length $\ge3$;
\item the subdivision is dual to an irreducible plane tropical curve in the following sense: first,
we take the embedded plane tropical curve dual to the subdivision (we show later that it does exist), then
resolve all self-intersection points dual to the parallelograms; the resulting graph must be connected.
\end{itemize} We convert an admissible subdivision into a {\it marked admissible subdivision}: if it contains
a polygon $P_i^\pm$, $i\ge0$, different from a triangle or a parallelogram, we mark this polygon, otherwise we mark one of the edges of the subdivision of lattice length $\ge3$.

To a pair $(\pi_0,\Sigma)$, where $\pi_0$ is a $\lambda$-path of length $n$,
$\Sigma$ a marked admissible subdivision of $\Delta$ obtained in the previous procedure starting with $\pi_0$, we assign the following multiplicity
$\mu(\pi_0,\Sigma)$:
\begin{itemize}\item if $\Sigma$ consists only of triangles and parallelograms, we set
$\mu(\pi_0,\Sigma)$ to be the product of lattice areas of all the triangles times $\frac{\theta(m)}{m}$, where
$m$ is the lattice length of the marked edge, $\theta$ introduced in Definition \ref{dcfc};
\item if $\Sigma$ contains a nondegenerate quadrangle $Q=\conv\{v_{k-1},v_k,v_{k+1},v'_k\}$ without parallel sides
(in the notation of item (c) above), then $\mu(\pi_0,\Sigma)$ equals the product of lattice areas of all triangles
times the lattice area of the triangle $\conv\{v_{k-1},v_k,v_{k+1}\}$;
\item if $\Sigma$ contains a trapeze $Q=\conv\{v_{k-1},v_k,v_{k+1},v'_k\}$ with bases of lengths
$m>m_1$, then $\mu(\pi_0,\Sigma)$ equals the product of lattice areas of all triangles
times $(m+m_2)$ if $v_k$ is incident to the long parallel side, or
$\frac{(m+m_2)m_1}{m}$ otherwise (here $m_2=m-m_1$).
\end{itemize}

\begin{theorem}\label{t4}
Let $\Delta^*\subset\Z^2\setminus\{0\})$ be a nondegenerate, primitive, balanced multiset
satisfying the hypotheses of Lemma \ref{luc1} and such that the induced convex lattice polygon $\Delta\subset\R^2_\Z$ is {\it h-transversal}.
Let $0\le g<p_a(\Delta)-1$, $n=|\partial\Delta\cap\Z^2|+g-2$, and $\lambda$ be given by (\ref{elf}). Then
$$\deg V_{\Delta,g}(A_2)=\sum_{(\pi_0,\Sigma)}\mu(\pi_0,\Sigma)\ ,$$
where $\pi_0$ ranges over all $\lambda$-paths in $\Delta$ of length $n$, and $\Sigma$ runs over all marked admissible subdivisions of $\Delta$ arising from $\pi_0$.
\end{theorem}

{\bf Proof.} {\it Step 1: Preliminary observations.} The condition of the primitivity and $h$-transversality is equivalent to the following one:

\smallskip
{\it (O1) Each vector $\oa\in\Delta^*$ either equals $(\pm1,0)$, or satisfies $\pr_v\oa=\pm1$, where $\pr_v$ is the
projection on the vertical axis in $\R^2_\Z$.}
\smallskip

Fix a configuration $\bx$ of $n$ points on the line $L_{1,-\eps}=\{y=-\eps x\}$:
\begin{equation}
\bx_i=(M_i,-M_i\eps),\ i=1,...,n,\quad 0<M_1\ll M_2\ll...\ll M_n\ .\label{econfig}\end{equation}
These points are in a tropically general position (see \cite[Theorem 2]{M}). If $(\Gamma,\bp,h)$ is a cuspidal
tropical curve of degree $\Delta^*$ and genus $g$ such that $h(\bp)=\bx$, then by Lemma \ref{ltrop2}, $(\Gamma,\bp,h)$ is
regular, and the points $\bx_i$, $i=1,...,n$, lie on edges of the embedded curve $h_*(\Gamma)$. Furthermore, by Lemma \ref{lmod}, we can suppose that, if $(\Gamma,\bp,h)$ has a cuspidal fragment of type B, then its flat trivalent vertex is incident to at most one end of $\Gamma$. In particular, it follows that all ends of $h_*(\Gamma)$ are of weight $1$ and are in one-to-one correspondence with the vectors of $\Delta^*$.

\smallskip

{\it Step 2: The case of a tropical curve $(\Gamma,\bp,h)$ with a cuspidal fragment of type C, D, or E.} In this situation, the reduced curve $(\Gamma_{red},\bp_{red},h_{red})$ is trivalent of degree $\Delta^*$ with $b_1(\Gamma_{red})=g-1$. Since
$n=|\Delta^*|+(g-1)-1$, we get to the setting of the original Mikhalkin's lattice path algorithm \cite[Section 7.2]{M}, which
reduces the enumeration of all the considered plane tropical curves to the enumeration of appropriate subdivisions of the Newton polygon into triangles and parallelograms (dual to the embedded tropical curves $h_*(\Gamma)$).

Now we make further important observations (cf. \cite[Section 2.5.6]{IMS}, \cite[Section 7.2]{M}
and \cite[Proof of Proposition 2.6]{IKS1}):

\smallskip
{\it (O2) The edges of $h_*(\Gamma)\setminus\bx$ are regularly oriented so the projections of their directing vectors
to the vector $(\eps,1)$, orthogonal to the line $L_{1,-\eps}$, are positive or negative according as the corresponding edges of $h_*(\Gamma)\setminus\bx$ lie above or below the line $L_{1,-\eps}$.}

This immediately follows from the balancing condition and the fact that the orienting vectors of the edges emanating from the fixed points $\bx$ satisfy the asserted condition (see Figure \ref{fig-lp}(a)). It also follows that the vertical projections of the orienting vectors are nonnegative or nonpositive according as the corresponding edges of $h_*(\Gamma)\setminus\bx$ lie above or below the line $L_{1,-\eps}$. The latter remark together with the balancing condition implies the observation

\smallskip
{\it (O3) For any vertex of $h_*(\Gamma)$, dual to a triangle or a parallelogram in the subdivision of the Newton polygon, $$\max|\pr_v\oa|\le\max|\pr_v\ob|\ ,$$
where $\oa$ runs over the directing vectors of incoming edges and $\ob$ runs over the directing vectors
of outgoing edges, which is straightforward from the balancing condition (see Figure \ref{fig-lp}(a)).}\smallskip

\begin{figure}
\setlength{\unitlength}{1.0mm}
\begin{picture}(130,40)(0,0)
\thicklines

\put(0,10){\line(1,0){55}}\put(60,10){\line(1,0){20}}
\put(85,10){\line(1,0){20}}\put(110,10){\line(1,0){20}}

\thinlines
\put(5,10){\vector(1,1){10}}\put(25,10){\vector(-1,1){10}}
\put(15,20){\vector(0,1){10}}\put(30,10){\vector(1,1){10}}
\put(50,10){\vector(-1,1){10}}\put(40,20){\vector(-1,1){10}}
\put(40,20){\vector(1,1){10}}\put(65,15){\vector(1,1){5}}
\put(75,15){\vector(-1,1){5}}\put(70,20){\vector(0,1){10}}
\put(70,20){\vector(0,-1){5}}\put(95,14){\vector(0,1){5}}
\put(95,20){\vector(0,1){10}}\put(95,20){\vector(-1,-1){5}}
\put(100,15){\vector(-1,1){5}}\put(115,15){\vector(1,1){5}}
\put(125,15){\vector(-1,1){5}}\put(120,20){\vector(2,1){10}}
\put(120,20){\vector(-1,2){5}}

\dottedline{1}(60,10)(65,15)\dottedline{1}(70,10)(70,15)
\dottedline{1}(75,15)(80,10)\dottedline{1}(85,10)(90,15)
\dottedline{1}(95,10)(95,15)\dottedline{1}(100,15)(105,10)
\dottedline{1}(110,10)(115,15)\dottedline{1}(125,15)(130,10)

\put(4,9){$\bullet$}\put(24,9){$\bullet$}\put(29,9){$\bullet$}\put(49,9){$\bullet$}
\put(25,0){(a)}\put(80,0){(b)}\put(117,0){(c)}

\end{picture}
\caption{Construction of a tropical curve via the lattice path algorithm}\label{fig-lp}
\end{figure}

An immediate consequence of (O1) and (O3) is

\smallskip{\it (O4) Each triangle in the dual to $h_*(\Gamma)$ subdivision of the Newton polygon $\Delta$
has a vertical side. Each edge of $\Sigma$ of lattice length $\ge2$ is vertical.}
\smallskip

\begin{lemma}\label{O5} The curve $(\Gamma,\bp,h)$ cannot have a cuspidal fragment of type D or E.
\end{lemma}

{\bf Proof.}
For the cuspidal fragment of type E,
the lattice triangle $T$ dual to the
elliptic vertex must contain at least two interior integral points, since otherwise an elliptic limit curve cannot be singular. It follows that $T$ necessarily has a side with the
horizontal projection $\ge2$.

For the cuspidal fragment of type D, consisting of a trivalent vertex and a four-valent vertex joined by two edges, let $T$ be the lattice triangle dual to the
four-valent vertex of $\Gamma_c$. The limit curve $C_T$ in $\Tor(T)$ has a singular point, center of two local branches on the toric divisor $\Tor(e)$, where $e$ is the side of $T$ dual to the $h$-image of the flat cycle of $\Gamma$ (cf. Section \ref{sec-ptl}). It follows that the lattice length of $e$ is at least $2$. In case $e$ is vertical, the horizontal projection of $T$ is $\ge2$, since otherwise the limit curve $C_T$ would be smooth. Hence, $T$ always has a side with horizontal projection $\ge2$.

For the cuspidal fragment of type D, consisting of two four-valent vertices joined by two edges, denote by $T_1,T_2$ the lattice triangles dual to these four-valent vertices, an by $e$ their common side dual to the $h$-image of the flat cycle. Each of the limit curves $C_{T_1}\subset\Tor(T_1)$ and $C_{T_2}\subset\Tor(T_2)$ is rational with two local branches centered at $\Tor(e)$. Thus, $|e|_\Z\ge2$. Suppose that the edge $E$ is vertical. Then we claim that the horizontal projection of either $T_1$, or $T_2$ is $\ge2$. Indeed, otherwise the curves $C_{T_1}$ and $C_{T_2}$ must be smooth, each having two smooth local branches centered at two distinct points $z_1,z_2\in\Tor(E)$. It follows that the limit of the cuspidal singularity in the family of curves (\ref{euc2}) must be either $z_1$, or $z_2$. However, this contradicts the fact that the modification of the tropical limit at $z_1$ and $z_2$ yields extra rational nodal curves,
see \cite[Section 3.5]{Shustin2005}.

Thus, in each of the above cases we always have an edge of the tropical curve with vertical projection $\ge2$. Combining this conclusion with observations (O1), (O3), and (O4), we derive the statement of the lemma.
\proofend

{\it Step 3: The case of a tropical curve $(\Gamma,\bp,h)$ with a cuspidal fragment of type A or B.} In this situation, the embedded curve $h_*(\Gamma)$ contains
a four-valent vertex $V$ dual a nondegenerate quadrangle different from a parallelogram, while the other tiles of the subdivision of $\Delta$ are triangles and parallelograms.

\smallskip
{\it Claim:} The observation (O2) is valid in the considered case.
\smallskip

To prove the claim, we argue on the contrary. Suppose that the four-valent vertex $V$ lies above the line $L_{1,-\eps}$. Introduce coordinates $x',y'$ in $\R^2$ so that $L_{1,-\eps}$ is the $x'$-axis, and the $y'$-axis is oriented by the vector $(\eps,1)$. Suppose that there exist edges of $h_*(\Gamma)\setminus\bx$ whose oriented vectors have negative $y'$-coordinate. Among the initial vertices of these edges pick the vertex $V'$ with the maximal $y'$-coordinate. This vertex cannot be dual to a triangle or a parallelogram, since otherwise, by the balancing condition, there should be an incoming edge, whose directing vector has a negative $y'$-coordinate contrary to the choice of $V'$. Thus, $V'=V$, and the two incoming edges incident to $V$ are directed by vectors with positive $y'$-coordinates. It follows that the fourth edge incident to $V$ is directed by a vector with a positive $y'$-coordinate and, furthermore, the vertical projection of the directing vector is $\ge2$ (see an illustration in Figure \ref{fig-lp}(b)). Propagating this edge further, we encounter only vertices which are $y'$-higher than $V$, hence dual to triangles and parallelograms, where two incoming edges having positive $y'$-coordinates, and finally we end up
with an infinite edge directed by a vector with the vertical projection $\ge2$ contrary to the condition (O1). Noticing that the behavior of $h_*(\Gamma)$ under the line $L_{1,-\eps}$ is as described in the preceding Step 2, we complete the proof of the claim.

One of the consequences is that $h_*(\Gamma)\cap L_{1,-\eps}=\bx$, and hence by duality we obtain

\begin{lemma}\label{O6} The lattice path is connected.\end{lemma}

This, in particular, justifies the fact that the lattice path algorithm described in the beginning of Section \ref{sec-lp} starts with a connected $\lambda$-path $\pi_0$ (see stage (i) of the algorithm).
Another consequence of the above Claim is that, in the stage (iii) of the lattice path algorithm, each intermediate lattice path $\pi^-_i$ (or $\pi^-_i$) appears to be a $\lambda$-path.

\smallskip
{\bf(4)} Finally, the formulas for the multiplicities of the marked admissible subdivisions
are straightforward consequences of formula (\ref{epw}). We only notice that the middle orientation in Figure \ref{ftrap1}(g) is forbidden by (O2), and, if the long base of the trapeze is a part of $\pi_0$, then
the multiplicity of the subdivision corresponds either to the multiplicity of
the tropical curve with a fragment shown in Figure \ref{ftrap1}(l), or to the sum of multiplicities of the curves containing the fragments shown in Figures \ref{ftrap1}(j,k) (of course, with the same result).
\proofend

\section{Enumeration of real cuspidal curves}\label{sec-ex}

We start with a tropical proof of the following classically known formula. The elements of the proof
will be essential for the proof of a lower bound for the maximal number of real plane curves with several real cusps occurring in generic real linear systems of an appropriate dimension (Theorem \ref{tex} below).
Denote by $V_d(A_2)$ the family of irreducible plane curves of degree $d$ having an ordinary cusp as their only singularity. It is an irreducible quasiprojective variety of codimension $2$ in $|{\mathcal O}_{\PP^2}(d)|$.

\begin{lemma}\label{l-unicusp}
$$\deg V_d=12(d-1)(d-2).$$
\end{lemma}

{\bf Proof.} We accept the notation from Section \ref{sec-lp}. In the considered case, each lattice path $\pi_0$ has vertices in all but two integral points in the Newton triangle ${\mathcal T}_d=\conv\{(0,0),(d,0),(0,d)\}$ (we call these two integral points
{\it missing points}). It is not difficult to see that if the missing points $(i',j'),(i'',j'')$ satisfy the following conditions:
\begin{itemize}\item either $|i'-i''|\ge2$, or $|i'-i''|=1$ and $0<j'<d-i'$, $0<j''<d-i''$,
\item or $i'=i''$ and $|j'-j''|\ge2$, or $i'=i''=0$, or $j'=j''=0$,
\item or $j'=d-i'$ and $j''=d-i''$,
\item or $j'=0$ and $j''=d-i''$,\end{itemize}
then the lattice path does not induce an admissible subdivision of ${\mathcal T}_d$. Indeed, in these cases one encounters either a pair of vertical segments of length $2$, or a point on the boundary of the Newton triangle that is not a vertex of the subdivision. In both the situations, the subdivision is not dual to any cuspidal tropical curve of a given degree.

Suppose that $d\ge 5$. The possible combinations of missing points and admissible subdivisions are shown in Figure \ref{fex2} (where the upper incline line is a part of the segment $[(d,0),(0,d)]$, the missing points are designated by bullets,
the lattice path segments are marked by arrows, and the part of ${\mathcal T}_d$ outside the designated fragments is covered by lattice triangles of lattice area $1$):
\begin{enumerate}
\item[(a)] If the missing points are $(i,j),(i,j+1)$, where
$1\le i\le d-3$ and $1\le j\le d-i-2$ (see Figure \ref{fex2}(a)), then we have a unique subdivision of ${\mathcal T}_d$ with an edge $[(i,j-1),(i,j+2)]$ dual to a cuspidal fragment of type C. The multiplicity of the subdivision equals $6$ (see formulas (\ref{ecfc1}), (\ref{ecfc2}) and Theorem \ref{t2}), and the total contribution of all such subdivisions equals $6\cdot\frac{(d-3)(d-2)}{2}=3(d-3)(d-2)$.
\item[(b)] If the missing points are $(i,d-i),(i,d-i-1)$, $1\le i\le d-2$ (see Figure \ref{fex2}(b)), then
an admissible subdivision contains a trapeze $$\conv\{(i,d-i),(i,d-i-2),(i-1,j),(i-1,j+1)\}\ ,$$ where
$0\le j\le d-i$. Any such subdivision has multiplicity $3$ (see the middle formula in (\ref{emuc}) and Theorem \ref{t2}), and their total contribution equals $3\cdot\frac{(d-2)(d+3)}{2}=\frac{3}{2}d^2+\frac{3}{2}d-9$.
\item[(c)] If the missing points are $(i,0),(i,1)$, where $1\le i\le d-4$, we have several types of admissible subdivisions:
    \begin{itemize}\item either containing a trapeze spanned by the points $(i,0),(i,2),
    (i+1,j),(i+1,j+1)$ as $3\le j\le d-i-2$ (see Figure \ref{fex2}(c)); the multiplicity of such a subdivision is $3$ as in item (b);
    \item or containing the quadrangle $\conv\{(i-1,0),(i,0),(i,2),(i+1,3)\}$
    (see Figure \ref{fex2}(d)); the multiplicity of this subdivision is
    $1$ (see Lemma \ref{lquad1}(iii) and Theorem \ref{t2}));
    \item or containing the quadrangle $\conv\{(i-1,0),(i,0),(i+1,2),(i+1,3)\}$
    (see Figure \ref{fex2}(e)); the multiplicity of this subdivision is $2$ (cf. the preceding case); \item or containing the trapeze $\conv\{(i-1,0),(i,0),(i+1,1),(i+1,2)\}$; the multiplicity of this subdivision (see Figure \ref{fex2}(f)) is $6$ according to the last formula in (\ref{emuc}) and Theorem \ref{t2}.\end{itemize} The total contribution of these four types of admissible subdivisions amounts to
    $$3\cdot\frac{(d-5)(d-4)}{2}+(1+2+6)(d-4)=\frac{3}{2}d^2-\frac{9}{2}d-6\ .$$
\item[(d)] If the missing points are $(i,0),(i,1)$, where $i=d-3$ or $d-2$, we have subdivisions containing
    either the trapeze $\conv\{(d-4,0),(d-3,0),(d-2,1,d-2,2)\}$ (see Figure \ref{fex2}(g)), or the quadrangle $\conv\{(d-3,0),(d-2,0),(d-2,2),(d-1,1)\}$ (see Figure \ref{fex2}(h)), and having multiplicity $6$ or $3$, respectively.
\item[(e)] The remaining options for the pair of missing points are $(i,j),(i+1,d-i-1)$, where $1\le i\le d-2$, $1\le j\le d-i-1$ (see Figure \ref{fex2}(i)), and $(i,j),(i-1,0)$, where $1\le i\le d-2$, $1\le j\le d-i-1$ (see Figure \ref{fex2}(j)); in both cases the multiplicity of the subdivision equals $6$, which results in the total contribution $6(d-1)(d-2)$.
\end{enumerate}
The sum of the above contributions equals $12(d-1)(d-2)$.

If $d=4$, then the subdivisions mentioned in item (c) do not occur. However, their contribution counted formally amounts to $(3d^2/2-9d/2-6)\big|_{d=4}=0$, and hence the above enumeration yields the correct answer.

If $d=3$, then all admissible subdivisions and their weights are shown in Figure \ref{fex10}. The sum of weights is $24$ as required. \proofend

\begin{figure}
\setlength{\unitlength}{1mm}
\begin{picture}(125,130)(0,0)
\thinlines
\put(0,110){\vector(1,0){20}}\put(30,110){\vector(1,0){20}}
\put(60,110){\vector(1,0){20}}\put(90,110){\vector(1,0){20}}
\put(0,60){\vector(1,0){20}}\put(30,60){\vector(1,0){20}}
\put(60,60){\vector(1,0){30}}\put(100,60){\vector(1,0){25}}
\put(0,10){\vector(1,0){20}}\put(30,10){\vector(1,0){20}}
\put(5,110){\line(1,4){5}}\put(5,110){\line(1,1){5}}
\put(10,130){\line(1,0){5}}\put(10,115){\line(1,3){5}}
\put(35,125){\line(1,0){5}}\put(35,120){\line(1,1){5}}\put(35,115){\line(1,0){5}}
\put(40,115){\line(0,1){10}}\put(40,115){\line(1,1){5}}
\put(65,110){\line(1,2){10}}\put(70,110){\line(1,3){5}}
\put(70,110){\line(0,1){10}}
\put(95,110){\line(1,2){5}}\put(100,120){\line(1,1){5}}\put(100,110){\line(1,3){5}}
\put(100,110){\line(1,2){5}}\put(100,110){\line(1,1){5}}\put(100,125){\line(1,1){5}}
\put(100,130){\line(1,0){5}}
\put(5,60){\line(1,2){5}}\put(10,70){\line(1,1){5}}\put(5,60){\line(2,3){10}}
\put(10,60){\line(1,2){5}}\put(10,60){\line(1,1){5}}\put(10,75){\line(1,1){5}}
\put(10,80){\line(1,0){5}}
\put(35,60){\line(1,2){5}}\put(40,70){\line(1,1){5}}\put(35,60){\line(2,3){10}}
\put(35,60){\line(1,1){10}}\put(40,60){\line(1,1){5}}\put(40,75){\line(1,1){5}}
\put(40,80){\line(1,0){5}}
\put(65,60){\line(1,2){5}}\put(65,60){\line(1,1){10}}\put(70,60){\line(1,1){5}}
\put(75,65){\line(1,0){5}}
\put(110,60){\line(1,1){5}}
\put(5,10){\line(1,3){5}}\put(5,10){\line(1,1){10}}\put(10,25){\line(1,0){5}}
\put(15,20){\line(0,1){5}}
\put(35,10){\line(0,1){5}}\put(35,15){\line(1,2){5}}\put(35,10){\line(1,1){5}}
\put(40,15){\line(1,2){5}}\put(40,25){\line(1,0){5}}\put(35,20){\line(1,2){5}}
\put(35,25){\line(1,1){5}}\put(35,30){\line(1,0){5}}

\thicklines
\put(5,140){\line(1,-1){15}}\put(30,135){\line(1,-1){20}}
\put(65,145){\line(1,-1){15}}\put(95,140){\line(1,-1){15}}
\put(5,90){\line(1,-1){15}}\put(35,90){\line(1,-1){15}}
\put(65,80){\line(1,-1){20}}\put(105,75){\line(1,-1){15}}
\put(5,35){\line(1,-1){15}}\put(30,40){\line(1,-1){20}}

\put(10,135){\vector(0,-1){5}}\put(10,130){\vector(0,-1){15}}\put(10,115){\vector(0,-1){5}}
\put(10,110){\vector(1,4){5}}
\put(35,130){\vector(0,-1){5}}
\put(35,125){\vector(0,-1){5}}\put(35,120){\vector(0,-1){5}}\put(35,115){\vector(0,-1){5}}
\put(35,110){\vector(1,1){5}}\put(40,115){\vector(0,-1){5}}\put(40,110){\vector(1,2){5}}
\put(70,140){\vector(0,-1){5}}\put(70,135){\vector(0,-1){5}}\put(70,130){\vector(0,-1){5}}
\put(70,125){\vector(0,-1){5}}\put(70,120){\vector(1,3){5}}\put(75,135){\vector(0,-1){5}}
\put(75,130){\vector(0,-1){5}}\put(75,125){\vector(0,-1){5}}\put(75,120){\vector(0,-1){5}}
\put(75,115){\vector(0,-1){5}}
\put(100,135){\vector(0,-1){5}}\put(100,130){\vector(0,-1){5}}\put(100,125){\vector(0,-1){5}}
\put(100,120){\vector(1,2){5}}
\put(105,130){\vector(0,-1){5}}\put(105,125){\vector(0,-1){5}}\put(105,120){\vector(0,-1){5}}
\put(105,115){\vector(0,-1){5}}
\put(10,85){\vector(0,-1){5}}\put(10,80){\vector(0,-1){5}}\put(10,75){\vector(0,-1){5}}
\put(10,70){\vector(1,2){5}}
\put(15,80){\vector(0,-1){5}}\put(15,75){\vector(0,-1){5}}\put(15,70){\vector(0,-1){5}}
\put(15,65){\vector(0,-1){5}}
\put(40,85){\vector(0,-1){5}}\put(40,80){\vector(0,-1){5}}\put(40,75){\vector(0,-1){5}}
\put(40,70){\vector(1,2){5}}
\put(45,80){\vector(0,-1){5}}\put(45,75){\vector(0,-1){5}}\put(45,70){\vector(0,-1){5}}
\put(45,65){\vector(0,-1){5}}
\put(10,85){\vector(0,-1){5}}\put(10,80){\vector(0,-1){5}}\put(10,75){\vector(0,-1){5}}
\put(65,60){\vector(1,3){5}}\put(70,75){\vector(0,-1){5}}\put(70,70){\vector(1,0){5}}
\put(75,70){\vector(0,-1){5}}\put(75,65){\vector(0,-1){5}}\put(75,60){\vector(1,1){5}}
\put(80,65){\vector(0,-1){5}}\put(80,60){\vector(1,0){5}}
\put(105,60){\vector(1,2){5}}\put(110,70){\vector(1,-1){5}}\put(115,65){\vector(0,-1){5}}
\put(115,60){\vector(1,0){5}}
\put(10,30){\vector(0,-1){5}}\put(10,25){\vector(0,-1){10}}\put(10,15){\vector(0,-1){5}}
\put(10,10){\vector(1,2){5}}\put(15,20){\vector(0,-1){5}}\put(15,15){\vector(0,-1){5}}
\put(35,35){\vector(0,-1){5}}\put(35,30){\vector(0,-1){5}}\put(35,25){\vector(0,-1){5}}
\put(35,20){\vector(0,-1){5}}\put(35,15){\vector(1,3){5}}\put(40,30){\vector(0,-1){5}}
\put(40,25){\vector(0,-1){10}}\put(40,15){\vector(0,-1){5}}\put(40,10){\vector(1,3){5}}

\put(8,100){\rm (a)}\put(38,100){\rm (b)}\put(68,100){\rm (c)}\put(98,100){\rm (d)}
\put(8,50){\rm (e)}\put(38,50){\rm (f)}\put(73,50){\rm (g)}\put(108,50){\rm (h)}
\put(8,0){\rm (i)}\put(38,0){\rm (j)}\put(73,0){\rm (k)}\put(110,0){\rm (l)}
\put(9,106){$i$}\put(39,106){$i$}\put(69,106){$i$}\put(99,106){$i$}
\put(9,56){$i$}\put(39,56){$i$}\put(84,56){$d$}\put(119,56){$d$}
\put(9,6){$i$}\put(39,6){$i$}
\put(9,119){$\bullet$}\put(9,124){$\bullet$}
\put(39,119){$\bullet$}\put(39,124){$\bullet$}
\put(69,109){$\bullet$}\put(69,114){$\bullet$}
\put(99,109){$\bullet$}\put(99,114){$\bullet$}
\put(9,59){$\bullet$}\put(9,64){$\bullet$}
\put(39,59){$\bullet$}\put(39,64){$\bullet$}
\put(69,59){$\bullet$}\put(69,64){$\bullet$}
\put(109,59){$\bullet$}\put(109,64){$\bullet$}
\put(9,19){$\bullet$}\put(14,24){$\bullet$}
\put(34,9){$\bullet$}\put(39,19){$\bullet$}

\thicklines
\put(62,25){\line(1,0){25}}\put(62,25){\line(-1,-1){5}}\put(62,25){\line(-1,1){5}}
\put(72,25){\line(0,-1){5}}\put(72,25){\line(-1,1){5}}
\put(87,25){\line(1,1){5}}\put(87,25){\line(1,-1){5}}

\put(100,25){\line(1,-1){10}}\put(110,25){\line(1,-1){5}}\put(110,15){\line(1,1){5}}
\put(115,20){\line(2,1){10}}\put(110,15){\line(0,-1){5}}\put(115,20){\line(0,-1){10}}

\put(76,24){$\bullet$}\put(61,22){$v_1$}\put(71,27){$v_3$}\put(76,22){$v$}\put(84,27){$v_2$}
\put(114,24){$\bullet$}\put(99,27){$v_1$}\put(108,27){$v_3$}\put(114,27){$v$}\put(123,27){$v_2$}

\thinlines
\dashline{1}(100,25)(125,25)\dottedline{1}(115,20)(115,25)

\end{picture}
\caption{Enumeration of complex and real unicuspidal curves}\label{fex2}
\end{figure}

\begin{figure}
\setlength{\unitlength}{1.2mm}
\begin{picture}(90,65)(-10,5)
\thinlines
\put(0,10){\line(1,0){24}}\put(0,10){\line(0,1){24}}\put(24,10){\line(-1,1){24}}
\put(30,10){\line(1,0){24}}\put(30,10){\line(0,1){24}}\put(54,10){\line(-1,1){24}}
\put(60,10){\line(1,0){24}}\put(60,10){\line(0,1){24}}\put(84,10){\line(-1,1){24}}
\put(0,44){\line(1,0){24}}\put(0,44){\line(0,1){24}}\put(24,44){\line(-1,1){24}}
\put(30,44){\line(1,0){24}}\put(30,44){\line(0,1){24}}\put(54,44){\line(-1,1){24}}
\put(60,44){\line(1,0){24}}\put(60,44){\line(0,1){24}}\put(84,44){\line(-1,1){24}}
\put(8,10){\line(0,1){16}}\put(0,18){\line(1,1){8}}\put(0,26){\line(1,0){8}}
\put(38,10){\line(-1,1){8}}\put(30,26){\line(1,0){8}}\put(38,10){\line(0,1){16}}
\put(68,10){\line(-1,1){8}}\put(68,10){\line(-1,2){8}}\put(68,10){\line(0,1){16}}
\put(0,60){\line(1,0){8}}\put(30,52){\line(1,1){8}}\put(30,60){\line(1,0){8}}
\put(38,44){\line(1,1){8}}\put(60,52){\line(1,1){8}}\put(60,60){\line(1,0){8}}
\put(76,44){\line(0,1){8}}

\thicklines
\put(0,34){\vector(0,-1){8}}\put(0,26){\vector(0,-1){8}}\put(0,18){\vector(0,-1){8}}
\put(0,10){\vector(1,0){8}}\put(8,10){\vector(1,1){8}}\put(16,18){\vector(0,-1){8}}
\put(16,10){\vector(1,0){8}}
\put(30,34){\vector(0,-1){8}}\put(30,26){\vector(0,-1){8}}\put(30,18){\vector(0,-1){8}}
\put(30,10){\vector(1,0){8}}\put(38,10){\vector(1,1){8}}\put(46,18){\vector(0,-1){8}}
\put(46,10){\vector(1,0){8}}
\put(60,34){\vector(0,-1){8}}\put(60,26){\vector(0,-1){8}}\put(60,18){\vector(0,-1){8}}
\put(60,10){\vector(1,0){8}}\put(68,10){\vector(1,1){8}}\put(76,18){\vector(0,-1){8}}
\put(76,10){\vector(1,0){8}}
\put(0,68){\vector(0,-1){8}}\put(0,60){\vector(0,-1){8}}\put(0,52){\vector(1,1){8}}
\put(8,60){\vector(0,-1){16}}\put(8,44){\vector(1,1){8}}\put(16,52){\vector(0,-1){8}}
\put(16,44){\vector(1,0){8}}
\put(30,68){\vector(0,-1){8}}\put(30,60){\vector(0,-1){8}}\put(30,52){\vector(0,-1){8}}
\put(30,44){\vector(1,2){8}}\put(38,60){\vector(1,-1){8}}\put(46,52){\vector(0,-1){8}}
\put(46,44){\vector(1,0){8}}
\put(60,68){\vector(0,-1){8}}\put(60,60){\vector(0,-1){8}}\put(60,52){\vector(0,-1){8}}
\put(60,44){\vector(1,2){8}}\put(68,60){\vector(0,-1){16}}\put(68,44){\vector(1,0){8}}
\put(76,44){\vector(1,0){8}}

\put(8,4){$\mu=3$}\put(38,4){$\mu=3$}\put(68,4){$\mu=3$}
\put(8,38){$\mu=6$}\put(38,38){$\mu=3$}\put(68,38){$\mu=6$}

\end{picture}
\caption{Enumeration of cuspidal cubics}\label{fex10}
\end{figure}

\begin{theorem}\label{tex}
For any fixed $r\ge 1$ and $d\ge2r+3$, there exists a generic real $2r$-dimensional linear subsystem in $|{\mathcal O}_{\PP^2}(d)|$ that contains at least $c_r(d)$ real curves with $r$ real cusps as their only singularities, where $c_r(d)$ is a positive function satisfying
$$c_r(d)=\frac{(3d^2)^r}{r!}+O(d^{2r-1})\quad\text{as}\quad d\to\infty\ .$$
\end{theorem}

\begin{remark}
Note that $c_r(d)$ differs (asymptotically) by a constant factor from the total number $\frac{1}{r!}(12d^2)^r+O(d^{2r-1})$ of complex $r$-cuspidal curves in a generic
$2r$-dimensional linear subsystem of $|{\mathcal O}_{\PP^2}(d)|$. The asymptotic formula for the count of complex $r$-cuspidal curves was communicated to us by D.Kerner \cite{Ker}. We shortly sketch the idea of his proof. Since $\deg V_d(A_2)=12(d-1)(d-2)$, the lift of $V_d(A_2)\subset\PP^{d(d+3)/2}$ to the incidence variety $\widetilde V_d(A_2)\subset \PP^{d(d+3)/2}\times\PP^2$ represents the cohomology class $$12(d-1)(d-2)\alpha^2\beta^2+a\alpha^3\beta+b\alpha^4,$$ where $\alpha\in H^2(\PP^{d(d+3)/2},\Z)$, $\beta\in H^2(\PP^2,\Z)$ are the standard generators. Correspondingly, the family of $r$-cuspidal curves of degree $d$ with numbered singular points (which is an $r!$-sheeted covering of $V_d(rA_2)\subset\PP^{d(d+3)/2}$) lifts to a variety in $\PP^{d(d+3)/2}\times(\PP^2)^r$ that represents the cohomology class
$$\prod_{i=1}^r(12(d-1)(d-2)\alpha^2\beta_i^2+a\alpha^3\beta_i+b\alpha^4)$$ corrected by reducing parasitic diagonal terms. In this cohomology class we are only interested in the coefficient of the monomial $\alpha^{2r}\beta_1^2...\beta_r^2$. It remains to observe that the correction coming from the parasitic diagonal terms amounts to $O(d^{2r-2})$. Thus, the aforementioned asymptotic follows. We notice also that a similar statement for nodal singularities can be deduced from the known formulas for degrees of Severi varieties \cite{Qv}.
\end{remark}

{\bf Proof.} {\bf(1)} Suppose that $r=1$. Introduce a configuration of $n=\frac{d(d+3)}{2}-2$ points \begin{equation}\bw=\{(t^{-M_i},t^{M_i\eps})\}_{i=1,...,n},\quad t>0\ ,\label{econfig1}\end{equation}
which tropicalizes to configuration
(\ref{econfig}).
We intend to show that the two-dimensional linear system of curves of degree $d$ passing through $\bw$ contains at least $3d^3+O(d)$ real cuspidal curves. To this end, we go through the computations in the proof of Lemma \ref{l-unicusp} and check how many real solutions occur among all complex ones. Since we are interested only in the leading term of the asymptotics, we focus only on the subdivisions as shown in Figures \ref{fex2}(a,b,c,i,j). The question reduces to counting real solutions in the computation of limit curves associated with a given subdivision and of limit curves associated with modifications, and in equations for conditions to pass through fixed points $\bw$.

Observe that all the limit curves associated with the considered subdivisions are real. Indeed, up to a constant factor, all the coefficients $a_\omega$, $\omega\in{\mathcal T}_d\cap\Z^2$ can be found from the condition to pass through the configuration $\bw$ (see details in \cite[Section 3.7, formula (3.7.27)]{Shustin2005} and computations in \cite[Section 2.5.7]{IMS}) and from the known structure of the limit curves associated with triangles and quadrangles, and all these conditions reduce to linear equations.

Now, we case by case analyze modifications of the tropical limit and modifications of the fixed point conditions. First, the coordinates of the fixed points (\ref{econfig1}) determine the coefficients $a_{ij}$ of the polynomials defining the limit curves at the points $(i,j)$ on the lattice path so that
the signs of $a_{ij}$ alternate.
\begin{itemize}\item Consider the case shown in Figure \ref{fex2}(a). An equation for the modified limit curve exposed in Lemma \ref{ledge2} is obtained after a complex coordinate change. In the real setting, we obtain an equation
$$y^2-2yf(x)+\frac{a_{i+1,d-i-1}}{a_{i-1,0}}=y^2-2yf(x)+(-1)^{d-i}=0\ .$$ The condition to have a cusp yields two real cubic polynomials $f(x)$ iff $d-i$ is even (cf. the first paragraph of the proof of Lemma \ref{ledge2}). The condition to pass through a point of $\bw$ amounts to taking the cubic root (for the cuspidal fragment of type B see equations (\ref{efixed}), for other instances see \cite[Equation (2.19) in Section 2.5.9]{IMS}). Hence, the number of real solutions is the sixth part of the number of complex ones, i.e., $\frac{1}{2}d^2+O(d)$.
\item In the case shown in Figures \ref{fex2}(b,c), a non-linear algebraic equation pops up only
in computing the limit curve of the modification associated with a quadrangle as shown in Figure
\ref{ftrap1}(f) with $m=2$, $m_1=m_2=1$: here we again take the cubic root, and hence the third part of all complex solutions appears to be real, that is, $d^2+O(d)$.
\item Similarly to the previous item, in the cases shown in Figure \ref{fex2}(i,j), the third part of the limit curves in the modification is real. An analysis of the fixed point conditions is based on equations (\ref{efixed}), where the first two equations correspond to the marked cuspidal fragments depicted in Figure \ref{ftrap1}(j,k), while the last equation corresponds to the marked cuspidal fragment depicted in Figure \ref{ftrap1}(l). The images of the unmarked cuspidal fragment and of its tropical Mikhalkin's modification, see, \cite[Section 5.1]{MR} (dual to the subdivision shown in
    Figure \ref{ftrap1}(f)) are depicted in Figure \ref{fex2}(k,l). Notice that the point $v$ divides the segment $[v_3,v_2]$ as $1:2$.
    Restriction (\ref{econfig}) to the exponents $M_i$, $i=1,...,n$, in the coordinates of the points
    (\ref{econfig1}) implies that the image of the marked point in the segment
    is located much closer to the vertex $v_2$ than to $v_3$ (in the notation of Figure \ref{ftrap1}(k,l)) in case of Figure \ref{fex2}(i), and vice versa in case of Figure \ref{fex2}(j). Thus, in the latter case, we are in a position shown in Figure \ref{ftrap1}(j,k), where both the fixed point conditions are linear, and hence here we encounter $d^2+O(d)$ real solutions.
In case of Figure \ref{fex2}(i), the fixed point condition amounts to solving the last equation in (\ref{efixed}), which in the considered case reads
    \begin{equation}a_{i,j+1}(\eta(t)-\zeta)^2+a_{i-1,0}t^{\delta}+\text{h.o.t.}=0\ ,\label{econfig2}\end{equation}
    with given $\zeta$ and $\delta$ and the unknown $\eta(t)$. The sign alternation rule pointed above yields that $a_{i-1,0}/a_{i,j+1}=(-1)^{d-i-j}$. It means that asymptotically in half of the cases as shown in Figure \ref{fex2}(i) the fixed point condition results in two real solutions, whereas in the other half in a pair of complex conjugate solutions. This finally gives us $\qquad$
    $\frac{1}{2}d^2+O(d)$ real modified tropical limits.
\end{itemize}
Summing up the results of the above computation, we obtain the lower bound $3d^2+O(d)$ for the number of the real unicuspidal curves in the constructed two-dimensional linear system.

\smallskip
{\bf(2)} Now we extend the bound obtained for $r=1$ to the case of $r\ge2$ in the way proposed in \cite{MMSS}. Introduce a configuration of $n=\frac{d(d+3)}{2}-2r$ points $$\bw=\{(t^{-M_i},t^{M_i\eps})\}_{i=1,...,n},\quad t>0\ ,$$
which tropicalize to configuration (\ref{econfig}). Then choose a sequence of integers $i_1,...,i_r$ between $1$ and $d-2$ so that
\begin{equation}|i_{s'}-i_{s''}|\ge2\quad\text{for all}\quad 1\le s'<s''\le r\ .\label{econfig3}\end{equation} For a given $i_s$, in the part of the Newton triangle
${\mathcal T}_d$ cut off by the strip $i_s-1\le x\le i_s+1$ we define the lattice path chosen arbitrarily among the
patterns shown in Figure \ref{fex2}(a,b,c,i,j). Observe that the choice of the lattice path in one strip does not affect the choice in any other strip. In the components of the complement to the strips, we define the lattice
$\lambda$-path going through all integral points. Since the elements of the subdivision of ${\mathcal T}_d$ in different strips, which are dual to the cuspidal fragments of the underlying tropical curve, are pairwise disjoint, we can smoothly extend the patchworking Theorem \ref{t2} to the considered situation.

The computations in Step 1 of the proof applied to the lattice paths described in the previous paragraph allow one to count all possible real modified tropical limits emerging from these lattice paths for each sequence $1\le i_1,...,i_r\le d-3$ subject to restrictions (\ref{econfig2}). Namely, when $i_1$ runs from $1$ to $d-2$, the fragments of the lattice path in the strips $i_1-1\le x\le i_1+1$ yield at least $3d^2+O(d)$ real modified tropical limits. Then we inductively allow $i_s$ to run from $1$ to $d-2$ obeying the restrictions $|i_s-i_{s'}|\ge2$ for all $s'<s$, as $s=2,...,r$. Notice that the number of real modified tropical limits produced in the fragments of the lattice path in all the strips $i_s-1\le x\le i_s+1$ differs from the maximal value $3d^2+O(d)$ by at most $6(s-1)d$. Thus, summing up over all sequences $i_1,...,i_r$ and factorizing by permutations, we obtain the bound asserted in Step 2 of the theorem (cf. the same type of computations in \cite[Section 6]{MMSS}).
\proofend

\section*{Appendix}

\subsection*{A1. Multiplicity of a flat cycle (fragment D)}
{\bf (1)} Let $(\Gamma,\bp,h)$ have a cuspidal fragment of type D (flat cycle) $(\Gamma_c,\bp_c,h_c)$ consisting of a four-valent vertex $V_1$, a flat trivalent vertex and two edges joining them
(see Figure \ref{fig-cfd}(a), where $m_1,m_2,m$ denote
the weights of the edges incident to the flat trivalent vertex). Since $m\ge2$ and $\Delta$ is primitive,
all the edges incident to the flat trivalent vertex are bounded, and we extend $(\Gamma_c,\bp_c,h_c)$
to a fragment $(\Gamma',\bp',h')$
by adding one more trivalent vertex $V_2$ as shown in Figure \ref{fig-cfd}(b). The fragment of the
subdivision of $P(\Delta)$ dual to $h'_+(\Gamma')$ consists of two triangles $T_1,T_2$
sharing a common side $e$
(see Figure \ref{fig-cfd}(c)). Denote by $V_{1,red}$ the vertex of $h'_*(\Gamma')$ dual to the
triangle $T_1$. The limit curves $C_1\subset\Tor(T_1)$, $C_2\subset\Tor(T_2)$ are rational,
$C_2$ is peripherally smooth and unibranch, while $C_1$ has two local branches centered at the same point of
$\Tor(e)$ and it is unibranch along the two other toric divisors. The following
lemma describes the geometry of $C_1$.

\begin{figure}
\setlength{\unitlength}{1.0mm}
\begin{picture}(130,150)(0,0)
\thicklines
\put(0,130){\line(1,1){10}}\put(0,150){\line(1,-1){10}}
\put(10,141){\line(1,0){20}}\put(10,139){\line(1,0){20}}
\put(30,140){\line(1,0){10}}

\put(50,130){\line(1,1){10}}\put(50,150){\line(1,-1){10}}
\put(60,141){\line(1,0){20}}\put(60,139){\line(1,0){20}}
\put(80,140){\line(1,0){10}}\put(90,140){\line(1,1){10}}\put(90,140){\line(1,-1){10}}

\put(110,140){\line(1,1){10}}\put(110,140){\line(1,-1){10}}
\put(120,130){\line(0,1){20}}\put(120,130){\line(1,1){10}}\put(120,150){\line(1,-1){10}}

\put(5,75){\line(1,1){10}}\put(5,75){\line(1,0){30}}
\put(15,85){\line(2,3){10}}\put(25,100){\line(2,-5){10}}

\put(50,75){\line(4,5){20}}\put(70,100){\line(2,-5){10}}
\put(50,75){\line(1,0){30}}

\put(95,75){\line(1,1){10}}\put(95,75){\line(1,0){30}}
\put(105,85){\line(2,3){10}}\put(115,100){\line(2,-5){10}}
\put(105,85){\line(1,-1){10}}

\put(20,20){\vector(-1,-1){10}}\put(10,30){\vector(1,-1){10}}
\put(35,21){\vector(-1,0){15}}\put(35,19){\vector(-1,0){15}}
\put(35,20){\vector(1,0){15}}\put(50,20){\vector(-1,0){15}}
\put(50,20){\vector(1,-1){10}}\put(60,30){\vector(-1,-1){10}}

\put(10,35){\vector(1,1){10}}\put(10,55){\vector(1,-1){10}}
\put(20,46){\vector(1,0){15}}\put(20,44){\vector(1,0){15}}
\put(35,45){\vector(1,0){15}}\put(50,45){\vector(1,-1){10}}
\put(60,55){\vector(-1,-1){10}}

\put(90,45){\vector(-1,-1){10}}\put(80,55){\vector(1,-1){10}}
\put(105,46){\vector(-1,0){15}}\put(105,44){\vector(-1,0){15}}
\put(120,45){\vector(-1,0){15}}\put(130,35){\vector(-1,1){10}}
\put(130,55){\vector(-1,-1){10}}

\put(90,20){\vector(-1,-1){10}}\put(80,30){\vector(1,-1){10}}
\put(90,21){\vector(1,0){15}}\put(90,19){\vector(1,0){15}}
\put(105,21){\vector(-1,0){15}}\put(105,19){\vector(-1,0){15}}
\put(105,20){\vector(1,0){15}}\put(120,20){\vector(1,-1){10}}
\put(130,30){\vector(-1,-1){10}}

\thinlines
\put(5,75){\vector(0,1){35}}\put(5,75){\vector(1,0){35}}
\put(50,75){\vector(0,1){35}}\put(50,75){\vector(1,0){35}}
\put(95,75){\vector(0,1){35}}\put(95,75){\vector(1,0){35}}

\dottedline{1}(15,75)(15,85)\dottedline{1}(25,75)(25,100)
\dottedline{1}(5,85)(15,85)\dottedline{1}(5,100)(25,100)

\dottedline{1}(70,75)(70,100)\dottedline{1}(50,100)(70,100)

\dottedline{1}(105,75)(105,85)\dottedline{1}(95,85)(105,85)
\dottedline{1}(115,75)(115,100)\dottedline{1}(95,100)(115,100)

\put(9,139){$\bullet$}\put(29,139){$\bullet$}
\put(59,139){$\bullet$}\put(79,139){$\bullet$}
\put(19,44){$\bullet$}\put(34,44){$\bullet$}
\put(89,44){$\bullet$}\put(104,44){$\bullet$}
\put(19,19){$\bullet$}\put(34,19){$\bullet$}
\put(89,19){$\bullet$}\put(104,19){$\bullet$}

\put(39,18.7){\Large$*$}\put(97,17.7){\Large$*$}\put(97,19.7){\Large$*$}

\put(10,143){$V_1$}\put(60,143){$V_1$}\put(87,142){$V_2$}
\put(115,139){$T_1$}\put(122,139){$T_2$}
\put(15,102){$Q$}\put(63,102){$T_0$}\put(1,99){$m$}\put(0,84){$m_1$}
\put(46,99){$m$}\put(91,99){$m$}\put(90,84){$m_1$}
\put(4,71){$0$}\put(14,71){$1$}\put(24,71){$2$}\put(34,71){$3$}
\put(49,71){$0$}\put(69,71){$2$}\put(79,71){$3$}
\put(94,71){$0$}\put(104,71){$1$}\put(114,71){$2$}\put(124,71){$3$}
\put(103,78){$Q'$}\put(113,85){$Q''$}
\put(19,120){(a)}\put(69,120){(b)}\put(118,120){(c)}
\put(19,65){(d)}\put(64,65){(e)}\put(109,65){(f)}\put(67,0){(g)}
\put(18,136){$m_1$}\put(18,143){$m_2$}\put(34,142){$m$}

\end{picture}
\caption{Cuspidal flat cycle and its modification, I}\label{fig-cfd}
\end{figure}

\medskip
\noindent
{\bf Lemma A1.} 
{\it Let ${\mathcal T}$ be a nondegenerate lattice triangle with sides $\sigma_1,\sigma_2,\sigma_3$. Suppose that
$m=\|\sigma_1\|_\Z=m_1+m_2$, $m_2\ge m_1\ge1$.
Let $M\subset{\mathcal M}_{0,4}(\Tor({\mathcal T}),{\mathcal L}_{\mathcal T})$ be the family of
isomorphism classes of
maps $\bn:\PP^1\to\Tor({\mathcal T})$ of $\PP^1$ with
four distinct marked points $p_1,p'_1,p_2,p_3\in\PP^1$ such that
$$\bn(p_1)=\bn(p'_1)\in\Tor^*(\sigma_1),\quad \bn(p_i)
\in\C^*\subset\Tor(\sigma_i),\ i=2,3\ ,$$
and
$$\bn^*\Tor(\sigma_1)=m_1p_1+m_2p'_1,\quad\bn^*\Tor(\sigma_i)=
\|\sigma_i\|_\Z\cdot p_i,\ i=2,3\ .$$
The following holds:

(1) If $\|{\mathcal T}\|_\Z=m$ then $M=\emptyset$; if
$\|{\mathcal T}\|_\Z=rm$, $r\ge2$, then $M$ is isomorphic to the union of disjoint copies of $(\C^*)^2$ that are in one-to-one correspondence with the set
\begin{equation}\begin{cases}\Omega_r:=\{\lambda\in\C\ :\ \lambda^r=1,\ \lambda\ne1\}\quad&\text{if}\ m_1<m_2,\\
\Omega_r/\conj\quad&\text{if}\ m_1=m_2=\frac{m}{2}.\end{cases}\label{ecfd1}\end{equation}
Except for the case
\begin{equation}\|\sigma_i\|_\Z\equiv0\mod2,\ i=1,2,3, \quad\text{and}\quad m_1=m_2=\frac{m}{2}\ ,
\label{en-ex}\end{equation} each irreducible component of $M$ parameterizes birational maps onto immersed curves. 
In the case (\ref{en-ex}), the components of $M$ associated with the pairs $(\lambda,\overline\lambda)\in\Omega_r/\conj$, $\lambda\ne-1$, parameterize birational maps as above,
while the component associated with $\lambda=-1$ parameterizes double ramified coverings $\bn:\PP^1\to C'\hookrightarrow\Tor({\mathcal T})$ with
$2C'\in|{\mathcal L}_{\mathcal T}|$ and $C'$ immersed, unibranch along the toric divisors,
and the ramification points located at
$C'\cap\Tor(\sigma_i)$, $i=2,3$.

(2) Under the condition
$\|{\mathcal T}\|_\Z\ge 2m$, given two points $z_i\in\Tor^*(\sigma_i)$, $z_j\in\Tor^*(\sigma_j)$,
$1\le i<j\le 3$, each component of $M$ transversally intersects in
${\mathcal M}_{0,4}(\Tor({\mathcal T}),{\mathcal L}_{\mathcal T})$
with the family $\{\bn(p_i)=z_i,\ \bn(p_j)=z_j\}$ in
\begin{equation}N=\begin{cases}\frac{1}{2}\cdot\|{\mathcal T}\|_\Z\big(\|\sigma_i\|_\Z\cdot\|\sigma_j\|_\Z\big)^{-1},\quad&\text{if}\ m_1=m_2,\ \lambda=-1,
\ \text{and}\\
&\|\sigma_2\|_\Z\equiv\|\sigma_3\|_\Z\equiv1\mod2,\\
\|{\mathcal T}\|_\Z\big(\|\sigma_i\|_\Z\cdot\|\sigma_j\|_\Z\big)^{-1},\quad&\text{otherwise}
\end{cases} \label{en-ex4}\end{equation}
distinct points.}
\medskip

{\bf Proof.} Applying a suitable automorphism of $\Z^2$, we can identify ${\mathcal T}$ with the triangle
$\conv\{(p,0),(p+m,0),(0,r)\}$, where $p\ge0$, $r>0$, $p+m>r$, and $\sigma_1=[(p,0),(p+m,0)]$. If $\|{\mathcal T}\|_\Z=m$, i.e.,
$r=1$, then $\bn^*\Tor(\sigma_1)$ is one point, and hence $M=\emptyset$.

Suppose that $r\ge2$. Then a map $\bn:\PP^1\to\Tor({\mathcal T})$ as asserted in Lemma can be given by
\begin{equation}x=a t^r,\quad y=b t^p(t-1)^{m_1}(t-\lambda)^{m_2},\quad t\in\C\ ,
\label{en-ex3}\end{equation} with some $\lambda\in\C\setminus\{0,1\}$ and arbitrary
$a,b\in\C^*$. Since $x(1)=x(\lambda)$, we get $\lambda=\exp\frac{2\pi k\sqrt{-1}}{r}$ with $1\le k\le r-1$, and for each value of
$k$, the family of such parameterizations is isomorphic to $(\C^*)^2$.
If $m_1<m_2$, then the parameterizations are in one-to-one correspondence with the elements of $M$.
If $m_1=m_2$ and $\lambda\ne-1$, then the parameterizations associated with the data
$(a,b,\lambda)$ and $(a,b\lambda^{p+m},\overline\lambda)$ define the same element of $M$. If $m_1=m_2$ and $\lambda=-1$, but $p$ is odd, then the parameterizations associated with the data
$(a,b,-1)$ and $(a,-b,-1)$ define the same element of $M$, and here $(\C^*)^2/\{(a,b)\sim(a,-b)\}\simeq
(\C^*)^2$. Hence, the formula (\ref{ecfd1}).

If either (\ref{en-ex}) does not hold, or (\ref{en-ex}) holds but $\lambda\ne-1$, then by (\ref{en-ex3})
the map $\bn$ is birational onto its image. If
(\ref{en-ex}) holds (in particular, $r,p,m$ are even) and $\lambda=-1$, then formula (\ref{en-ex3}) turns into
$$x=a (t^2)^{r/2},\quad y=b (t^2)^{p/2}(t^2-1)^{m/2},\quad t\in\C\ ,$$ which yields the double covering as asserted in the lemma.

For the second part of the lemma, we introduce
$$d_1=\gcd(p,r),\ p'=\frac{p}{d_1},\ r'=\frac{r}{d_1},\quad d_2=\gcd(p+m,r),\ p''=\frac{p}{d_2},\ r''=\frac{r}{d_2}\ .$$
Then the conditions $\bn(p_i)=z_i$, $\bn(p_j)=z_j$ amount to the following systems
of equations in unknowns $a,b$
$$\begin{cases}&a=\alpha,\\ &a^{p'}=\beta b^{r'},\end{cases}\qquad\begin{cases}&a=\alpha,\\ &a^{p''}=\beta b^{r''},
\end{cases}\qquad\begin{cases}&a^{p'}=\alpha b^{r'},\\
&a^{p''}=\beta b^{r''},\end{cases}$$ for $(i,j)=(1,2)$, $(1,3)$, or $(2,3)$, respectively, with some
$\alpha,\beta\in\C^*$. In each case, the number of solutions $(a,b)$ is given by
the second value in (\ref{en-ex4}), while in the case $m_1=m_2$, $\lambda=-1$, and $p\equiv1\mod2$, one
has to identify solutions $(a,b)$ and $(a,-b)$.
The transversality of the intersection follows from the fact that each of the above systems has only simple solutions
in $(\C^*)^2$.
\proofend

We say that a nondegenerate lattice triangle $T$ is {\it special} with respect to its side $\sigma$ if
\begin{itemize}\item $m:=\|\sigma\|_\Z$ is even;
\item by an automorphism of $\Z^2$ the triangle $T$ can be brought to the form
$\conv\{(0,r),(p,0),(p+m,0\}$, where the segment $[(p,0),(p+m,0)]$ is the image of the side $\sigma$;
\item it holds 
\begin{equation}\Omega_{r,p,m}:=\left\{\lambda\in\Omega_r\setminus\{-1\}\ :\ \lambda^{p+m/2}=(-1)^{m/2}\right\}\ne\emptyset\ .\label{especial}\end{equation}
\end{itemize}
Note that the parameter $p$ is uniquely defined modulo $r$, and hence the set $\Omega_{r,p,m}$ does not depend on the choice of the form of the triangle $T$ in the second condition.

\medskip
{\bf Lemma A2.} {\it Let $(\Gamma,\bp,h)$ be a plane cuspidal $n$-marked
tropical curve of a nondegenerate, primitive degree $\Delta$ and genus $g$, where $n=|\Delta|+g-2$,
$h(\bp)=\bx$ a configuration of $n$ distinct points in $\R^2$ in general position, and $\bx=\val(\bw)$ with $\bw$ a configuration of $n$ points in $(\K^*)^2$. Suppose that $(\Gamma,\bp,h)$ has a cuspidal tropical
fragment $(\Gamma_c,\bp_c,h_c)$ of type D and that the triangle $T_1$ is not special with respect to
its side $e=T_1\cap T_2$ (see Figure \ref{fig-cfd}(c)). Then the number of curves $C\in V_{\Delta,g}(A_2)$, passing through $\bw$ and tropicalizing into $(\Gamma,\bp,h)$ equals (cf. (\ref{epw})
$$\mu_c(\Gamma_c,\bp_c,h_c)\cdot\prod_{V\in\Gamma^0\setminus\Gamma_c}\mu(V)\ ,$$
where the first factor is as follows (in the notations and conventions of Lemma A1):
\begin{enumerate}\item[(i)] if $m_1<m_2$, then
\begin{equation}\mu_c(\Gamma_c,\bp_c,h_c)=\frac{\mu(V_{1,red})(\mu(V_{1,red})-m)(m+m_2)m_1}{m^2}\ ;
\label{ecfd2}\end{equation}
\item[(ii)] if $m_1=m_2$ and two sides of $T_1$ have odd length 
then \begin{equation}\mu_c(\Gamma_c,\bp_c,h_c)=\frac{3}{8}\mu(V_{1,red})(\mu(V_{1,red})-m)\ ;\label{ecfd1x}\end{equation}
\item[(iii)] if $m_1=m_2$ and all three sides of $T_1$ have even length, then
\begin{equation}\mu_c(\Gamma_c,\bp_c,h_c)=\frac{3}{4}\mu(V_{1,red})(\mu(V_{1,red})-m-2)\ .\label{ecfd2x}\end{equation}
\end{enumerate}}
\medskip

{\bf Proof.} We, first, compute the number of refined tropical limits associated with a given cuspidal tropical curve with a cuspidal fragment of type D and with the given configuration of points $\bw\subset(\K^*)^2$. Then we prove a patchworking statement that to each of the constructed tropical limits assigns a curve $C\in V_{\Delta,g}(A_2)$.

\smallskip

{\bf(1)} Let us perform a local modification
at the point $z=C_1\cap C_2$ along the edge $e=T_1\cap T_2$ following the recipe of Section \ref{nod-mod},
and we obtain a fragment of the modified subdivision inscribed either into the non-convex
quadrangle $Q$ with vertices $(0,0),(3,0),(1,m_1),(2,m)$ if $m_1<m_2$, or into the triangle $T=\conv\{(0,0),
(3,0),(2,m)\}$ if $m_1=m_2=\frac{m}{2}$ (see Figure \ref{fig-cfd}(d,e)).

\smallskip
{\bf(1i)} Suppose that $m_1<m_2$. It follows from Theorem \ref{tuc1} that the union of the limit curves
corresponding to the subdivision of $Q$ induced by the modified polynomial must be of arithmetic genus zero
and it must contain a local singular branch. Getting rid of the monomial $xy^{m_1-1}$ in the modified polynomial and applying the argument used in Section \ref{sec-cfb}, we derive that the only possible subdivision of $Q$ is as shown in Figure \ref{fig-cfd}(f), while the limit curve $C'\subset
\Tor(Q')$ is nodal, rational, and $C''\subset\Tor(Q'')$ is rational with a unique singular branch in
$\Tor^*(Q'')$. Observe that the curves $C_1,C_2$ determine the intersection points
of $C',C''$ with the toric divisors corresponding to the inclined sides of $Q$. Taking into account
possible orientations on $\Gamma'\setminus\bp'$ induced by the regular orientation of
$\Gamma'$ in case $\bp_c=\emptyset$ (two upper graphs in Figure \ref{fig-cfd}(g)) and applying
\cite[Lemma 3.9]{Shustin2005} and Lemma A1(1), we obtain that
the number of the tuples $(C_1,C_2,C',C'')$ matching the initial data equals
$$\mu(V_2)\widetilde m^{-1}(r-1)\cdot\frac{\mu(V_{1,red})}{m}\cdot\left|\begin{matrix}
-1 & -m_2\\ 1 & m\end{matrix}\right|\cdot m_1$$
\begin{equation}=\mu(V_2)\widetilde m^{-1}\mu_c(\Gamma_c,\bp_c,h_c)\ ,\label{ecfd3}\end{equation} where $\widetilde m$ is the product of the weights of the outer edges of $(\Gamma',\bp',h')$ oriented
inward, and $\mu_c(\Gamma_c,\bp_c,h_c)$ is given by (\ref{ecfd2}).
If $\bp_c\ne\emptyset$, i.e., there is a marked point on one of the edges incident to
the flat trivalent vertex (see two lower graphs in Figure \ref{fig-cfd}(g)), then we take into account
the marked point condition as in \cite[Section 2.3.9]{IMS} and also sum up the multiplicities
arising when the marked point lies on the upper and lower edge of the flat cycle. Then the final result coincides with (\ref{ecfd2}).

\smallskip
\smallskip{\bf(1ii)} Suppose now that $m_1=m_2=\frac{m}{2}$ and two sides of $T_1$ have odd length.
Then $C_1$ is reduced (see Lemma A1). In addition, the two local branches of $C_1$ centered at a point of $\Tor(e)$ intersect each other with multiplicity $\frac{m}{2}$. Indeed,
the corresponding map germs
$\bn:(\C,1)\to\Tor({\mathcal T})$ and $\bn:(\C,\lambda)\to\Tor({\mathcal T})$ are given by
\begin{equation}\begin{cases}x=a+ar(t-1)+O((t-1)^2),&\\
y=b(1-\lambda)^{m/2}(t-1)^{m/2}+O((t-1)^{m/2+1})&\end{cases}
\quad\text{with}\ t-1\in(\C,0)\ ,\label{en-ex1}\end{equation} and
\begin{equation}\begin{cases}x=a+ar\lambda^{r-1}(t-\lambda)+O((t-\lambda)^2),&\\
y=b\lambda^p(\lambda-1)^{m/2}(t-\lambda)^{m/2}+O((t-\lambda)^{m/2+1})&\end{cases}
\quad\text{with}\ t-\lambda\in(\C,0)\ ,\label{en-ex2}\end{equation}
which after elimination of $t$ turns into
$$y=b\left(\frac{1-\lambda}{ra}\right)^{m/2}(x-a)^{m/2}+O((x-a)^{m/2+1})\quad\text{and}$$
$$y=b\lambda^p\left(\frac{\lambda(\lambda-1)}{ra}\right)^{m/2}(x-a)^{m/2}+O((x-a)^{m/2+1})
\ .$$
The equality of the coefficients of $(x-a)^{m/2}$ in the above two formulas yields either relation (\ref{especial})
for $\lambda\ne\pm1$, or the relation $r\equiv p\equiv0\mod2$ for $\lambda=-1$, the latter condition being equivalent to (\ref{en-ex}).
Since both the conditions (\ref{especial}) and (\ref{en-ex}) are excluded, we get that the two considered local branches of $C_1$ intersect with multiplicity $m/2$.

As we noticed in the proof of Lemma A1, the refinement of the tropical limit at the intersection point $C_1\cap C_2$ yields the triangle $T_0$ (see Figure \ref{fig-cfd}(e)). The curves $C_1,C_2$ impose the three preassigned intersection points of the limit curve $C_0$ of the refinement with the toric divisors $\Tor(e_1),\Tor(e_2)$, where $e_1=[(0,0),(2,m)]$, $e_2=[(2,m),(3,0)]$:
\begin{itemize}\item one point on the toric divisor
$\Tor(e_2)$, \item two distinct points on the toric divisor
$\Tor(e_1)$. \end{itemize} The limit curve $C_0$ must be of arithmetic genus zero and it intersects the toric divisors wih multiplicity $1$ at each of the above three preassigned points. Furthermore, it must have a singular local branch (outside the toric divisors $\Tor(e_1),\Tor(e_2)$).

Getting rid of the monomial $x^2y^{m-1}$ by an extra coordinate change (cf. Section \ref{sec-cfb}), we exclude the segment $[(2,0),(2,m)]$ as
an element of a possible subdivision of $T$. Since the induced convex piecewise linear function is linear along the edge $e_1$, its midpoint cannot be a vertex of a subdivision of $T_0$.
It follows that no any subdivision of $T_0$ into proper subpolygons is possible, when requiring limit curves with a cusp. Thus, $T_0$ is an entire piece of a subdivision. The corresponding limit curves $C_0\subset\Tor(T_0)$ are rational and admit a parametrization in the form
\begin{equation}x=\alpha t^{m/2}(t-1)^{m/2},\quad y=\beta\frac{S(t)}{t(t-1)},\quad \alpha,\beta\in\C^*\ ,
\label{euc2a}\end{equation} where $S(t)$ is a monic cubic polynomial. Since $x'(t)=0$ implies $t=\frac{1}{2}$, the singular branch of $C_0$ is unique, and by a shift $y\mapsto y+b$ we can place its center to the line $y=0$. This means $y(\frac{1}{2})=y'(\frac{1}{2})=0$, or, equivalently $S(t)=(t-\gamma)(t-\frac{1}{2})^2$. In a neighborhood of $\frac{1}{2}$, we have
$$x=m2^{m-3}\left(t-\frac{1}{2}\right)^2+\text{h.o.t.},\ y=(4\gamma-2)\left(t-\frac{1}{2}\right)^2-4\left(t-\frac{1}{2}\right)^3+\text{h.o.t.}$$
which yields the type $A_2$ of the singular branch. The point conditions along $\Tor(e_1)\cup\Tor(e_2)$ amount to a system of equations $$\alpha\beta^{m/2}\gamma^{m/2}=\xi_1,\quad\alpha\beta^{m/2}(1-\gamma)^{m/2}=\xi_2,\quad
\alpha=\xi_3\beta^m\ ,$$ where $\xi_1,\xi_2,\xi_3\ne0$, $\xi_1\ne\xi_2$. It follows that the system has $\frac{3}{4}m^2$ solutions $(\alpha,\beta,\gamma)$.

Assume that $r$ is odd. Then in view of Lemma A1 and the above computations, we obtain that the number of tuples $(C_1,C_2,C_0)$ matching the initial data equals
$$\mu(V_2)\widetilde m^{-1}\cdot\frac{r-1}{2}\cdot\frac{\mu(V_{1,red})}{m}\cdot\frac{3}{4}m^2=\mu(V_2)\widetilde m^{-1}\mu_c(\Gamma_c,\bp_c,h_c)\ ,$$
where $\widetilde m$ is the product of the weights of the outer edges of $(\Gamma',\bp',h')$ oriented
inward, and $\mu_c(\Gamma_c,\bp_c,h_c)$ is given by (\ref{ecfd1x}).

Assume that $r$ is even, but $p$ is odd. By Lemma A1, we have $\frac{r-2}{2}$ elements
$(\lambda,\overline\lambda)\in\Omega_r/\conj\setminus\{-1\}$ and thus,
$$\mu(V_2)\widetilde m^{-1}\cdot\frac{r-2}{2}\cdot\frac{\mu(V_{1,red})}{m}\cdot\frac{3}{4}m^2
=\mu(V_2)\widetilde m^{-1}\cdot\frac{3}{8}\mu(V_{1,red})(\mu(V_{1,red})-2m)$$ corresponding tuples $(C_1,C_2,C_0)$. For $\lambda=-1$, we get
$$\mu(V_2)\widetilde m^{-1}\cdot\frac{\mu(V_{1,red})}{2m}\cdot\frac{3}{4}m^2$$ tuples $(C_1,C_2,C_0)$.
Summing up, we obtain $\mu(V_2)\widetilde m^{-1}\cdot\mu_c(\Gamma_c,\bp_c,h_c)$ with the last factor given by (\ref{ecfd1x}).

\smallskip
{\bf(1iii)} Suppose now that $m_1=m_2$ and all the sides of $T_1$ have even length. If $\lambda\in\Omega_r\setminus\{-1\}$, then the argument in the preceding item yields that the number of the corresponding triples $(C_1,C_2,C_0)$ equals
$$\mu(V_2)\widetilde m^{-1}\cdot\frac{r-2}{2}\cdot\frac{\mu(V_{1,red})}{m}\cdot\frac{3}{4}m^2
=\mu(V_2)\widetilde m^{-1}\cdot\frac{3}{8}\mu(V_{1,red})(\mu(V_{1,red})-2m)\ .$$

If $\lambda=-1$, then $C_1$ is nonreduced (see Lemma A1(1)). That is, the curves $C_1,C_2$ determine one preassigned point on each of the toric divisors $\Tor(e_1)$, $\Tor(e_2)$.
As in the preceding item, we can assume that $C_0$ is rational with Newton triangle $T_0$.
Note that $C_0$ cannot be unibranch at the point on the toric divisor $\Tor(e_1)$. This can be showed by applying literally the argument from the proof of Lemma \ref{lmod}. Hence $C_0$ has two local branches at this point (transversally intersecting $\Tor(e_1)$), and a singular local branch centered outside $\Tor(e_1)\cup\Tor(e_2)$ (and which we can shift to the line $y=0$). Its parametrization is then as follows:
$$x=\alpha t^{m/2}(t-1)^{m/2},\quad y=\beta\frac{(t-1/2)^2(t-\gamma)}{t(t-1)}\ ,$$ and the point conditions along
$\Tor(e_1),\Tor(e_2)$ read
$$\alpha\beta^{m/2}(-\gamma)^{m/2}=\alpha\beta^{m/2}(1-\gamma)^{m/2}=\xi_1,\quad \alpha=\xi_2\beta^m$$ with some $\xi_1,\xi_2\in\C^*$. The system has $\frac{3}{4}m(m-2)$ solutions $(\alpha,\beta,\gamma)$, which due to the symmetry
$(\alpha,\beta,\gamma,t)\leftrightarrow(\alpha,-\beta,1-\gamma,1-t)$ define $\frac{3}{8}m(m-2)$ curves $C_0$.

Combining this with the computation in the preceding paragraph, we obtain formula (\ref{ecfd2x}).

\medskip
{\bf(2)} The proof of the patchworking statement, which recovers algebraic curves $C\in V_{\Delta,g}(A_2)$ passing through $\bw$ out of the refined tropical limits, literally follows the proof of Theorem \ref{t2}, which in turn is based on \cite[Theorem 2.4]{Sh} as the main ingredient. We only comment on one case that differs from those considered in Theorem \ref{t2}. It is the case mentioned in Lemma A1(1), where the limit curve $C_1$ is nonreduced, i.e., it is the image of a double covering
\begin{equation}\bn:\PP^1\to C'\hookrightarrow\Tor(T_1)\label{edouble}\end{equation}
ramified at two points on toric divisors of $\Tor(T_1)$, while $C_1=2C'$.
Then it is more convenient to work with the parameterized tropical limit. Furthermore, since the other limit curve are birational images of $\PP^1$, they can be treated as in the proof of Theorem \ref{t2}. For the parameterized limit curve
(\ref{edouble}), we only have to confirm the trasversality condition, but they coincide with the conditions of the transversality of intersections in Lemma A1(2).
\proofend

\medskip\noindent
{\bf Lemma A3.} {\it Let $(\Gamma,\bp,h)$ be a plane cuspidal $n$-marked
tropical curve of a nondegenerate, primitive degree $\Delta$ and genus $g$, where $n=|\Delta|+g-2$,
$h(\bp)=\bx$ a configuration of $n$ distinct points in $\R^2$ in general position, and $\bx=\val(\bw)$ with $\bw$ a configuration of $n$ points in $(\K^*)^2$. Suppose that $(\Gamma,\bp,h)$ has a cuspidal tropical
fragment $(\Gamma_c,\bp_c,h_c)$ of type D and that the triangle $T_1$ is special with respect to
its side $e=T_1\cap T_2$ (see Figure \ref{fig-cfd}(c)). Then the number of curves $C\in V_{\Delta,g}(A_2)$, passing through $\bw$ and tropicalizing into $(\Gamma,\bp,h)$ equals (cf. (\ref{epw}))
\begin{equation}\mu_c(\Gamma_c,\bp_c,h_c)\cdot\prod_{V\in\Gamma^0\setminus\Gamma_c}\mu(V)\ ,\label{eD}\end{equation}
where the first factor $\mu_c(\Gamma_c,\bp_c,h_c)$ is a function of the pair $(T_1,e)$ considered up to automorphisms of the lattice $\Z^2$.
}
\medskip

{\bf Proof.}
Without loss of generality, we can suppose that the fragment $(\Gamma',\bp',h')$ does not contain marked points, i.e., $\bp'=\emptyset$ (cf. item (1i) in the proof of Lemma A2).

The case of a nonreduced curve $C_1$ is covered in item (1iii) of the proof of Lemma A2. Thus, assume that $C_1$ is reduced. In this situation, the case of $\lambda\in\Omega_r\setminus\Omega_{r,p,m}$ is also covered in items (1ii) and (1iii) of the proof of Lemma A2.
Suppose that $\lambda\in\Omega_{r,p.m}$. We claim that the two local branches of $C_1$ centered at the point on $\Tor(e)$ intersect each other with multiplicity $\frac{m}{2}+1$. Indeed, let us refine formulas (\ref{en-ex1}), (\ref{en-ex2}) assuming for simplicity $\alpha=\beta=1$
(equivalently, making change of the coordinates $x=\alpha x'$, $y=\beta y'$):
$$\begin{cases}x=1+rt'+\frac{r(r-1)}{2}(t')^2+O((t')^3),&\\
y=(1-\lambda)^{m/2}(t')^{m/2}+(1-\lambda)^{m/2-1}(p+m/2-p\lambda)(t')^{m/2+1}&\\
\qquad\qquad+O((t')^{m/2+2}),&\end{cases}$$ where $t'=t-1$, and
$$\begin{cases}x=1+r\lambda^{r-1}t''+\frac{r(r-1)}{2}\lambda^{r-2}(t'')^2+O((t'')^3),&\\
y=\lambda^p(\lambda-1)^{m/2}(t'')^{m/2}&\\
\qquad+\lambda^{p-1}(\lambda-1)^{m/2-1}(\lambda(p+m/2)-p)(t'')^{m/2+1}+O((t'')^{m/2+2}\emph{}),&\end{cases}$$ where $t''=t-\lambda$. Eliminating $t'$ and $t''$, we respectively obtain
$$y=\left(\frac{1-\lambda}{r}\right)^{m/2}(x-1)^{m/2}+\frac{m(1-\lambda)^{m/2-1}}{4r^{m/2+1}}
(2-(r-1)(1-\lambda))(x-1)^{m/2+1}$$
$$+O((x-1)^{m/2+2})\ ,$$
$$y=\lambda^p\left(\frac{\lambda(\lambda-1)}{r}\right)^{m/2}(x-1)^{m/2}\qquad\qquad\qquad\qquad\qquad
\qquad\qquad\qquad\qquad\qquad$$
$$\qquad+\lambda^{p+m/2}\frac{(\lambda-1)^{m/2-1}}{r^{m/2+1}}\left(\left(p+\frac{m}{2}\right)
\lambda-p-\frac{m(\lambda-1)(r-1)}{4}\right)
(x-1)^{m/2+1}$$ $$\qquad+O((x-1)^{m/2+2})\ .$$ In view of $\lambda^{p+m/2}=(-1)^{m/2}$ the coefficients of $(x-1)^{m/2}$ coincide in both the expressions, but the equality of the coefficients of $(x-1)^{m/2+1}$ leads to the relation $\lambda=\frac{p-m/2}{p+m/2}$, which contradicts the assumption $\lambda\in\Omega_{r,p,m}$.

Recall that we look for algebraic curves over the non-Archimedean field $\K$, which equivalently can be expressed in the form of the diagram (see, for example, \cite[Section 2.3]{Shustin2005})
\begin{equation}\begin{matrix}{\mathcal C} & \hookrightarrow & {\mathfrak X}\\
\downarrow & & \downarrow\\ (\C,0) & = & (\C,0)\end{matrix}\label{efamilies}\end{equation}
where \begin{itemize}\item ${\mathfrak X}\to(\C,0)$ is a flat family of surfaces with the general fibre ${\mathfrak X}=\Tor(P(\Delta))$ and the central fibre ${\mathfrak X}_0$ being a reduced, reducible surface, whose components are toric surfaces associated with the polygons of the subdivision of $P(\Delta)$;\item ${\mathcal C}\to(\C,0)$ is an inscribed flat family of curves with the general fibre
${\mathcal C}_t\in V_{\Delta,g}(A_2)$ and the central fibre a reducible curve; in particular, $C_1$, $C_2$ are the components of ${\mathcal C}_0$ located in the components
$\Tor(T_1)$, $\Tor(T_2)$ of ${\mathfrak X}_0$, respectively.\end{itemize}
The singularity of ${\mathcal C}_0$ at the (unique) point $z$ on $\Tor(e)$ is, in fact,
a planar curve singularity combined of three smooth local branches so that two of them, ${\mathcal P}_1$ and ${\mathcal P}_2$, intersect with multiplicity $\frac{m}{2}+1$, and the third one ${\mathcal P}_3$ intersects each of ${\mathcal P}_1,{\mathcal P}_2$ with multiplicity $\frac{m}{2}$. The family ${\mathcal C}\to(\C,0)$ induces a one-parameter deformation of this singularity $({\mathcal C}_0,z)$, whose general member is represented by a map of a sphere with three holes to $\Tor(P(\Delta))$ such that its image has $\frac{3}{2}m-2$ nodes and one cusp. Note that, in the base of a versal deformation of $({\mathcal C}_0,z)$, the closure ${\mathcal V}$ of the locus consisting of the parameters that correspond to maps of a sphere with three holes to a neighborhood of $z$ with the image possessing $\frac{3}{2}m-2$ nodes and a cusp is the germ of an analytic subvariety of codimension $\frac{3}{2}m$.

On the other hand, the transversality conditions, imposed on the curves $C_1,C_2$ by the patchworking statement \cite[Theorem 2.4]{Sh} (see also \cite[Theorem 5]{Shustin2005}) allow one only a restricted class of deformations of the singularity $({\mathcal C}_0,z)$. Namely, assume that the regular orientation of $\Gamma'$ is as shown in the upper right picture in Figure \ref{fig-cfd}(g) The case of the orientation shown in he upper left picture can be treated in the same manner). The aforementioned transversality means the following:
\begin{itemize}\item the conditions imposed on $C_2$ by the rationality, by tangency at fixed points to the toric divisors in $\Tor(T_2)$ different from $\Tor(e)$, and by the mobile tangency to $\Tor(e)$ are independent (transversal),
\item the conditions imposed on $C_1$ by the rationality, one fixed and one mobile tangency to the toric divisors in $\Tor(T_1)$ different from $\Tor(e)$, and by intersecting $\Tor(e)$ with multiplicity $m$ at the fixed point $z\in\Tor(e)$ are independent (transversal). \end{itemize} These independence requirements hold, since, in terms of the deformation of maps $$\bn_1:\PP^1\to C_1\hookrightarrow\Tor(T_1),\quad \bn_2:\PP^1\to C_2\hookrightarrow\Tor(T_2),$$ they can be reduced to the vanishing $H^1(\PP^1,{\mathcal O}_{\PP^1}(-1))=0$. Furthermore, one can express the considered independence in terms of the coefficients of the defining polynomials. Namely, after the refinement at the point $z$ (see Figure \ref{fig-cfd}(e)), we have that the coefficients of the $\frac{3}{2}m$ monomials $xy^i$, $0\le i<\frac{m}{2}$, $x^2y^i$, $0\le i\le m+1$, $i\ne m-1$, corresponding to some integral points in $T_0$, can be varied independently, while the other coefficients (including those outside $T_0$) are analytic functions of the aforementioned coefficients and of the parameter $\tau$ in the families (\ref{efamilies}). That is, in the base of a versal deformation of the singularity $({\mathcal C}_0,z)$ (which we can suppose to contain the above coefficients as parameters), we have an intersection of the locus ${\mathcal V}$ of codimension $\frac{3}{2}m$ with an analytic subvariety of dimension $\frac{3}{2}m$. The number of intersections (in general, counting multiplicities) we denote $\Phi(r,p,m,\lambda)$, which is a function of the combinatorial parameters $r$, $p\mod r$, $m$ of the pair $(T_1,e)$ and of $\lambda\in\Omega_{r,p,m}/\conj$.

We cannot provide an explicit formula for $\Phi(r,p,m,\lambda)$, thus, leaving the factor $\mu_c(\Gamma_c,\bp_c,h_c)$ in (\ref{eD}) as an unknown function of the pair $(T_1,e)$.
\proofend

\subsection*{A2. Multiplicity of an elliptic trivalent vertex (fragment E)}
Suppose that the cuspidal fragment $(\Gamma_c,h_c)$ of $(\Gamma,\bp,h)$ is a non-flat trivalent
vertex of genus $1$.

\medskip\noindent
{\bf Lemma A4.} {\it
Let ${\mathcal T}$ be a nondegenerate lattice triangle,
$M\subset|{\mathcal L}_{\mathcal T}|$ the family of elliptic
curves that are unibranch along the toric divisors and have at least one singular local branch
in the big torus $\Tor^*({\mathcal T})$.

(1) If $|\Int({\mathcal T})\cap\Z^2|\le 1$ then $M=\emptyset$.

(2) If
$|\Int({\mathcal T})\cap\Z^2|\ge 2$, then $M$ is either empty, or is isomorphic to several
disjoint copies of $(\C^*)^2$ parameterizing curves, which are smooth along the toric divisors and which have exactly one singular local branch in
$\Tor^*({\mathcal T})$; furthermore, this singular branch has order $2$.

(3) Under the preceding assumption, choose two sides $\sigma_1,\sigma_2$ of ${\mathcal T}$ and fix points
$z_i\in\Tor^*(\sigma_i)$, $i=1,2$. Then $M$ transversally intersects with the linear system
$|{\mathcal L}_{\mathcal T}(-z_1-z_2)|\subset|{\mathcal L}_{\mathcal T}|$.}
\medskip

{\bf Proof.} Note that if $|\Int({\mathcal T})\cap\Z^2|\le 1$ then either there are no
elliptic curves in $|{\mathcal L}_{\mathcal T}|$, or all elliptic curves are smooth. Thus, $M=\emptyset$.

 Suppose that $|\Int({\mathcal T})\cap\Z^2|\ge2$ and $M\ne\emptyset$. In view of the $(\C^*)^2$-action, each component of
$M$ has dimension $\ge2$. Let us show that a generic curve $C$ of a component $M_0$ of $M$ has
in $\Tor({\mathcal T})$ either a unique singular branch
of order $\le3$, or two singular branches both of order $2$. Indeed, let $m_1,...,m_s\ge2$ ($s\ge1$) be orders of all singular branches
of $C$. Then we apply
\cite[Inequality (5)]{IKS}, which in our situation reads as follows (here we denote by $\sigma_1,\sigma_2,\sigma_3$ the sides of ${\mathcal T}$):
$$\sum_{i=1}^3\|\sigma_i\|_\Z\ge\sum_{i=1}^3(\|\sigma_i\|_\Z-1)+\sum_{i=1}^s(m_i-1)+1
\quad\Longrightarrow\quad\sum_{i=1}^s(m_i-1)\le2,$$
and hence the claim follows.

Now suppose that $C$ has a singular branch of order $3$ (which then must be in $\Tor^*({\mathcal T})$),
or two singular branches of order $2$ (and at least one of them in $\Tor^*({\mathcal T})$). Since $C\in M_0$ is generic, the
germ $(M_0,C)$ is equisingular, thus, is contained in the germ $(M_0^{ec},C)$ of the equiclassical
stratum (see details in \cite{DH,Shustin1998}). Then, in the same manner as in the proof of
Lemma \ref{ledge1} and additionally using \cite[Inequality (21)]{GuS} for the case of a singular branch on a toric divisor,
we conclude that the tangent cone to the considered equiclassical stratum can be identified with
$H^0(\widehat C,{\mathcal O}_{\widehat C}(D-D_0-D_{ec}))$, where
where $\bn:\widehat C\to C$ is the
normalization, $D=\bn^*(c_1({\mathcal O}_C\otimes{\mathcal L}_\Delta))$, $D_0$ is the double-point divisor,
$D_{ec}$ is the pull-back of the centers of singular local branches of $C$ and the intersection points of $C$ with the
toric divisors, counted with the total multiplicity
$$\deg D_{ec}=\sum_i(m_i-1)+\sum_{i=1}^3(\|\sigma_i\|_\Z-1)\ .$$
Since $$\deg D_0=2\sum_{z\in\Sing(C)}\delta(C,z)=C^2+CK_{\Tor({\mathcal T})}=\deg c_1({\mathcal O}_C\otimes{\mathcal L}_\Delta)
+K_{\Tor({\mathcal T})}C\ ,$$
$$\text{and}\quad\sum_{i=1}^3\|\sigma_i\|_\Z=-K_{\Tor({\mathcal T})}C\ ,$$ we have
$$\deg(D-D_0-D_{ec})=-K_{\Tor({\mathcal T})}C-\sum_i(m_i-1)-\sum_{i=1}^3(\|\sigma_i\|_\Z-1)$$
$$=3=\sum_{i=1}^s(m_i-1)>0=2g(\widehat C)-2\ ,$$ and hence
\begin{equation}H^1(\widehat C,{\mathcal O}_{\widehat C}(D-D_0-D_{ec}))=0\label{etrans}\end{equation}
and
$$h^0(\widehat C,{\mathcal O}_{\widehat C}(D-D_0-D_{ec}))=\deg(D-D_0-D_{ec})-g(\widehat C)+1=3-\sum_{i=1}(m_i-1)\ .$$
Due to the fact that the germ of $M$ at $C$ is at least two-dimensional, we derive that
$$s=1,\quad m_1=2,\quad \text{and}\quad \dim M=2\ ,$$ that is $M$ is the union of disjoint smooth orbits of the
$(\C^*)^2$-action.

The claim (3) follows from the fact that the stabilizer of any pair of points, one in
$\Tor^*(\sigma_1)$ and the other in $\Tor^*(\sigma_2)$, is finite.
\proofend

Recall that a deformation of an isolated
plane curve singularity is called {\it equiclassical} if it preserves both the total
$\delta$- and total $\kappa$-invariant (see \cite[Page 433, item (3)]{DH}).
The following statement is contained in \cite[Section 5]{DH}.

\medskip\noindent{\bf Lemma A5.} {\it
Let $k\ge2$, and let ${\mathcal B}_{2k}(S,z)\subset{\mathcal O}_{\C^2,z}$, ${\mathcal B}_{2k}\simeq(\C^{2k},0)$ be a versal deformation base of a plane curve singularity $(S,z)$ of type
$A_{2k}$. Denote by ${\mathcal B}_{2k}^{ec}(S,z)\subset{\mathcal B}_{2k}(S,z)$ the equiclassical stratum. Then ${\mathcal B}_{2k}^{ec}(S,z)$ is an irreducible germ of a complex variety of codimension $k-1$, its tangent cone $T{\mathcal B}_{2k}^{ec}(S,z)$ can be identified with the linear space
$$J^{ec}_{S,z}:=\{\varphi\in{\mathcal B}_{2k}(S,z)\ :\ \ord\varphi\big|_{S,z}\ge2k+1\}$$
of codimension $k-1$, a generic element of ${\mathcal B}_{2k}^{ec}(S,z)$
has $k-1$ nodes and one ordinary cusp.}
\medskip

We define the {\it multiplicity} $\mt({\mathcal B}_{2k}^{ec}(S,z))$ of the equiclassical stratum as follows. Let $\Lambda_t$, $t\in(\C,0)$, be a family of affine subspaces of ${\mathcal B}_{2k}(S,z)$ of dimension $k-1$ such that $\Lambda_0$ intersects
$T{\mathcal B}_{2k}^{ec}(S,z)$ transversally at the origin, and $\Lambda_t$, $t\ne0$, intersects
${\mathcal B}_{2k}^{ec}(S,z)$ in its generic elements. Set
$$\mt({\mathcal B}_{2k}^{ec}(S,z))=|\Lambda_t\cap{\mathcal B}_{2k}^{ec}(S,z)|\quad\text{for}\ t\ne0\ .$$
According to \cite[Proposition 5.2]{Shloc}, $\mt({\mathcal B}_{2k}^{ec}(S,z))=k$ for all $k\ge1$.

Define the multiplicity of the considered cuspidal tropical fragment $\mu_c(\Gamma_c,h_c)$ to be
zero if $|\Int({\mathcal T}\cap\Z^2)|\le1$, and, in case $|\Int({\mathcal T}\cap\Z^2)|\ge2$,
\begin{equation}\mu_c(\Gamma_c,h_c)=\sum_{C\in M^*}\mt({\mathcal B}_{2k}^{ec}(C,z))\ ,
\label{ecfe1}\end{equation}
where $M^*$ is the (finite) set of elements of $M$ given by polynomials having equal coefficients at
the vertices of ${\mathcal T}$, $(C,z)$ is a singular local branch of $C\in M^*$, and the type of $(C,z)$ is $A_{2k}$.

\medskip\noindent
{\bf Remark A6.} {\it
From Lemma \ref{ledge2}, one can extract an explicit formula for
$\mu_c(\Gamma_c,h_c)$ in the particular case of $(\Gamma_c,h_c)$ dual to the triangle
$\quad\quad\quad\quad$ ${\mathcal T}=\conv\{(0,0),(0,2),(m,1)\}$.}

\subsection*{A3. Patchworking: full version}
We present here an extension of Theorem \ref{t2} to arbitrary cuspidal plane tropical curves.

\medskip\noindent
{\bf Theorem A7.}
{\it Let $\Delta\subset\Z^2\setminus\{0\})$ be a nondegenerate, primitive, balanced multiset
satisfying the hypotheses of Lemma \ref{luc1}, $0\le g<p_a((P(\Delta))-1$,
$n=|\Delta|+g-2$, and $\bw$ a configuration of $n$ distinct points in $(\K^*)^2$ such that
$\bx=\val(\bn(\bw))$ is a set
of $n$ distinct point in $\R^2$ in general position. Let $(\Gamma,\bp,h)$ be a plane cuspidal $n$-marked
tropical curve of degree $\Delta$ and genus $g$ such that $h(\bp)=\bx$. Assume in addition that, in case B, the curve $(\Gamma,\bp,h)$ is not exceptional.
Then the number of curves $C\in V_{\Delta,g}(A_2)$ passing through $\bw$ and tropicalizing to $(\Gamma,
\bp,h)$ equals the expression (\ref{epw}),
where $\mu_c(\Gamma_c,\bp_c,h_c)$ should be appropriately chosen from formulas
(\ref{emuc1}), (\ref{emuc}), (\ref{ecfc1}), (\ref{ecfc2}), (\ref{ecfd2}), (\ref{ecfd1x}), (\ref{ecfd2x}), (\ref{eD}), or
(\ref{ecfe1}).}
\medskip

{\bf Proof.} This statement is a combination of Theorem \ref{t2}, Lemmas A2, A3 and the following extra argument related to the cuspidal tropical fragments of type E.

Let the cuspidal fragment be of type E, and let $M$ be the family from Lemma A4.
Choose a curve $C\in M$. For each singular point $z\in\Sing(C)$ and each intersection point
$z=C\cap D_i$ with the toric divisor $D_i=\Tor(\sigma_i)$, $i=1,2,3$, take the linear space
${\mathcal B}_z={\mathcal O}_{\Tor({\mathcal T}),z}/{\mathfrak m}_z^d$, where ${\mathfrak m}_z\subset{\mathcal O}_{\Tor({\mathcal T}),z}$ is the maximal ideal and $d\gg0$.
Furthermore, for each point $z\in\Sing(C)$ away from the singular branch of $C$ take the subspace
${\mathcal B}_z^{eg}\subset{\mathcal B}_z$ parameterizing local equigeneric deformations (i.e.,
preserving the $\delta$-invariant, see \cite[Page 433, item (2)]{DH}), for the point
$z_c\in\Sing(C)$, the center of the singular branch of $C$, take the subspace
${\mathcal B}_{z_c}^{ec}\subset{\mathcal B}_{z_c}$ parameterizing local equiclassical deformations,
for each point $z_i=C\cap D_i$, $i=1,2,3$, take the subspace ${\mathcal B}_{z_i}^{tan}\subset{\mathcal B}_{z_i}$
parameterizing deformations preserving the intersection multiplicity at the intersection point
with $D_i$ (while the intersection point may move along $D_i$). We have a natural embedding
of the germ of
$|{\mathcal L}_{\Tor({\mathcal T})}|$ at $C$ into $\prod_z{\mathcal B}_z$, where $z$ runs over
$\Sing(C)\cup\{z_1,z_2,z_3\}$. We claim that the image of that germ
intersects transversally in $\prod_z{\mathcal B}_z$ with the product
$$\prod_{z\in\Sing(C)\setminus\{z_c\}}{\mathcal B}_z^{eg}\times{\mathcal B}_{z_c}^{ec}\times\prod_{i=1}^3
{\mathcal B}_{z_i}^{tan}\ .$$ Indeed, this claim amounts to the $h^1$-vanishing (\ref{etrans}) established in the
proof of Lemma A4.
Note that the germ of the family $M$ at $C$ is the preimage of the considered intersection in
$|{\mathcal L}_{\mathcal T}|$. Combining this with the statement of Lemma A4(3), we obtain that
the choice of an admissible enhanced tropical limit (as in the proof of Theorem \ref{t2}),
containing a curve $C\in M$, together with the choice in the conditions to pass through the configuration $\bw$ yields $\mt({\mathcal B}_{2k}^{ec}(C,z))$ families $\{C^{(t)}\}_{t\in(\C,0)\setminus\{0\}}$ of
curves of genus $G$ with nodes and one cusp as required.
\proofend

\end{document}